\documentclass[11 pt]{amsart}
\usepackage{standalone}
\usepackage[utf8]{inputenc}
\usepackage{amssymb}
\usepackage{amsmath}
\usepackage{amsthm}
\usepackage{amsfonts}
\usepackage{comment}
\usepackage{young}
\usepackage{tikz, ifthen}
\usepackage[nomessages]{fp}
\usepackage{mathtools}
\usetikzlibrary{calc}
\usetikzlibrary{chains}
\usepackage{booktabs,caption}
\usepackage{dynkin-diagrams}
\usepackage{hyperref}
\usepackage{subcaption} 
\input xy
\xyoption{all}

\newcommand{\A}{\mathcal{A}}
\newcommand{\mO}{\mathcal{O}}
\newcommand{\des}{\text{des}}

\newtheorem{thm}{Theorem}[section]
\newtheorem{theorem}[thm]{Theorem}
\newtheorem{prop}[thm]{Proposition}
\newtheorem{proposition}[thm]{Proposition}
\newtheorem{lemma}[thm]{Lemma}

\newtheorem{cor}[thm]{Corollary}
\newtheorem{corollary}[thm]{Corollary}
\newtheorem{defn}[thm]{Definition}
\newtheorem{rem}[thm]{Remark}
\newtheorem{example}[thm]{Example}
\newtheorem{question}[thm]{Question}

\usepackage{tkz-graph}
\usepackage{tikz-cd}

\newcommand{\RR}{\mathbb{R}}
\newcommand{\CC}{\mathbb{C}}
\newcommand{\ZZ}{\mathbb{Z}}

\newcommand{\mfS}{\mathfrak{S}}

\DeclareMathOperator{\SL}{SL}

\DeclareMathOperator{\flip}{flip}
\DeclareMathOperator{\flatt}{flat}
\DeclareMathOperator{\Cat}{Cat}
\DeclareMathOperator{\Inc}{Br}
\newcommand{\iso}{\cong}
\newcommand{\mcD}{\mathcal{D}}
\newcommand{\mcB}{\mathcal{B}}
\newcommand{\redt}{\color{red}}

\newcommand{\newword}[1]{\textbf{\emph{#1}}}
\newcommand*{\defeq}{\stackrel{\text{def}}{=}}		

\newcommand{\vcenteredinclude}[1]{\begingroup
\setbox0=\hbox{#1}%
\parbox{\wd0}{\box0}\endgroup}

\tikzset{node distance=2em, ch/.style={circle,draw,on chain,inner sep=2pt},chj/.style={ch,join}, line width=1pt,baseline=-1ex}

\newcounter{x}
\newcounter{y}
\newcounter{z}

\newcommand*\cubecolors[1]{%
  \ifcase#1\relax
  \or\colorlet{cubecolor}{cyan}%
  \or\colorlet{cubecolor}{green}%
  \or\colorlet{cubecolor}{yellow}%
  \or\colorlet{cubecolor}{pink}%
  \or\colorlet{cubecolor}{purple}%
  \or\colorlet{cubecolor}{blue}%
  \or\colorlet{cubecolor}{white}%
    \else
    \colorlet{cubecolor}{black}%
  \fi
}
\newcommand\yaxis{180}
\newcommand\zaxis{-27}
\newcommand\xaxis{90}

\newcommand\topside[3]{
  \fill[fill=cubecolor, draw=black,shift={(\xaxis:#1)},shift={(\yaxis:#2)},
  shift={(\zaxis:#3)}] (0,0) -- (1,0) -- (0.5,0.25) --(-0.5,0.25)--(0,0);
}

\newcommand\leftside[3]{
  \fill[fill=cubecolor, draw=black,shift={(\xaxis:#1)},shift={(\yaxis:#2)},
  shift={(\zaxis:#3)}] (0,0) -- (0,-1) -- (-0.5,-0.75) --(-0.5,0.25)--(0,0);
}

\newcommand\rightside[3]{
  \fill[fill=cubecolor, draw=black,shift={(\xaxis:#1)},shift={(\yaxis:#2)},
  shift={(\zaxis:#3)}] (0,0) -- (1,0) -- (1,-1) --(0,-1)--(0,0);
}

\newcommand\topsidedashed[3]{
  \fill[fill=cubecolor, draw=black,shift={(\xaxis:#1)},shift={(\yaxis:#2)},
  shift={(\zaxis:#3)}] [dashed](0,0) -- (1,0) -- (0.5,0.25) --(-0.5,0.25)--(0,0);
}

\newcommand\leftsidedashed[3]{
  \fill[fill=cubecolor, draw=black,shift={(\xaxis:#1)},shift={(\yaxis:#2)},
  shift={(\zaxis:#3)}] [dashed](0,0) -- (0,-1) -- (-0.5,-0.75) --(-0.5,0.25)--(0,0);
}

\newcommand\rightsidedashed[3]{
  \fill[fill=cubecolor, draw=black,shift={(\xaxis:#1)},shift={(\yaxis:#2)},
  shift={(\zaxis:#3)}] [dashed](0,0) -- (1,0) -- (1,-1) --(0,-1) -- (0,0);
}

\newcommand\cube[3]{
  \topside{#1}{#2}{#3} \leftside{#1}{#2}{#3} \rightside{#1}{#2}{#3}
}

\newcommand\cubedashed[3]{
  \topsidedashed{#1}{#2}{#3} \leftsidedashed{#1}{#2}{#3} \rightsidedashed{#1}{#2}{#3}
}

\newcommand\ppAff[2]{
 \setcounter{x}{0}
 \foreach \a in {#2} {
    \addtocounter{x}{1}
    \setcounter{y}{-1}
    \foreach \b in \a {
      \addtocounter{y}{1}
      \ifthenelse{\b=0}{\addtocounter{y}{0}}{
        \FPeval{\result}{clip(#1-\the\value{y}-1)}
        \cubecolors{\b}
        \cube{\value{x}}{\value{y}}{0};
        \FPeval{\result}{clip(#1-\the\value{y}-1)}
        \draw[draw=black,shift={(\xaxis:\value{x})},shift={(\yaxis:\value{y})},
  shift={(\zaxis:0)}] (0.5,-0.5) node {\textsf{\result}};}
    }
  }
}

\newcommand\planepartitionD[2]{
 \setcounter{x}{0}
 \foreach \a in {#2} {
    \addtocounter{x}{1}
    \setcounter{y}{-1}
    \cubecolors{\value{x}}
    \foreach \b in \a {
      \addtocounter{y}{1}
      \setcounter{z}{-1}
      \foreach \c in {0,...,\b} {
        \addtocounter{z}{1}
      \ifthenelse{\c=0}{\setcounter{z}{-1},\addtocounter{y}{0}}{
        \FPeval{\newz}{clip(0.55*\the\value{z})}
        \cube{\value{x}}{\value{y}}{\newz};
        \FPeval{\result}{clip(#1-\the\value{y}-2*\the\value{z})}
        \draw[draw=black,shift={(\xaxis:\value{x})},shift={(\yaxis:\value{y})},
  shift={(\zaxis:\newz)}] (0.5,-0.5) node {\textsf{\result}};}
      }
    }
  }
}

\newcommand\ppAfff[2]{
 \setcounter{x}{0}
 \foreach \a in {#2} {
    \addtocounter{x}{1}
    \setcounter{y}{-1}
    \foreach \b in \a {
      \addtocounter{y}{1}
      \ifthenelse{\b=0}{\addtocounter{y}{0}}{
        \FPeval{\result}{clip(#1-\the\value{y}-1)}
        \cubecolors{\b}
        \ifthenelse{\result=0\OR\value{y}=0}
        {
        \cubedashed{\value{x}}{\value{y}}{0};
        \FPeval{\result}{clip(#1-\the\value{y}-1)}
        \draw[draw=black,shift={(\xaxis:\value{x})},shift={(\yaxis:\value{y})}, shift={(\zaxis:0)}] (0.5,-0.5) node {\textsf{0}};}
        {
        \cube{\value{x}}{\value{y}}{0};
        \FPeval{\result}{clip(#1-\the\value{y}-1)}
        \draw[draw=black,shift={(\xaxis:\value{x})},shift={(\yaxis:\value{y})},shift={(\zaxis:0)}] (0.5,-0.5) node {\textsf{\result}};}}
    }
  }
}
\allowdisplaybreaks[2] \textwidth15.1cm \textheight22cm \headheight12pt \oddsidemargin.4cm
\title{Coxeter groups and Billey-Postnikov decompositions}
\author{Suho Oh}
\email{s\_o79@txstate.edu}

\author{Edward Richmond}
\email{edward.richmond@okstate.edu}

\begin{document}

\begin{abstract}
In this chapter, we give an overview of Billey-Postnikov (BP) decompositions which have become an important tool for understanding the geometry and combinatorics of Schubert varieties.  BP decompositions are factorizations of Coxeter group elements with many nice properties in relation to Bruhat partial order.  They have played an important role in the classification and enumeration of smooth Schubert varieties.  They have also been used in the study of inversion hyperplane arrangements and permutation pattern avoidance.  We survey many of these applications.
\end{abstract}

\maketitle

\tableofcontents

\section{Introduction}

Coxeter groups and their combinatorics play a vital role in the study of Lie groups, flag varieties, and Schubert varieties.  Of particular importance are the length function and Bruhat partial order (also called Bruhat-Chevalley order) on a Coxeter group.  Geometrically, the length of an element equals the dimension of the corresponding Schubert variety, and the Bruhat order gives closure relations of Schubert cells in a flag variety.   Let $W$ be a Coxeter group, and for any $w\in W$, we define the Poincar\'e polynomial 
\[P_w(q):=\sum_{u\leq w} q^{\ell(w)}\]
where $\ell:W\rightarrow\ZZ_{\geq 0}$ denotes the length function and $\leq$ denotes Bruhat order.  If $W$ is the Weyl group of a complex reductive Lie group, then $P_w(q^2)$ is the topological Poincar\'e polynomial with respect to singular cohomology of the corresponding Schubert variety $X(w)$.  The central problem we address is to characterize when the polynomial $P_w(q)$ factors in some nice or natural way.  

When $W=\mfS_n$ is the permutation group, Gasharov proved in \cite{Ga98} that the polynomial $P_w(q)$ is a product of factors of the form $(1+q+\cdots +q^r)$ if and only if $w$ avoids the permutation patterns 3412 and 4231. 
 This is a generalization of the well known result that the Poincar\'e series on the full permutation group is
\begin{equation}\label{eqn:poincare_series_perms}
\sum_{w\in\mfS_n} q^{\ell(w)}=\prod_{k=0}^{n-1}(1+q+\cdots+q^k).
\end{equation}
Chevalley \cite{Ch55} and Solomon \cite{So66} prove that the Poincar\'e series of any finite Coxeter group has a similar factorization as in Equation \eqref{eqn:poincare_series_perms}, so we can ask if there is an analogue of Gasharov's factorization theorem for elements of other groups.  In \cite{LS90}, Lakshmibai and Sandya prove that a permutation $w\in\mfS_n$ avoiding the patterns 3412 and 4231 is equivalent to the Schubert variety $X(w)$ being smooth, and hence Gasharov's factorization theorem is a property of ``smooth" permutations.  As it turn out, nice Poincar\'e polynomial factorizations hold for not only smooth permutations, but ``smooth" elements of any Coxeter group of finite Lie-type.  We summarize these properties of (rationally) smooth elements in Section \ref{S:rat_smooth_intro}.  Of critical importance to these results is the notion of what is now called a Billey-Postnikov (BP) decomposition.  BP decompositions are certain factorizations of group elements with many nice properties in relation to intervals in Bruhat order.  The name ``BP" comes from the paper \cite{BP05} where Billey and Postnikov use these decompositions to give a root-theoretic criterion for the rational smoothness of Schubert varieties.  However, the idea of these decompositions dates back to earlier works such as \cite{Bi98} and \cite{Ga98}. BP decompositions have had numerous applications and have appeared in the study of fiber bundle structures of Schubert varieties \cite{RS16, RS17, RS18, Az23}, inversion hyperplane arrangements \cite{OPY08, OY10,SL15,MSH19} and permutation pattern avoidance \cite{AR18, GG20, Az23}.  The purpose of this chapter is to provide an overview of BP decompositions and their applications.  

We structure this chapter as follows.  In Section \ref{S:background}, we review the basic properities of Coxeter groups and define BP decompositions.  In Section \ref{S:rat_smooth_intro} we given an overview of the how BP decompositions are used to study rationally smooth elements of Coxeter groups.  In Section \ref{S:hyperplanes}, we discuss the applications of BP decompositions to inversion hyperplane arrangments.  In Section \ref{S:Geometry} we state and prove how BP decomposition correspond to fiber bundle structures on Schubert varieties.  In Section \ref{s:iterated_BPs}, we look at iterated BP decompositions and how they are modeled by staircase diagrams.  One application of staircase diagrams is that they can used to enumerate smooth and rationally smooth Schubert varieties.  In Section \ref{S:BP_pattern_avoidance}, we focus on permutation groups and discuss how BP decompositions are connected to permutation pattern avoidance.  Finally, in Section \ref{S:future} we state some open questions and possible future directions for the study of BP decompositions.  

Our hope is that this chapter will provide readers insights on the nature of BP decompositions and how they are applied.  Since this is a survey article, many statements will be given without proof.  If a proof is not provided, then we will either have brief outlines of the proof or have references provided.  Many results and concepts will have illustrating examples.

\section{Background on Coxeter groups}\label{S:background}
We review several foundational properties of Coxeter groups.  For more details, we refer readers to \cite{BB05}.  Let $W$ be a Coxeter group with simple generating set $S$.  In other words, $S$ is a finite set and $W$ is the group generated by $S$ where for any $s,t\in S$, we have a relation
\[(st)^{m_{st}}=e\]
for some $m_{st}\in \ZZ_{>0}\cup\{\infty\}$ where $m_{st}=1$ if and only if $s=t$.  We say an expression of $w\in W$ in the simple generators
\[w=s_1\cdots s_k\] is \newword{reduced} if $w$ cannot be expressed in fewer generators.  The length of any reduced expression is unqiue, so we define the function $\ell:W\rightarrow \ZZ_{\geq 0}$ which maps $w\in W$ to the length of any reduced expression.  We call the value $\ell(w)$ the \newword{length} of $w$.

Let $w=s_1\cdots s_k$ be a reduced expression.  We say that $u\leq w$ in the \newword{Bruhat order} if there exists a subsequence $(i_1,\ldots, i_j)\subseteq(1,\ldots, k)$ such that $u=s_{i_1}\cdots s_{i_j}$ is a reduced expression for $u$.  We remark that this definition is known as the sub-word property as Bruhat order has other equivalent definitions (see \cite[Theorem 2.2.2]{BB05}).

One particularly important family of Coxeter groups are the permutation groups on integers $\{1,\ldots,n\}$.  We will use the notation $W=\mfS_n$ when focusing on permutations.  As a Coxeter group, $\mfS_n$ has a simple generating set $S=\{s_1,\ldots,s_{n-1}\}$ where $s_i$ corresponds the simple transposition that swaps $i$ and $(i+1)$.  These generators satisfy the Coxeter relations 
\[s_i^2=(s_is_j)^2=(s_is_{i+1})^3=e\quad\text{for all $|i-j|>1$.} \]
The permutation group $\mfS_n$ is referred to as the Coxeter group of type $A_{n-1}$.
\begin{example}\label{Ex:type_A2}
The permutation group $\mfS_3$ is generated by $S=\{s_1, s_2\}$ and the Bruhat order is given by the following Hass\'e diagram:

$$\begin{tikzpicture}[scale=0.5]
  \node (max) at (0,4) {$s_1s_2s_1$};
  \node (a) at (-2,2) {$s_2s_1$};
  \node (b) at (2,2) {$s_1s_2$};
  \node (c) at (-2,0) {$s_1$};
  \node (d) at (2,0) {$s_2$};
  \node (min) at (0,-2) {$e$};
  \draw (d) -- (min) -- (c) -- (a) -- (max) -- (b) -- (d) -- (a);
  \draw[preaction={draw=white, -,line width=6pt}] (c) -- (b);
\end{tikzpicture}$$
\end{example}

\begin{figure}[htbp]
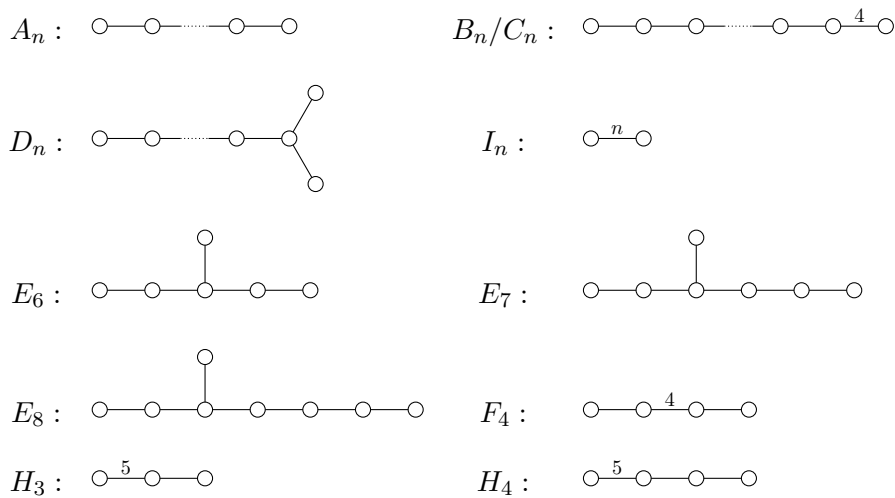

\begin{tabular}{clcl}
$A_n:$ & $\dynkin [scale=2, Coxeter] A{oo...oo}$ & $B_n/C_n:$ & $\dynkin [scale=2, Coxeter]B{ooo...ooo}$ \\  & \\
$D_n:$ & $\dynkin  [scale=2, Coxeter]D{oo...oooo}$ & $I_n:$ & $\dynkin  [scale=2, Coxeter,gonality=n]I{oo}$\\ & \\
$E_6:$ & $\dynkin  [Coxeter, mark=o, scale=2]E6$ & $E_7:$ & $\dynkin  [Coxeter, mark=o, scale=2]E7$\\  & \\
$E_8:$ & $\dynkin  [Coxeter, mark=o, scale=2]E8$& $F_4:$ & $\dynkin  [Coxeter, mark=o, scale=2]F4$\\  & \\
$H_3:$ & $\dynkin  [Coxeter, mark=o, scale=2]H3$& $H_4:$ & $\dynkin  [Coxeter, mark=o, scale=2]H4$\\ 
\end{tabular}
\caption{Coxeter diagrams of finite type.}\label{fig:Coxeter diagrams}
\end{figure}

Irreducible Coxeter groups of finite type are classified into four infinite families and six additional types.  This classification is commonly given in terms of Coxeter diagrams (or Dynkin-Coxeter diagrams).  The Coxeter diagram of a Coxeter group is a labeled graph with vertex set $S$ and edges $(s,t)$ labeled by the value $m_{st}$ under the conventions that we draw no edge if $m_{st}=2$ and an unlabelled edge if $m_{st}=3$.  See Figure \ref{fig:Coxeter diagrams} for the complete classification of irreducible finite Coxeter groups and note that Coxeter groups of types $H_2$ and $G_2$ are the dihedral groups $I_5$ and $I_6$ respectively. 

\smallskip

If $W$ is the Weyl group of finite dimensional Lie group, then we will refer to these Coxeter groups as ``finite Lie type" (these groups are also called finite crystalographic Coxeter groups).  In the classification found in Figure \ref{fig:Coxeter diagrams}, these are Coxeter groups of types $A_n, B_n, C_n, D_n, E_{6,7,8}, F_4$ and $G_2=I_6$.

In this section, unless otherwise specified, $W$ is an arbitrary Coxeter group with generating set $S$.  For any $u\leq w\in W$, we use the notation $$[u,w]:=\{z\in W\ |\ u\leq z\leq w\}$$
to denote intervals in Bruhat order.  For any $w\in W$, define the \newword{Poincar\'e polynomial} as the rank generating function on the lower order ideal $[e,w]$:
$$P_w(q):=\sum_{z\in [e,w]} q^{\ell(z)}.$$ 
If $X(w)$ denotes the Schubert variety corresponding to $w\in W$, then the Poincar\'e polynomial recovers the Hilbert-Poincar\'e series on singular cohomology:
$$P_w(q^2)=\sum_{k} \dim(H^k(X(w)))\, q^k.$$
For example, if we take $w=s_1s_2s_1$ as in Example \ref{Ex:type_A2}, then 
\[P_{w}(q) = 1+2q+2q^2+q^3.\]
In this case, $X(s_1s_2s_1)$ is the full flag variety of type $A_2$.  See Section \ref{S:Geometry} for a more detailed description of Schubert varieties. 

\begin{example}\label{example:free1231}
Let $W=\langle s_1,s_2,s_3\ |\ s_i^2=e\rangle$ be the the free Coxeter group on three generators and $w=s_1s_2s_3s_1.$
Then the Bruhat interval $[e,w]$ is:
$$\begin{tikzpicture}[scale=0.5]
  \node (max) at (0,4) {$s_1s_2s_3s_1$};
  \node (a1) at (-2,2) {$s_1s_2s_3$};
  \node (a2) at (-6,2) {$s_1s_2s_1$};
  \node (a3) at (6,2)  {$s_1s_3s_1$};
  \node (a4) at (2,2)  {$s_2s_3s_1$};
  \node (b1) at (-8,0) {$s_1s_2$};
  \node (b2) at (4,0) {$s_1s_3$};
  \node (b3) at (-4,0)  {$s_2s_1$};
  \node (b4) at (0,0)  {$s_2s_3$};
  \node (b5) at (8,0)  {$s_3s_1$};
  \node (c1) at (-0,-2) {$s_1$};
  \node (c2) at (-4,-2)  {$s_2$};
  \node (c3) at (4,-2)  {$s_3$};
  \node (min) at (0,-4) {$e$};
  \draw(a3)--(b2)--(a1)--(b4)--(a4)--(b3)--(c1)--(b2) --(c3) --(b5) -- (a3) -- (max) --(a2) --(b3) -- (c2) -- (min) -- (c1) -- (b1) -- (a1) -- (max) -- (a4) -- (b5) -- (c3) -- (min);
  \draw (c3)--(b4)--(c2)--(b1)--(a2);
  \draw (c1) -- (b5);
\end{tikzpicture}$$
and the Poincar\'{e} polynomial is
$$P_w(q)=1+3q+5q^2+4q^3+q^4.$$
\end{example}

Observe that the group inverse map $w \rightarrow w^{-1}$ is an automorphism of Bruhat order in the sense that 
$$u \leq w \leftrightarrow u^{-1} \leq w^{-1}.$$
One consequence is
 $$P_w(q) = P_{w^{-1}}(q)$$
 for any $w\in W$.
Taking $w = s_1s_2$ from Example \ref{Ex:type_A2}, we see
$$P_{s_1s_2}(q)=P_{s_2s_1}(q)= 1+2q+q^2.$$

\subsection{Quotients and parabolic decompositions}

In this section, we discuss parabolic quotients of Coxeter groups. 
 First, we say that a product $w=xy$ is a \newword{reduced factorization} if $\ell(w)=\ell(x)+\ell(y)$.  Let $J \subseteq S$ and let $W_J$ denote the subgroup of $W$ generated by the set $J$. Subgroups of form $W_J$ are called \newword{parabolic subgroups} of $W$. Each left coset $wW_J$ has a unique representative of minimal length. The set of minimal coset representatives can be defined as \[W^J:=\{w \in W \ |\ ws > w \text{ for all } s \in J\}.\]
The next proposition is from \cite[Proposition 2.4.4]{BB05}.
\begin{proposition}\label{prop:parabolic_decomposition}
Let $J \subseteq S$. Then the following hold:
\begin{enumerate}
    \item Every $w \in W$ has a unique factorization $w = vu$ such that $v \in W^J$ and $u \in W_J$.
    \item The decomposition $w=vu$ is a reduced factorization.  In other words, $$\ell(w)=\ell(v)+\ell(u).$$
\end{enumerate}
\end{proposition}

We call the decomposition $w=vu$ in Proposition \ref{prop:parabolic_decomposition} the \newword{parabolic decomposition with respect to $J$}.  We remark that each $w\in W$ also has a ``left-sided" parabolic decomposition $w=uv$ where $v$ denotes a minimal length representative of the right coset $W_Jw$.  If needed, we denote this set of minimal length representative by $^{J}W$.  However, the convention we take is that parabolic decompositions will be ``right-sided" decompositions $w=vu$ where $v$ is the minimal element of $wW_J$ and $u\in W_J$.

One consequence of Proposition \ref{prop:parabolic_decomposition} is that the coset decomposition of the group
\[W=\bigsqcup_{v\in W^J} vW_J\]
respects length in the sense that 
\[\sum_{w\in W} q^{\ell(w)}=\left(\sum_{v\in W^J} q^{\ell(v)}\right)\cdot \left(\sum_{u\in W_J} q^{\ell(u)}\right).\]
It is natural to ask if the analogous coset decomposition of the interval 
$$[e,w]=\bigsqcup_{v\in W^J}[e,w]\cap vW_J$$ gives a similar factorization of the Poincar\'e polynomial $P_w(q)$.  To make this question precise, we define relative Bruhat intervals and relative Poincar\'e polynomials.   For any $J\subset S$ and Bruhat interval $[u,v]$, define the \newword{Bruhat interval relative to $J$} as $$[u,v]^J:=[u,v]\cap W^J.$$  For any $v\in W^J$, we define the \newword{Poincar\'e polynomial relative to $J$} as
$$P^J_v(q):=\sum_{z\in [e,v]^J} q^{\ell(z)}.$$
\begin{defn}
Let $J\subseteq S$.  We say the parabolic decomposition $w=vu$ such that $v\in W^J$ and $u\in W_J$ is a \newword{Billey-Postinkov (BP) decomposition with respect to $J$} if the Poincar\'e polynomial factors as
$$P_w(q) = P_v^{J}(q)\cdot P_u(q).$$
\end{defn}
Proposition \ref{prop:parabolic_decomposition} implies that $w=vu$ is a BP decomposition  with respect to $J$ if and only if there is a graded poset isomorphism:
\[[e,v]^J\times [e,u]\simeq [e,w]\]
where $(v',u')\mapsto v'u'$.

\begin{example}\label{example:12321BP}
Let $w=s_1s_2s_3s_2s_1\in \mfS_3$ and Let $J=\{s_1,s_3\}$.  Then 
\[w=vu=(s_1s_3s_2)(s_1s_3)\]
is a BP decomposition with respect to $J$. Here we have 
\[[e,v]^J=\{e, s_2, s_1s_2, s_3s_2, s_1s_3s_2\}\quad\text{and}\quad [e,u]=\{e,s_1,s_3,s_1s_3\}\] with 
\[P_v^J(q)=1+q+2q^2+q^3\quad\text{and}\quad P_u(q)=1+2q+q^2.\]
The Poincar\'e polynomial 
\[P_w(q)=(1+q+2q^2+q^3)(1+2q+q^2)=1+3q+5q^2+6q^3+4q^4+q^5.\]
In Figure \ref{fig:12321_13}, we assign each coset a different color and see that the interval $$[e,v]^J\times [e,u]\simeq [e,w].$$
\begin{figure}[htbp]
    \centering
\begin{tikzpicture}[scale=0.5]

 \node[black!20!orange] (4231) at (0,8) {$s_1s_2s_3s_2s_1$};

 \node[black!20!blue] (4132) at (-7,6) {$s_3s_2s_3s_1$};
 \node[black!20!orange] (4213) at (2,6) {$s_3s_1s_2s_1$};
 \node[black!20!orange] (2431) at (-2,6) {$s_1s_2s_3s_2$};
 \node[black!20!yellow] (3241) at (7,6) {$s_1s_2s_3s_1$}; 

 \node[black!20!blue] (1432) at (-10,4) {$s_3s_2s_3$};
 \node[black!20!blue] (4123) at (-6,4) {$s_3s_2s_1$};
 \node[black!20!orange] (2413) at (-2,4) {$s_3s_1s_2$};
 \node[black!20!green] (3142) at (2,4) {$s_2s_3s_1$};
 \node[black!20!yellow] (3214) at (10,4) {$s_1s_2s_1$};
 \node[black!20!yellow] (2341) at (6,4) {$s_1s_2s_3$};

 \node[black!20!blue] (1423) at (-8,2) {$s_3s_2$};
 \node[black!20!green] (1342) at (-4,2) {$s_2s_3$};
  \node[black!20!red] (2143) at (0,2) {$s_3s_1$};
 \node[black!20!green] (3124) at (4,2) {$s_2s_1$};
 \node[black!20!yellow] (2314) at (8,2) {$s_1s_2$};

 \node[black!20!red] (1243) at (-5,0) {$s_3$};
 \node[black!20!green] (1324) at (0,0) {$s_2$};
 \node[black!20!red] (2134) at (5,0) {$s_1$};

  \node[black!20!red] (1234) at (0,-2) {$e$};
  
\draw[thick] (1234) -- (1243);
\draw[dashed] (1234) -- (1324);
\draw[thick] (1234) -- (2134);
\draw[dashed] (1243) -- (1342);
\draw[dashed] (1243) -- (1423);
\draw[thick] (1243) -- (2143);
\draw[thick] (1324) -- (1342);
\draw[dashed] (1324) -- (1423);
\draw[dashed] (1324) -- (2314);
\draw[thick] (1324) -- (3124);
\draw[thick] (2134) -- (2143);
\draw[dashed] (2134) -- (2314);
\draw[dashed] (2134) -- (3124);
\draw[dashed] (1342) -- (1432);
\draw[dashed] (1342) -- (2341);
\draw[thick] (1342) -- (3142);
\draw[thick] (1423) -- (1432);
\draw[dashed] (1423) -- (2413);
\draw[thick] (1423) -- (4123);
\draw[dashed] (2143) -- (2341);
\draw[dashed] (2143) -- (2413);
\draw[dashed] (2143) -- (3142);
\draw[dashed] (2143) -- (4123);
\draw[thick] (2314) -- (2341);
\draw[dashed] (2314) -- (2413);
\draw[thick] (2314) -- (3214);
\draw[thick] (3124) -- (3142);
\draw[dashed] (3124) -- (3214);
\draw[dashed] (3124) -- (4123);

\draw[dashed] (1432) -- (2431);
\draw[thick] (1432) -- (4132);

\draw[dashed] (2341) -- (2431);
\draw[thick] (2341) -- (3241);

\draw[thick] (2413) -- (2431);
\draw[thick] (2413) -- (4213);

\draw[dashed] (3142) -- (3241);
\draw[dashed] (3142) -- (4132);

\draw[thick] (3214) -- (3241);
\draw[dashed] (3214) -- (4213);

\draw[thick] (4123) -- (4132);
\draw[dashed] (4123) -- (4213);

\draw[thick] (2431) -- (4231);

\draw[dashed] (3241) -- (4231);

\draw[dashed] (4132) -- (4231);

\draw[thick] (4213) -- (4231);

\end{tikzpicture}
\caption{BP decomposoition $w=(s_1s_3s_2)(s_1s_3)$ with respect to $J=\{s_1,s_3\}$}
\label{fig:12321_13}
\end{figure}
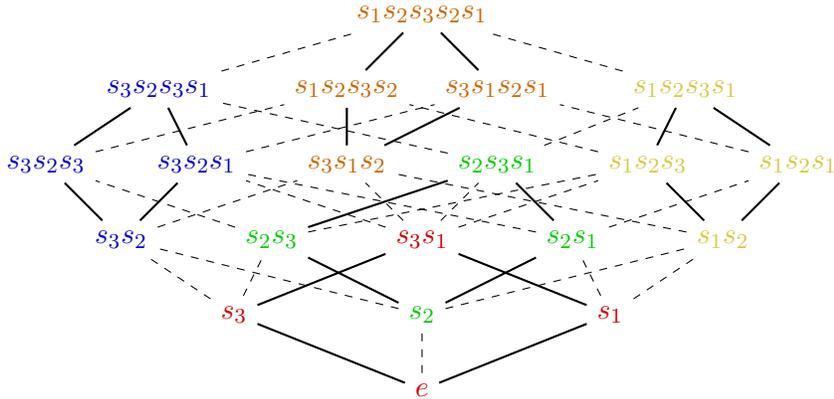
\end{example}

We remark that not all parabolic decompositions are BP decompositions.  In fact, BP decompositions are rather special and should not be expected in general. 

\begin{example}\label{example:12321not_BP}
Let $w=s_1s_2s_3s_2s_1\in \mfS_3$ and let $J=\{s_1,s_2\}$.  Then the parabolic decomposition
\[w=vu=(s_1s_2s_3)(s_2s_1)\]
is not a BP decomposition. Here we have 
\[[e,v]^J=\{e, s_3, s_2s_3, s_1s_2s_3\}\quad\text{and}\quad [e,u]=\{e,s_1,s_2,s_2s_1\}\] with 
\[P_v^J(q)=1+q+q^2+q^3\quad\text{and}\quad P_u(q)=1+2q+q^2.\]
The Poincar\'e polynomial 
\[P_w(q)\neq (1+q+q^2+q^3)(1+2q+q^2).\]
In Figure \ref{fig:12321_12}, we assign each coset in $[e,w]$ a different color and observe $[e,v]^J\times [e,u]$ and $[e,w]$ are not poset isomoprhic.
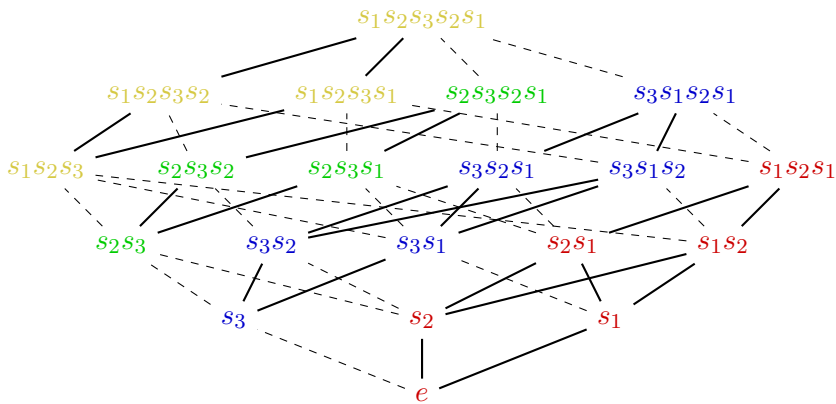
\begin{figure}[htbp]
    \centering
\begin{tikzpicture}[scale=0.5]

 \node[black!20!yellow] (4231) at (0,8) {$s_1s_2s_3s_2s_1$};

 \node[black!20!green] (4132) at (2,6) {$s_2s_3s_2s_1$};
 \node[black!20!blue] (4213) at (7,6) {$s_3s_1s_2s_1$};
 \node[black!20!yellow] (2431) at (-7,6) {$s_1s_2s_3s_2$};
 \node[black!20!yellow] (3241) at (-2,6) {$s_1s_2s_3s_1$}; 

 \node[black!20!green] (1432) at (-6,4) {$s_2s_3s_2$};
 \node[black!20!blue] (4123) at (2,4) {$s_3s_2s_1$};
 \node[black!20!blue] (2413) at (6,4) {$s_3s_1s_2$};
 \node[black!20!green] (3142) at (-2,4) {$s_2s_3s_1$};
 \node[black!20!red] (3214) at (10,4) {$s_1s_2s_1$};
 \node[black!20!yellow] (2341) at (-10,4) {$s_1s_2s_3$};

 \node[black!20!blue] (1423) at (-4,2) {$s_3s_2$};
 \node[black!20!green] (1342) at (-8,2) {$s_2s_3$};
  \node[black!20!blue] (2143) at (0,2) {$s_3s_1$};
 \node[black!20!red] (3124) at (4,2) {$s_2s_1$};
 \node[black!20!red] (2314) at (8,2) {$s_1s_2$};

 \node[black!20!blue] (1243) at (-5,0) {$s_3$};
 \node[black!20!red] (1324) at (0,0) {$s_2$};
 \node[black!20!red] (2134) at (5,0) {$s_1$};

  \node[black!20!red] (1234) at (0,-2) {$e$};
  
\draw[dashed] (1234) -- (1243);
\draw[thick] (1234) -- (1324);
\draw[thick] (1234) -- (2134);
\draw[dashed] (1243) -- (1342);
\draw[thick] (1243) -- (1423);
\draw[thick] (1243) -- (2143);
\draw[dashed] (1324) -- (1342);
\draw[dashed] (1324) -- (1423);
\draw[thick] (1324) -- (2314);
\draw[thick] (1324) -- (3124);
\draw[dashed] (2134) -- (2143);
\draw[thick] (2134) -- (2314);
\draw[thick] (2134) -- (3124);
\draw[thick] (1342) -- (1432);
\draw[dashed] (1342) -- (2341);
\draw[thick] (1342) -- (3142);
\draw[dashed] (1423) -- (1432);
\draw[thick] (1423) -- (2413);
\draw[thick] (1423) -- (4123);
\draw[dashed] (2143) -- (2341);
\draw[thick] (2143) -- (2413);
\draw[dashed] (2143) -- (3142);
\draw[thick] (2143) -- (4123);
\draw[dashed] (2314) -- (2341);
\draw[dashed] (2314) -- (2413);
\draw[thick] (2314) -- (3214);
\draw[dashed] (3124) -- (3142);
\draw[thick] (3124) -- (3214);
\draw[dashed] (3124) -- (4123);

\draw[dashed] (1432) -- (2431);
\draw[thick] (1432) -- (4132);

\draw[thick] (2341) -- (2431);
\draw[thick] (2341) -- (3241);

\draw[dashed] (2413) -- (2431);
\draw[thick] (2413) -- (4213);

\draw[dashed] (3142) -- (3241);
\draw[thick] (3142) -- (4132);

\draw[dashed] (3214) -- (3241);
\draw[dashed] (3214) -- (4213);

\draw[dashed] (4123) -- (4132);
\draw[thick] (4123) -- (4213);

\draw[thick] (2431) -- (4231);

\draw[thick] (3241) -- (4231);

\draw[dashed] (4132) -- (4231);

\draw[dashed] (4213) -- (4231);

\end{tikzpicture}
\caption{The parabolic decomposition $w=(s_1s_3s_2)(s_1s_3)$ is not a BP decomposition.}
\label{fig:12321_12}
\end{figure}
\end{example}

One important observation contrasting Examples \ref{example:12321BP} and \ref{example:12321not_BP}, is that the interval cosets $[e,w]\cap vW_J$ all have the same shape as $[e,u]$ in Example \ref{example:12321BP} (cosets are distinguished by different colors in Figure \ref{fig:12321_13}), while they do not in Example \ref{example:12321not_BP}.  In particular, when comparing the identity coset $[e,w]\cap W_J$ with the ``top" coset $[e,w]\cap vW_J$, if $w=vu$ is a BP decomposition, then $u$ must be the maximal element of $[e,w]\cap W_J$.  In \cite{BP05}, it is shown this maximally condition is also sufficient for the existence of a BP decomposition.  We remark that van den Hombergh proved in \cite{vdHo74} that the set $[e,w]\cap W_J$ always has a unique maximal element.  This fact was proved separately by Billey, Fan, and Losonczy in \cite{BFL99}.  

\subsection{Characterizing BP decompositions}

Our next goal is to list a several combinatorial characterizations of a BP decomposition and prove they are equivalent.  For any $w\in W$, define the \newword{support of $w$} as the set
$$S(w):=\{s\in S\ |\ s\leq w\}.$$
The support of $w$ can be viewed as the set of generators needed to make any reduced expression of $w$.  We also define the \newword{left and right descent sets} of $w$ as
$$D_L(w):=\{s\in S\ |\ \ell(sw)\leq \ell(w)\}$$
and  $$D_R(w):=\{s\in S\ |\ \ell(ws)\leq \ell(w)\}.$$
These decent sets can be thought of as the set of generators appearing on the left (respectively right) of some reduced expression of $w$.  Observe that $D_L(w)=D_R(w^{-1})$.  For example, if $w=s_1s_2s_1s_3\in \mfS_3$,
then \[S(w)=\{s_1,s_2,s_3\}, \quad D_L(w)=\{s_1,s_2\},\quad\text{and}\quad D_R(w)=\{s_1, s_3\}.\]
The following characterization theorem appears in \cite[Proposition 4.2]{RS16}.

\begin{theorem}\label{thm:BP_characterization}
Let $J\subset S$ and let $w=vu$ be a parabolic decomposition with respect to $J.$  Then the following are equivalent:
\begin{enumerate}
    \item $w=vu$ is a BP decomposition.
    \item The map $[e,v]^J\times[e,u]\rightarrow [e,w]$ given by $(v',u')\mapsto v'u'$ is bijective.
    \item $u$ is maximal in $[e,w]\cap W_J.$
    \item $S(v)\cap J \subseteq D_L(u).$
\end{enumerate}
\end{theorem}

\begin{proof}
We prove the theorem by showing the equivalencies: $(1)\leftrightarrow (2)$, $(2)\leftrightarrow (3)$, and $(3)\leftrightarrow (4)$. 

\smallskip

Proof of $(1)\leftrightarrow (2)$:  By Proposition \ref{prop:parabolic_decomposition}, the multiplication map given in part (2) is length preserving and injective.  Hence, part (2) says the interval $[e,w]$ decomposes as a product of posets $[e,v]^J\times[e,u]$.  This implies that part (1) is equivalent to part (2).  

\smallskip

Proof of $(2)\leftrightarrow (3)$:  We first assume the multiplication map is surjective (and hence bijective).  Then $[e,u]=[e,w]\cap W_J$ and hence $u$ must be the maximal element in $[e,w]\cap W_J$.  Conversely, assume $u$ is maximal in $[e,w]\cap W_J$ and let $x\in [e,w]$.  Let $x=v'u'$ denote the parabolic decomposition of $x$ with respect to $J$.  By \cite[Proposition 2.5.1]{BB05}, since $x\leq w$, we have $v'\leq v$.  But $u$ is maximal in $[e,w]\cap W_J$ and hence $u'\leq u$.  Thus the multiplication map is surjective.  

\smallskip

Proof of $(3)\leftrightarrow (4)$:  First suppose that $u$ is maximal in $[e,w]\cap W_J$ and $s\in S(v)\cap J$.  Then $su\in[e,w] \cap W_J$ and by the maximality of $u$, we must have $su\leq u$.  Hence $s\in D_L(u)$.  Conversely, suppose that $S(v)\cap J\subseteq D_L(u)$ and hence we can write a reduced factorization $u=u_0u'$ where $u_0$ is the maximal element in $W_{S(v)\cap J}.$  Let $x\in [e,w]\cap W_J$.  Since $x\leq w=vu=(vu_0)u'$, we can write a reduced factorization $x=u_1u_2$ where $u_1,u_2\in W_J$ with $u_1\leq vu_0$ and $u_2\leq u'$.  In particular, we have $u_1\in W_{S(v)\cap J}$ and hence $x=u_1u_2\leq u_0u'=u$.  Thus $u$ is maximal in $[e,w]\cap W_J$.
\end{proof}

While parts (3) and (4) of Theorem \ref{thm:BP_characterization} seem like less conventional ways to describe a BP decomposition, we will see in the subsequent sections these characterizations are very useful when working with BP decompositions.

\section{Rationally smooth elements of Coxeter groups} \label{S:rat_smooth_intro} 

In this section, we discuss BP-decompositions and ``rationally smooth" elements of Coxeter groups.  

\begin{defn}\label{defn:rationally_smooth}
We say $w\in W$ is \textbf{rationally smooth} if the coefficients of Poincar\'e polynomial 
\[P_w(q)=\sum_{i=0}^{\ell(w)} a_i q^i\]
satisfy $a_i=a_{\ell(w)-i}$ for all $0\leq i\leq \ell(w)$.  In other words, $P_w(q)$ is a palindromic polynomial.

Similarly, we say $v\in W^J$ is \textbf{rationally smooth with respect to $J$} if $P_v^J(q)$ is a palindromic polynomial.
\end{defn}
For example, $w=s_1s_2s_1\in\mfS_n$ as in Example \ref{Ex:type_A2} is rationally smooth, but the $w=s_1s_2s_3s_1$ in Example \ref{example:free1231} is not rationally smooth.  Note that rational smoothness and rational smoothness with respect to $J$ do not imply each other.  

\begin{example}
Let $w=s_2s_1s_3s_2\in\mfS_4$, then $w$ is rationally smooth with respect to $J=\{s_1,s_3\}$, but is not rationally smooth (See Figure \ref{fig:Bruhat2132}).
Here we have 
\[P_w(q)=1+3q+5q^2+4q^3+q^4\quad\text{and}\quad P_w^J(q)=1+q+2q^2+q^3+q^4.\]
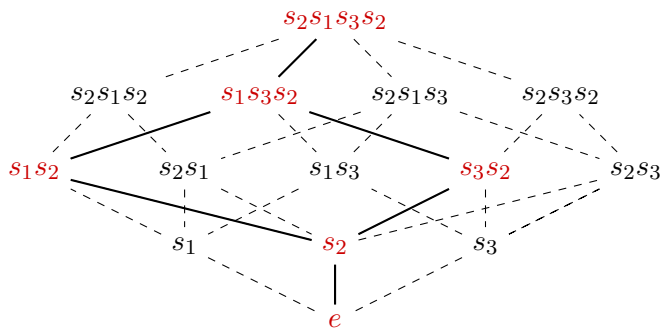
\begin{figure}[htbp]
\begin{tikzpicture}[scale=0.5]
  \node[black!20!red] (max) at (0,4) {$s_2s_1s_3s_2$};
  \node (a1)[black!20!red] at (-2,2) {$s_1s_3s_2$};
  \node (a2) at (-6,2) {$s_2s_1s_2$};
  \node (a3) at (6,2)  {$s_2s_3s_2$};
  \node (a4) at (2,2)  {$s_2s_1s_3$};
  \node[black!20!red] (b1) at (-8,0) {$s_1s_2$};
  \node[black!20!red] (b2) at (4,0) {$s_3s_2$};
  \node (b3) at (-4,0)  {$s_2s_1$};
  \node (b4) at (0,0)  {$s_1s_3$};
  \node (b5) at (8,0)  {$s_2s_3$};
  \node[black!20!red] (c1) at (-0,-2) {$s_2$};
  \node (c2) at (-4,-2)  {$s_1$};
  \node (c3) at (4,-2)  {$s_3$};
  \node[black!20!red] (min) at (0,-4) {$e$};
  \draw[dashed](a3)--(b2)--(a1)--(b4)--(a4)--(b3)--(c1)--(b2) --(c3) --(b5) -- (a3) -- (max) --(a2) --(b3) -- (c2) -- (min) -- (c1) -- (b1) -- (a1) -- (max) -- (a4) -- (b5) -- (c3) -- (min);
  \draw[dashed] (c3)--(b4)--(c2)--(b1)--(a2);
  \draw[dashed] (c1) -- (b5);
  \draw[thick] (min)--(c1)--(b1)--(a1)--(max);
  \draw[thick] (c1)--(b2)--(a1);  
\end{tikzpicture}
\caption{The Bruhat interval $[e,s_2s_1s_3s_2]$ with the subinterval $[e,s_2s_1s_3s_2]^J$, $J={\{s_1,s_3\}}$ highlighted in red.}\label{fig:Bruhat2132}
\end{figure}

We also see that if $w=s_1s_3s_2\in\mfS_4$, then $w$ is rationally smooth, but not rationally smooth with respect to $J=\{s_1,s_3\}$.  The polynomials
\[P_w(q)=1+3q+3q^2+q^3\quad\text{and}\quad P_w^J(q)=1+q+2q^2+q^3.\]
\end{example}

The term ``rationally smooth" is derived from the corresponding Schubert variety being rationally smooth in the geometric sense.  The geometric notion of rational smoothness was developed by Kazhdan and Lusztig \cite{KL79,KL80} where they show that a Schubert variety is rationally smooth if and only if certain Kazhdan-Lusztig polynomials are trivial.   It was proved by Carrell and Peterson in \cite{Ca94} that a Schubert variety $X(w)$ is rationally smooth if and only if $P_w(q)$ is a palindromic polynomial (this result also holds in the relative case with $X^J(w)$ and $P^J_w(q)$).  Note that Definition \ref{defn:rationally_smooth} well-defined for elements in any Coxeter group, even if there is no corresponding Schubert variety.  Any smooth variety is rationally smooth for topological reasons, however the converse is not true.  If $W$ is a simply-laced Coxeter of finite type (ie. type $A$, $D$ or $E$ in Figure \ref{fig:Coxeter diagrams}), then $X(w)$ is smooth if and only if it is rationally smooth.  This fact was proved by Deodhar in type $A$ \cite{De85} and then later in all simply-laced types by Carrell and Kuttler using ideas by Peterson in \cite{CK03}.

\subsection{BP decompositions of rationally smooth elements}
The next theorem connects BP decompositions to rationally smooth permutations and is a rephrasing of results due to Gasharov in \cite{Ga98} and, independently, due to Lascoux in \cite{La98}.
\begin{theorem}\label{thm:smooth_perm}
Let $w\in \mfS_n$.  Then $w$ is smooth if and only if either $w$ or $w^{-1}$ has a BP decomposition $vu$ with respect to $J=S\setminus\{s_{n-1}\}$ where
\[P_w(q)=(1+q+\cdots + q^{\ell(v)})\cdot P_u(q)\] and $u\in W_J\simeq \mfS_{n-1}$ is smooth.
\end{theorem}
In Theorem \ref{thm:smooth_perm}, the relative interval $[e,v]^J$ is a chain and hence relative Poincar\'e polynomial 
\[P^J_v(q)=1+q+\cdots + q^{\ell(v)}.\] 
Hence $v$ is rationally smooth with respect to $J$.  Polynomials of this form will come up frequently, so we use $q$-integer notation
\[[r]_q:=1+q+\cdots +q^{r-1}\]
for $r\in\ZZ_{>0}$.  Since $P_w(q)=P_{w^{-1}}(q)$, the reverse implication of Theorem \ref{thm:smooth_perm} follows from the fact products of $q$-integers are palindromic polynomials.  In Section \ref{S:BP_pattern_avoidance}, we provide a new proof of the forward direction of Theorem \ref{thm:smooth_perm} using ``split-pattern" avoidance which was developed by Alland and Richmond in \cite{AR18} to describe BP decompositions in $\mfS_n$ in one-line notation.

The following is a generalization of Theorem \ref{thm:smooth_perm} to Coxeter groups of finite Lie-type.

\begin{theorem}
\label{thm:BPgen}
Let $W$ be a Coxeter group of finite Lie-type and let $w\in W$ such that $|S(w)|\geq 2$.  Then $w$ is rationally smooth if and only if there is a leaf $s \in S(w)$ of the Coxeter diagram of $W_{S(w)}$ such that either $w$ or $w^{-1}$ has a BP decomposition $vu$ with respect to $J= S(w) \setminus\{ s\}$ where 
\begin{enumerate}
    \item $v$ is rationally smooth with respect to $J$ and 
    \item $u$ is rationally smooth.
\end{enumerate} 

Furthermore, $s\in S(w)$ can be chosen so that $v$ is either the maximal length element in $W_{S(v)}\cap W^J$, or one of the following holds:
\begin{enumerate}
    \item $W_{S(v)}$ is of type $B_n$ or $C_n$, with either
    \begin{enumerate}
        \item $J = S(w) \setminus \{s_1\}$, and $v = s_ks_{k+2}\dots s_ns_{n-1}\dots s_1$, for some $1 < k \leq n$.
        \item $J = S(w) \setminus \{s_n\}$ with $n \geq 2$ and $v = s_1 \dots s_n$
    \end{enumerate}
    \item $W_{S(v)}$ is of type $F_4$, with either
    \begin{enumerate}
        \item $J = S(w) \setminus \{s_1\}$ and $v = s_4s_3s_2s_1$
        \item $J = S(w) \setminus \{s_4\}$ and $v = s_1s_2s_3s_4$
    \end{enumerate}
    \item $W_{S(v)}$ is of type $G_2$, and $v$ is one of the elements
    $$ s_2s_1,\ s_1s_2s_1,\ s_2s_1s_2s_1,\ s_1s_2,\ s_2s_1s_2,\ s_1s_2s_1s_2.$$
\end{enumerate}
We use the standard conventions of Bourbaki \cite{Bo68} for the vertex labelling of Coxeter-Dynkin diagrams.
\end{theorem}

Note that the later part of Theorem \ref{thm:BPgen} only concerns Coxeter groups that are not simply laced.  For the classical types of $B,C$ and $D$, Theorem \ref{thm:BPgen} was proved by Billey in \cite{Bi98}.  The exceptional types were later proved in \cite{BP05} by Billey and Postnikov and in \cite{OY10} by Oh and Yoo.  We remark that if $w$ is rationally smooth of type $B/C$, then $P_w(q)$ factors into a product $q$-integers as in the type $A$ case \cite{Bi98}.  This is not necessarily true for Poincar\'e polynomials of rationally smooth elements of other types.  We also remark that, while the notion of rational smoothness is equivalent in types $B$ and $C$, the notion of smoothness is not since the set of smooth Schubert varieties is differs in these types.  We refer the reader to \cite{Bi98} and \cite{RS16,RS17} for the distinctions between smoothness in types $B$ and $C$.

We say a BP decomposition $w=vu$ with respect to $J$ is a \newword{Grassmannian BP decomposition} if $J$ is a maximal proper subset of $S(w)$.  Grassmannian BP decompositions are ``optimal" in the sense that they minimize the degree of the factor $P_v^J(q)$. Note that all the BP decompositions in Theorem \ref{thm:BPgen} are Grassmannian.  Moreover, since $J=S\setminus\{s\}$ where $s$ is leaf in the Coxeter diagram, the poset structure of the relative interval $[e,v]^J$ is less complex compared to when $s$ is not a leaf.  For example, in type $A$, the interval $[e,v]^J$ is always a chain of length $\ell(v)$ when $J=S\setminus\{s\}$ and $s$ is a leaf.  Grassmannian BP decompositions are discussed in more detail in Sections \ref{S:completeBP} and \ref{ss:split_patterns}.   

One issue with Theorem \ref{thm:BPgen} is the condition that ``either $w$ or $w^{-1}$ has a BP decomposition".  In Section \ref{S:Geometry}, we discuss how BP decompositions correspond to fiber bundle structures on Schubert varieties.  Since $X(w)$ is not isomorphic $X(w^{-1})$ in general, we would like an analogue of Theorem \ref{thm:BPgen} without the ``$w$ or $w^{-1}$" condition.  In \cite[Theorem 3.6]{RS16}, Richmond and Slofstra prove the following.

\begin{theorem}
\label{thm:BPexists}
Let $W$ be a Coxeter group of finite Lie-type.  If $w\in W$ is rationally smooth, then $w$ has a Grassmannian BP-decomposition with respect to $J=S(w)\setminus\{s\}$ for some $s\in S(w)$.
\end{theorem}

The sacrifice in Theorem \ref{thm:BPexists} is that we may not necessarily choose $s\in S(w)$ to be a leaf.  Theorem \ref{thm:BPexists} relies on Theorem \ref{thm:BPgen} and we give a brief outline of the proof in Section \ref{S:completeBP}.

\begin{example}
Let $w = s_2s_1s_3\in\mfS_4$. The support set $S(w)=\{s_1,s_2,s_3\}$ with Coxeter diagram
\begin{equation*}
  \begin{tikzpicture}[scale=.4]
    \draw[thick] (0,0)--(2,0)--(4,0);
    \draw[thick,fill=white] (0,0) circle (.3cm) node[label={[label distance=.1cm]-90:$s_1$}] { };
    \draw[thick,fill=white] (2,0) circle (.3cm) node[label={[label distance=.1cm]-90:$s_2$}]{ };
    \draw[thick,fill=white] (4,0) circle (.3cm) node[label={[label distance=.1cm]-90:$s_3$}]{ };
  \end{tikzpicture}
\end{equation*}
The Poincar\'e polynomial is \[P_{w}(q) = 1+3q+3q^2+q^3\] is palindromic and hence $w$ is rationally smooth. 
If $J=S(w)\setminus\{s_3\}=\{s_1,s_2\}$, then the parabolic decomposition with respect to $J$
\[w=vu=(s_2s_3)(s_1)\]
is not a BP decomposition. 
$$\begin{tikzpicture}[scale=0.5]
  \node[black!20!green] (a1) at (0,2) {$s_2s_1s_3$};
  \node[black!20!green] (b2) at (4,0) {$s_2s_3$};
  \node[black!20!red] (b3) at (-4,0)  {$s_2s_1$};
  \node[black!20!blue] (b4) at (0,0)  {$s_1s_3$};
  \node[black!20!red] (c1) at (-0,-2) {$s_2$};
  \node[black!20!red] (c2) at (-4,-2)  {$s_1$};
  \node[black!20!blue] (c3) at (4,-2)  {$s_3$};
  \node[black!20!red] (min) at (0,-4) {$e$};  
  \draw(min)--(c2)--(b3)--(a1)--(b4)--(c2);
  \draw (min)--(c1)--(b3);
  \draw (min)--(c3)--(b4)--(a1)--(b2)--(c3);
  \draw (c1)--(b2);
\end{tikzpicture}$$

Likewise $w=vu=(s_2s_3)(s_1)$ is not a BP decomposition with respect to $J=\{s_2,s_3\}$ and hence $w$ has no ``leaf removed" BP decomposition. 

If we take $J =S(w)\setminus\{s_2\}= \{s_1,s_3\}$, then \[w =vu= (s_2)(s_1s_3)\] is a Grassmannian BP decomposition.  Here we have 
\[P_{w}(q) = P^J_{s_2}(q) \cdot P_{s_1s_3}(q)=(1+q)\cdot (1+2q+q^2).\] 
$$\begin{tikzpicture}[scale=0.5]
  \node[black!20!blue] (a1) at (0,2) {$s_2s_1s_3$};
  \node[black!20!blue] (b2) at (4,0) {$s_2s_3$};
  \node[black!20!blue] (b3) at (-4,0)  {$s_2s_1$};
  \node[black!20!red] (b4) at (0,0)  {$s_1s_3$};
  \node[black!20!blue] (c1) at (-0,-2) {$s_2$};
  \node[black!20!red] (c2) at (-4,-2)  {$s_1$};
  \node[black!20!red] (c3) at (4,-2)  {$s_3$};
  \node[black!20!red] (min) at (0,-4) {$e$};  
  \draw(min)--(c2)--(b3)--(a1)--(b4)--(c2);
  \draw (min)--(c1)--(b3);
  \draw (min)--(c3)--(b4)--(a1)--(b2)--(c3);
  \draw (c1)--(b2);
\end{tikzpicture}$$
Observe that the inverse $w^{-1}=s_1s_3s_2$ does have a leaf removed BP decomposition with respect to both $J=\{s_1,s_2\}$ and $J=\{s_2,s_3\}$.
\end{example}







\subsection{Background on Permutations}
In this section, we focus on the permutation group $\mfS_n$.  For any $n\in \ZZ_{>0}$, let $[n]:=\{1,2,\ldots, n\}.$  Each permutation $w\in\mfS_n$ corresponds to a bijection $w:[n]\rightarrow [n]$ and has a unique presentation in one-line notation $w=w(1)w(2)\cdots w(n)$.  Under the Coxeter presentation of $\mfS_n$, the generators $s_i$ correspond to the simple transpositions swapping $i$ and $(i+1)$.  The right action of $s_i$ on the one-line notation of $w$ is given by swapping the $w(i)$ and $w(i+1)$ where the left action is given by swapping the position of the entries $i$ and $(i+1)$.  The length of a permutation can be calculated by counting inversions:
\[\ell(w)=\#\{(i,j)\ |\ i>j\ \text{and}\ w(i)<w(j)\}.\]
The Bruhat partial order is generated by the relations $w\leq w'$ where $w'$ is $w$ with two entries swapped and $\ell(w)<\ell(w')$.
For example, the Bruhat order on $\mfS_4$ is given in Figure \ref{fig:4321-bp}.
\begin{figure}[htbp]
    \centering
    \includegraphics{Figures/bruhat-4321-tikz.tex}
\caption{Bruhat order on $\mfS_4$.}
\label{fig:4321-bp}
\end{figure}

It will be common this chapter to state results for general Coxeter groups and then give more details in the case of permutations.  Sometimes it will be more convenient to use one-line notation over Coxeter theoretic reduced words to represent permutations.  In Section \ref{S:BP_pattern_avoidance}, we give a detailed overview of how BP-decompositions on permutations can be described using pattern avoidance.  Pattern avoidance has been a remarkable tool used to describe many properties of both permutations and Schubert varieties.  A survey of many of these results can be found in \cite{AB16}.

We conclude this section with a one-line notation version of Theorem \ref{thm:smooth_perm} that will be important in Section \ref{S:hyperplanes} on hyperplane arrangements.  Given $w \in \mfS_n$, we define the operation $\flatt(w,i)$ as the permutation in $\mfS_{n-1}$ obtained by removing the entry $w(i)$ from $w$ and then relabeling the remaining entries with $[n-1]$ while maintaining the relative order.  For example, $\flatt(635214,3)=54213$.

\begin{corollary}[\cite{Ga98},\cite{OPY08}]
\label{cor:BPonA}
Let $w \in \mfS_n$ be a smooth permutation and assume $w(d) = n$ and $w(n) = e$. Then at least one of the following two statements is true:
\begin{enumerate}
    \item $w(d) > w(d+1) > \dots > w(n)$, or
    \item $w^{-1}(e) > w^{-1}(e+1) > \dots > w^{-1}(n).$
\end{enumerate}
In both cases, the Poincar\'e polynomial factors as 
$$P_w(q) = [m+1]_q \cdot P_u(q),$$
where 
\begin{enumerate}
    \item $u = \flatt(w,d)$ and $m = n-d$ in the first case and 
    \item $u = \flatt(w,n)$ and $m = n-e$ in the second case.
\end{enumerate}
\end{corollary}

\begin{example}\label{example:2431}

Let $w=2431 = s_1s_2s_3s_2\in\mfS_4$. We have $d=2$ and $e=1$. Observe that $w(2) = 4 > w(3) = 3 > w(4) = 1$, so the first statement of Corollary \ref{cor:BPonA} holds. We have $u = \flatt(2431,2) = 231 = s_1s_2$ and $m = 4-2 = 2$. The Poincar\'e polynomial factors
\[P_{2431}(q) = (1+q+q^2)\cdot P_{231}(q)=(1+q+q^2)(1+2q+q^2).\] 
This factorization corresponds to the BP decomposition of \[w^{-1}=(s_2s_3)(s_2s_1)\] with respect to $J=\{s_1,s_2\}$.  In Figure~\ref{fig:2431-bp}, we highlight this decomposition in the Bruhat interval $[1234,2431]$.

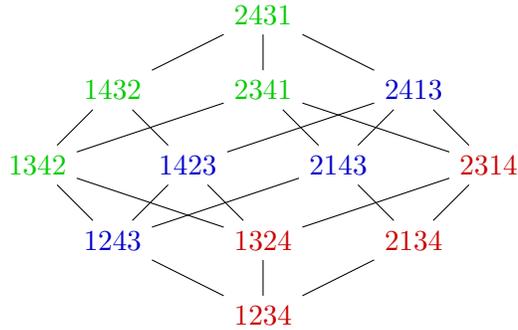
\begin{figure}[htbp]
\begin{tikzpicture}[scale=0.5]
  \node[black!20!green] (max) at (0,4) {$2431$};
  
  \node[black!20!green] (a1) at (-4,2) {$1432$};
  \node[black!20!green] (a2) at (0,2) {$2341$};
  \node[black!20!blue] (a3) at (4,2) {$2413$};
  
  \node[black!20!green] (b1) at (-6,0) {$1342$};
  \node[black!20!blue] (b2) at (-2,0) {$1423$};
  \node[black!20!blue] (b3) at (2,0)  {$2143$};
  \node[black!20!red] (b4) at (6,0)  {$2314$};
  
  \node[black!20!blue] (c1) at (-4,-2) {$1243$};
  \node[black!20!red] (c2) at (0,-2)  {$1324$};
  \node[black!20!red] (c3) at (4,-2)  {$2134$};
  
  \node[black!20!red] (min) at (0,-4) {$1234$};
  \draw(min)--(c1)--(b1)--(a1)--(max)--(a2)--(b1);
  \draw (min)--(c2)--(b2)--(a3)--(max);
  \draw (min)--(c3)--(b3)--(a3)--(b4)--(c3);  
  \draw (c2)--(b1); \draw (c2)--(b4); \draw (b3)--(c1)--(b2);
  \draw (a1)--(b2);\draw (b4)--(a2)--(b3);
\end{tikzpicture}
\caption{Bruhat order on $w=2431$.}
\label{fig:2431-bp}
\end{figure}
\end{example}









\section{Hyperplane arrangements}\label{S:hyperplanes}

In this section, we will give an overview of one of the major applications of the Billey-Postnikov decomposition, focusing on hyperplane arrangements. Let $W$ be a Coxeter group of finite Lie type.  For each $w \in W$, we will compare the Poincar\'e polynomial $P_w(q)$ with another polynomial, which comes from an associated hyperplane arrangement.


The poset structure we assign on the chambers of the hyperplane arrangements we study will be motivated from the weak Bruhat order on $W$.


\begin{defn}
Let $(W,S)$ be a Coxeter system and let $u,w \in W$.  The \textbf{right and left weak Bruhat orders} $\leq_R$ and $\leq_L$ are generated by the following cover relations.
\begin{enumerate}
    \item We have $u \leq_R w$ if $w = us$, for some $s\notin D_R(u)$. 
    \item We have $u \leq_L w$ if $w = su$, for some $s\notin D_L(u)$.
\end{enumerate}
\end{defn}

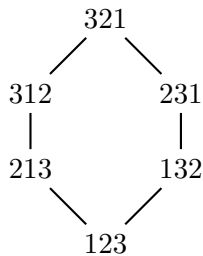
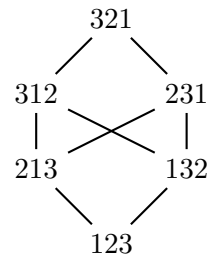
\begin{figure}[htbp]
    \centering
    \begin{subfigure}[b]{0.45\textwidth}
        \centering
        \centering
\begin{tikzpicture}[scale=0.5]

\node (321) at (0,6) {$321$};
\node (312) at (-2,4) {$312$};
\node (231) at (2,4) {$231$};
\node (132) at (2,2) {$132$};
\node (213) at (-2,2) {$213$};
\node (123) at (0,0) {$123$};

\draw[thick] (123) -- (213);
\draw[thick] (123) -- (132);
\draw[thick] (213) -- (312);
\draw[thick] (132) -- (231);
\draw[thick] (231) -- (321);
\draw[thick] (312) -- (321);

\end{tikzpicture}
        \caption{Left weak Bruhat order of $\mfS_3$}
        \label{fig:weak_bruhat}
    \end{subfigure}
    \hfill
    \begin{subfigure}[b]{0.45\textwidth}
        \centering
        \centering
\begin{tikzpicture}[scale=0.5]

\node (321) at (0,6) {$321$};
\node (312) at (-2,4) {$312$};
\node (231) at (2,4) {$231$};
\node (132) at (2,2) {$132$};
\node (213) at (-2,2) {$213$};
\node (123) at (0,0) {$123$};

\draw[thick] (123) -- (213);
\draw[thick] (123) -- (132);
\draw[thick] (213) -- (312);
\draw[thick] (132) -- (231);
\draw[thick] (231) -- (321);
\draw[thick] (312) -- (321);
\draw[thick] (231) -- (213);
\draw[thick] (132) -- (312);

\end{tikzpicture}
        \caption{Strong Bruhat order of $\mfS_3$}
        \label{fig:strong_bruhat}
    \end{subfigure}
    
    \caption{The left weak Bruhat order of $\mfS_3$ compared to the (strong) Bruhat order}
    \label{fig:weakvsstrong}
\end{figure}

An example of left weak Bruhat order of $\mfS_3$ is drawn in Figure~\ref{fig:weakvsstrong}. On the right side of the figure is the usual (strong) Bruhat order of $\mfS_3$. Notice that the set of elements is the same, and the rank of each element is the same between the two posets. This is true in general for any Coxeter group \cite{BB05}.

\subsection{Hyperplane arrangement of a Coxeter group}
Each finite Coxeter group $W$ is naturally associated with a hyperplane arrangement through a root system.  Let $\Phi$ be a finite collection of non-zero vectors in some Euclidean space $\RR^n$.  For each $\alpha\in\Phi$, we define the reflection $s_\alpha:\RR^n\rightarrow \RR^n$ by 
\[s_\alpha(x):=x-\frac{2(\alpha,x)}{(\alpha,\alpha)}\alpha.\]
We say $\Phi$ is a \newword{root system} of $W$ if the following hold:
\begin{enumerate}
    \item $\Phi\cap \RR\alpha=\{\alpha,-\alpha\}$ for each $\alpha\in\Phi$,
    \item $s_\alpha(\Phi)\subseteq \Phi$ for all $\alpha\in \Phi$, and 
    \item $W\simeq \langle s_\alpha\ |\ \alpha\in\Phi\rangle$.
\end{enumerate}
The vectors $\alpha\in\Phi$ are called \newword{roots}.  Let $S$ denote the set of simple generators of $W$.  For each $s\in S$, there is a unique root (up to sign) $\alpha_s$ which satisfies the condition that $s(\alpha_s)=-\alpha_s$.  We define a set of \newword{simple roots}
\[\Delta:=\{\alpha_s\ | \ s\in S\}\]
by selecting $\pm\alpha_s$ that all lie on one side of a suitably generic, but fixed hyperplane in $\RR^n$.  It can be shown that each $\alpha\in\Phi$ is either a totally positive or totally negative linear combination of simple roots, so we can decompose 
\[\Phi=\Phi^+\sqcup \Phi^-\]
into collections of positive and negative roots with respect to $\Delta$.  For more on root systems of Coxeter groups, see \cite[Chapter 1]{Hu90} or \cite[Chapter 4]{BB05}.  

To each $\alpha\in \Phi$, there is the corresponding hyperplane given by 
\[H_\alpha:=\{x\in \RR^n\ | \ (\alpha,x)=0\}.\]
Note that $H_\alpha=H_{-\alpha}$ and that if $x\in H_\alpha$, then $s_\alpha(x)=x$.  The collection of hyperplanes \[\A_W:=\{H_\alpha\ |\ \alpha\in \Phi^+\}\] is called the \newword{Coxeter arrangement of $W$}.  Since each $H_\alpha$ contains the origin in $\RR^n$, $\A_W$ is an example of a central hyperplane arrangement.


For each element $w \in W$, we can take a subset of hyperplanes in $\A_W$ corresponding to the inversions of $w$. 
 We define the \newword{inversion set} of $w$ as the set of positive roots
 $$\Phi_w := \{ \alpha\in \Phi^+\ |\ w(\alpha)\in \Phi^- \}.$$
and the \newword{inversion hyperplane arrangement}
$$\A_w:=\{H_\alpha\ |\ \alpha\in \Phi_w\}.$$
If $w_0$ denotes the longest element in $W$, then $\A_{w_0}=\A_W$.  For permutations, the root-theoretic inversion set $\Phi_w$ corresponds the usual inversion set of pairs:
\[\{(i,j)\in [n]^2\ |\ i<j\ \text{and}\ w(i) > w(j)\}.\]




Consider the permutation group $\mfS_n$ and let $\RR^n$ be a vectors space with coordinate basis $\{x_1,\ldots,x_n\}$.  The group $\mfS_n$ acts on $\RR^n$ by the standard permutation action the coordinate basis elements.  The set of vectors 
\[\Phi=\{x_i-x_j\ |\ i\neq j\}\]
is a root system of $\mfS_n$ with positive roots $\Phi^+=\{x_i-x_j\ |\ i>j\}$ and simple roots $\Delta=\{x_i-x_{i+1}\ |\ 1\leq i<n\}$.  If $\alpha=x_i-x_j$, then the hyperplane $H_\alpha$ is defined by the equation $x_i=x_j$.  The Coxeter arrangement of $\mfS_n$ is \[\A_{\mfS_n}=\{x_i=x_j\ | \ i\neq j\}.\]
For any $w\in \mfS_n$, we have the inversion arrangement
\[\A_{w}=\{x_i=x_j\ | \  i<j\ \text{and}\ w(i) > w(j) \}.\]
Note that the hyperplane $\sum_i x_i=0$ is invariant under the action of $\mfS_n$. Hence we can realize $\A_{\mfS_n}$ as a hyperplane arrangement in $$\RR^{n-1}\simeq \RR^n/(\sum_i x_i=0).$$  This reduction will help us visualize hyperplane arrangements in examples.

For any inversion arrangement $\A_w$, Let $r_0$ denote the \newword{fundamental chamber} which is defined as the set of points $x\in\RR^n$ such that  $(\alpha,x)>0$ for all $\alpha \in \Phi_w$. We define the \newword{distance enumerating polynomial} on $\A_w$ as
$$R_w(q) :=\sum_r q^{d(r_0,r)}$$ where the sum
is over all chambers of the arrangement $\A_w$ and $d(r_0,r)$ is the minimum number of hyperplanes separating
$r_0$ and $r$. 

\begin{example}
In Figure~\ref{fig:A2Rw}, we have the inversion arrangements $\A_{321}$ and $\A_{312}$. Since $321$ is the longest element of $\mfS_3$, we have $\A_{321} = \A_{\mfS_3}$ which consists of the three hyperplanes
\[\A_{\mfS_3}=\{x_1=x_2,\ x_1=x_2,\ x_2=x_3\}\]

Moreover, we can label each region with permutation in $\mfS_3$. The fundamental chamber is labelled with the identity permutation $123$. Starting from this, we measure the distance between each chamber and the fundamental chamber by counting the minimal number of hyperplanes needed to cross to reach the fundamental chamber.  We highlight this distance in blue. We can see that 
$$R_{321}(q) = 1+2q+2q^2 + q^3.$$ 

For $\A_{312}$, we remove the hyperplane $x_2 = x_3$ from $\A_{321}$. The fundamental chamber is the unique chamber that contains the identity label $123$.  Counting in a similar way, we obtain $$R_{312}(q) = 1+2q+q^2.$$
\begin{figure}[htbp]
    \centering
    \includegraphics[width=0.40\textwidth]{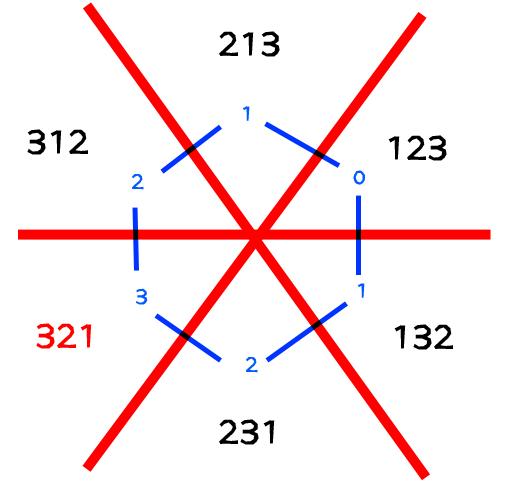}
    \includegraphics[width=0.40\textwidth]{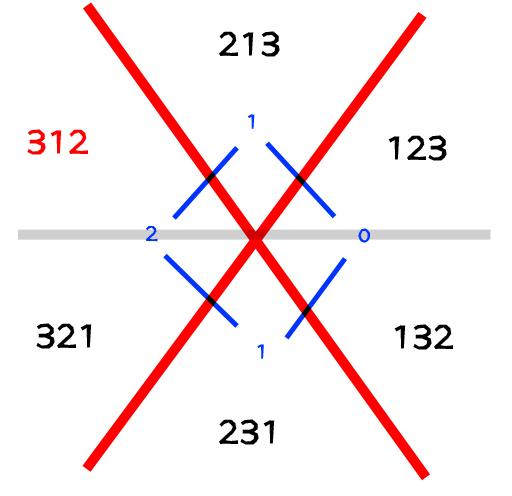}
\caption{Arrangements $\A_{321}$ and $\A_{312}$ and how to get the distance enumerating polynomial of $R_{321}$ and $R_{312}$.}
\label{fig:A2Rw}
\end{figure}
\end{example}


Our main goal of this section is to prove the following result found in \cite{OY10} and \cite{MSH19}:

\begin{theorem}\label{thm:hypmain}
Let $W$ be a Coxeter group of finite Lie-type. Then $w\in W$ is rationally smooth if and only if $P_w(q) = R_w(q)$.
\end{theorem}


\begin{example}\label{example:hyp:2431}
Let $w=4321$ denote the longest permutation in $\mfS_4$.  Then \[\A_{4321}=\{x_1=x_2,\ x_1=x_3,\ x_1=x_4,\ x_2=x_3,\ x_2=x_4,\ x_3=x_4 \}.\]
In Figure \ref{fig:A3fullarrangement}, we label the chambers by the values $d(r_0,r)$.  By the symmetry of the picture, we can see that $$R_{4321}=1+3q+5q^2+6q^3+5q^4+3q^5+q^6.$$

The rank generating function of $[1234,4321]$ in Figure~\ref{fig:4321-bp} is the same polynomial, so this verifies the fact that $P_{4321}(q) = R_{4321}(q)$.

If $w=2431$ (See Figure \ref{fig:A3fullarrangement}), then inversion arrangement
\[\A_{2431}=\{x_1 = x_4,\ x_2 = x_3,\ x_2 = x_4,\ x_3 = x_4\}\]
and $R_{2431}(q)=1+3q+4q^2+3q^3+q^4.$  From Figure~\ref{fig:2431-bp}, we see that $w$ is smooth and $P_{2431}(q) = R_{2431}(q)$.
\begin{figure}[htbp]
    \centering
    \includegraphics[width=0.35\textwidth]{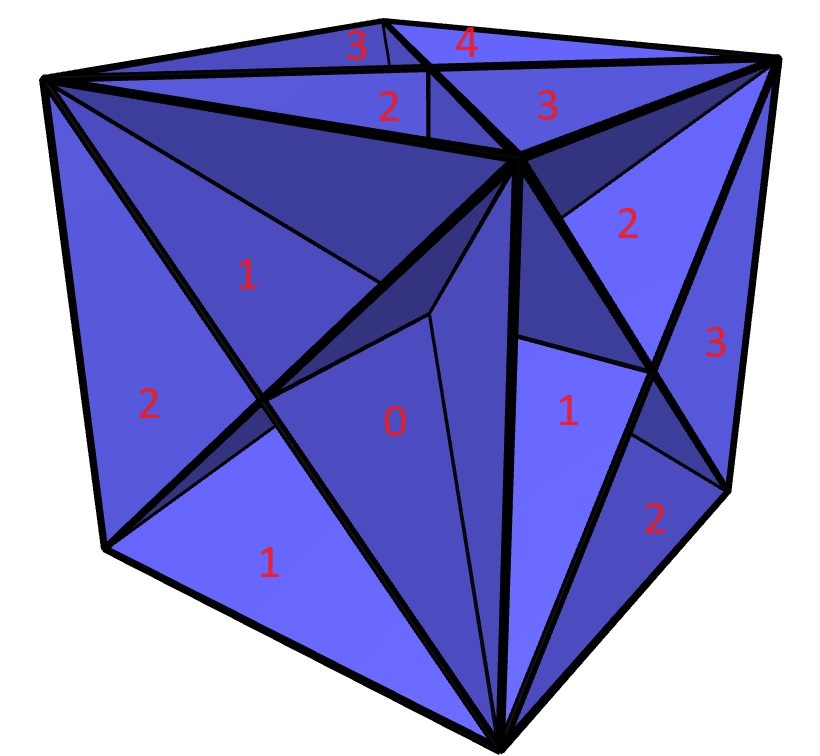}\hspace{0.4in}
     \includegraphics[width=0.35\textwidth]{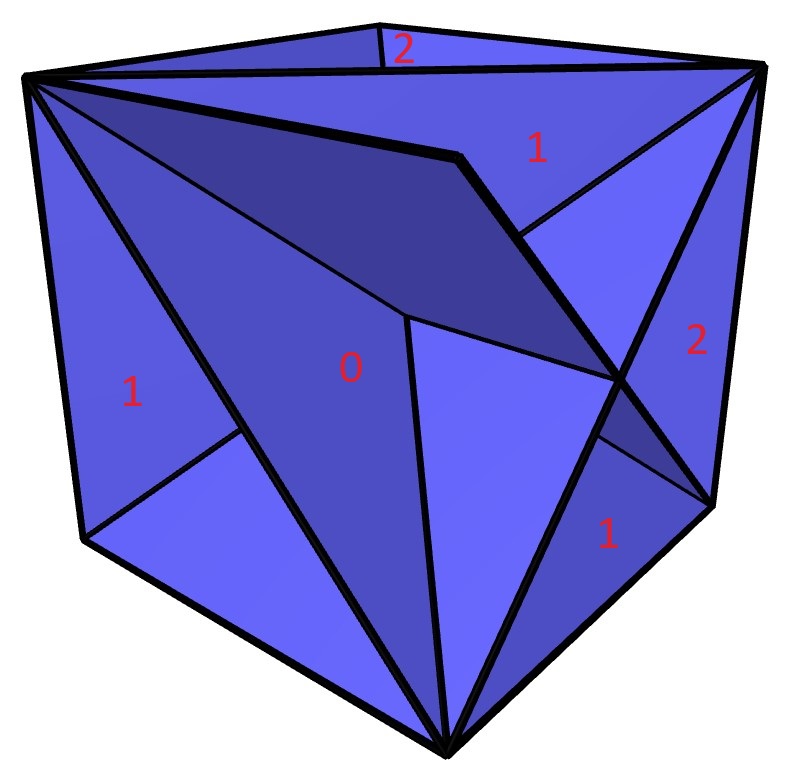}   
\caption{The inversion arrangements $\A_{4321}$ and $\A_{2431}$.}
\label{fig:A3fullarrangement}
\end{figure}

If $w=4231$, then $\A_{4231}=\A_{4321}\setminus\{x_2=x_3\}$.  
In Figure \ref{fig:A3arrangement4231}, we see 
\[R_{4231}(q)=1+4q+4q^2+4q^3+4q^5+q^5.\]
\begin{figure}[htbp]
    \centering
    \includegraphics[width=0.35\textwidth]{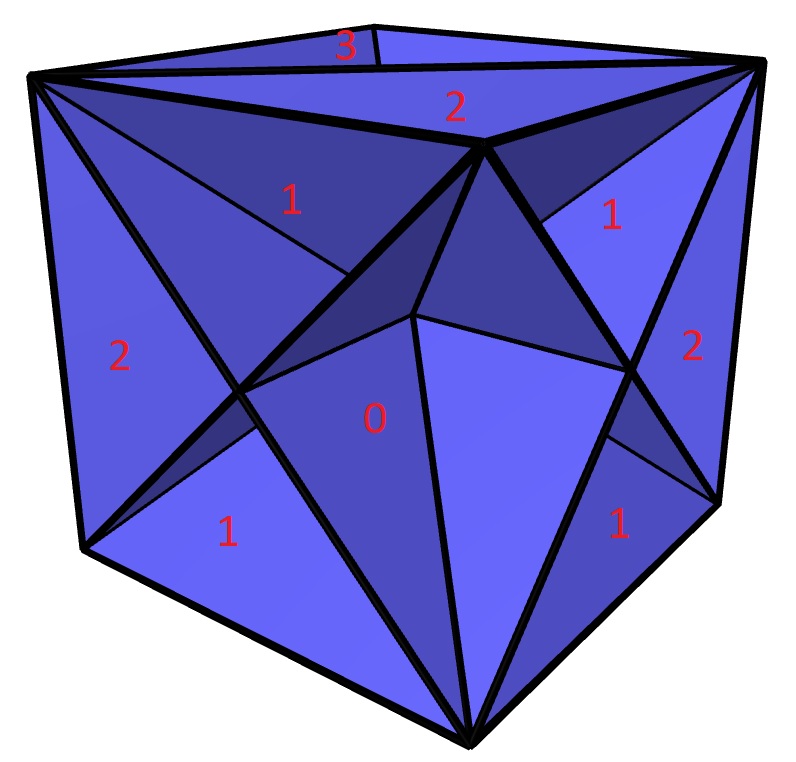}  
\caption{The inversion arrangement $\A_{4231}$.}
\label{fig:A3arrangement4231}
\end{figure}

In this case $w$ is not smooth (see Example \ref{example:12321BP}) and $$P_{4231}(q)= 1 + 3^q + 5q^2 + 6q^3 + 4q^4 + q^5.$$

\end{example}




\subsection{Inversion arrangements for permutations}
In this section, we provide a proof of Theorem \ref{thm:hypmain} for permutations. The main strategy is to show that $R_w(q)$ follows the exact same decomposition as $P_w(q)$ in Corollary \ref{cor:BPonA}.

For an undirected graph $G$ on vertex set $\{1,\ldots,n\}$, the \newword{graphical arrangement}, $\A_G$ is the hyperplane arrangement in $\RR^n$ with hyperplanes $x_i=x_j$ for all edges $(i,j)$ of $G$. In the case $G$ is a complete graph on $n$ vertices, we get the Coxeter arrangement $\A_{\mfS_n}$. 

Given a permutation $w \in \mfS_n$, we define its \newword{inversion graph} $G_w$ as an undirected graph on vertices $[n]=\{1,\ldots,n\}$ and edges $(i,j)$ whenever we have $i < j$ and $w(i) > w(j)$. Note that the inversion arrangement $\A_w$ is the graphical arrangement $\A_{G_w}$. 

An \newword{acyclic orientation} on $G$ is an assignment of directions to the edges of $G$ so that no directed cycles are formed. It is easy to see that the regions of $\A_G$ are in bijection with acyclic orientations of $G$. Indeed, if $\mO$ is an acyclic orientation of $G$, then we interpret each directed edge $i \rightarrow j$ as $x_i < x_j$. This corresponds to choosing a side of each hyperplane of $\A_G$, hence uniquely defining a region in $\A_G$.

From this observation, the distance enumerating polynomial $R_w(q)$ can be described in terms of acyclic orientations of the graph $G_w$. For an acyclic orientation $\mO$, let $\des(\mO)$ be the number of edges oriented oriented as $i \rightarrow j$ in $\mO$ where $i>j$ (so this corresponds to a descent of $w$). We define 
$$R_G(q) := \sum_{\mO} q^{\des(\mO)}.$$
It can be shown that $R_{G_w}(q)=R_w(q)$.  A \newword{clique} of $G$ is a subgraph of $G$ such that it is isomorphic to a complete graph. Given a graph $G$ and a vertex $k$ of $G$, let $G\setminus k$ denote the graph obtained by deleting $k$ and its adjacent edges in $G$.  The following lemma is from \cite{Bjorner1990}.

\begin{lemma}\label{lem:inversion_graph}
Suppose that a graph $G$ on vertex set $[n]$ has a vertex $k$ that satisfies the following two conditions:
\begin{enumerate}
    \item The neighbors of $k$ form a clique in $G$.
    \item Either all neighbors of $k$ are less than $k$ or all neighbors of $k$ are greater than $k$.
\end{enumerate}
Then $R_G(q) = [m+1]_q\cdot R_{G \setminus k}(q)$ where $m$ is the degree of the vertex $k$.
\end{lemma}

\begin{example}
Consider the inversion graph of $w=2431$
\[
\begin{tikzpicture}[scale=0.25]
\fill (0,0) circle (7pt) node[label={[label distance=.1cm]180:$4$}]{ };
\fill (0,4) circle (7pt) node[label={[label distance=.1cm]180:$1$}]{ };
\fill (4,4) circle (7pt) node[label={[label distance=.1cm]0:$2$}]{ };
\fill (4,0) circle (7pt) node[label={[label distance=.1cm]0:$3$}]{ };
\draw[thick] (0,0)--(4,4)--(4,0)--(0,0)--(0,4);
\end{tikzpicture}
\]
corresponding to inversions $(1,4),(2,3),(2,4),(3,4)$. The neighbor of $2$ is $\{3,4\}$, which is a clique, since $(3,4)$ is an edge.  Moreover all vertices in this clique is bigger than $2$.  By Lemma \ref{lem:inversion_graph}, we have $$R_{2431}(q) = (1+q+q^2)\cdot R_{231}(q).$$   
Comparing this to $P_{2431}(q)$ in Example \ref{example:2431}, we see these polynomials decompose in the exact same manner. 
\end{example}

\begin{proof}[Proof of Theorem \ref{thm:hypmain} for permutations]

Let $w\in \mfS_n$.  We need to check that the recursion in Lemma \ref{lem:inversion_graph} behaves exactly the same way as the recursion in Corollary~\ref{cor:BPonA}.  Using $d$ and $e$ as in Corollary~\ref{cor:BPonA}, we have $w(d) > w(d+1) > \dots > w(n)$ if and only if the neighbor of $n$ in $G_w$ forms a clique. The other case, $w^{-1}(e) > \dots > w^{-1}(n)$ happens if and only if the neighbor of $n$ in $G_{w^{-1}}$ forms a clique. From this observation that the two recurrences are the same whenever $w$ is a smooth permutation.
\end{proof}

\subsection{The general case}
We now consider the case when $W$ is a Coxeter group of finite Lie type.  
As with permutations, the main idea is to use BP-decompositions to show that the polynomials $P_w(q)$ and $R_w(q)$ follow the same recursion.  In type $A$, the analysis was simpler since we could always choose $v$ so it is the maximal length element of $W_{S(v)}\cap W^J$, in the context of Theorem~\ref{thm:BPgen}. 


For other types there are more cases of $v$ to consider, so the proof is more technical. Despite already having covered the type $A$ case, the examples we use in this section used to illustrate the ideas will also come from the type $A$ for simplicity.  We will be utilizing Theorem~\ref{thm:BPgen} as our main tool for decomposing the polynomials.  The results in this section are due to Mcalmon, Oh, and Yoo in \cite{MSH19} and Oh and Yoo in \cite{OY10}.




Let $\A$ be a central hyperplane arrangement with a fixed fundamental chamber $r_0$ and let $Q_{\A}$ denote the set of chambers of $\A$.  We define a poset structure on $Q_{\A}$ generated by the covering relations $r_1<r_2$ if chamber $r_1$ is adjacent to chamber $r_2$ and $d(r_1,r_0)=d(r_2,r_0)-1$.   

If $\A'$ is some subarrangement of $\A$ and $r\in Q_{\A'}$, we define the induced subposet $Q_{\A,\A',r}$ to be the subposet of $Q_{\A}$ obtained by restricting to the chambers of $\A$ contained in $r$.  We say that $\A$ is \newword{uniform} with respect to $\A'$ if for all chambers $r$ of $\A'$, the induced subposets $Q_{\A,\A',r}$ are all isomorphic.  In this case, we use $Q_{\A,\A'}$ to denote the poset.

\begin{example}
We consider the inversion arrangement of $w=4132$:
\[\A_{4132}=\{x_1 = x_2,\ x_1 = x_3,\ x_1 = x_4,\ x_3 = x_4\}\]

Now consider the hyperplane arrangement $\A_{3124}$ which is a subarrangement of $\A_{4132}$ by removing the hyperplanes $x_1 = x_4$ and $x_3 = x_4$.  In Figure \ref{fig:A3arrangement_subs}, we highlight $\A_{3124}$ in yellow.  
\begin{figure}[htbp]
    \centering
    \includegraphics[width=0.35\textwidth]{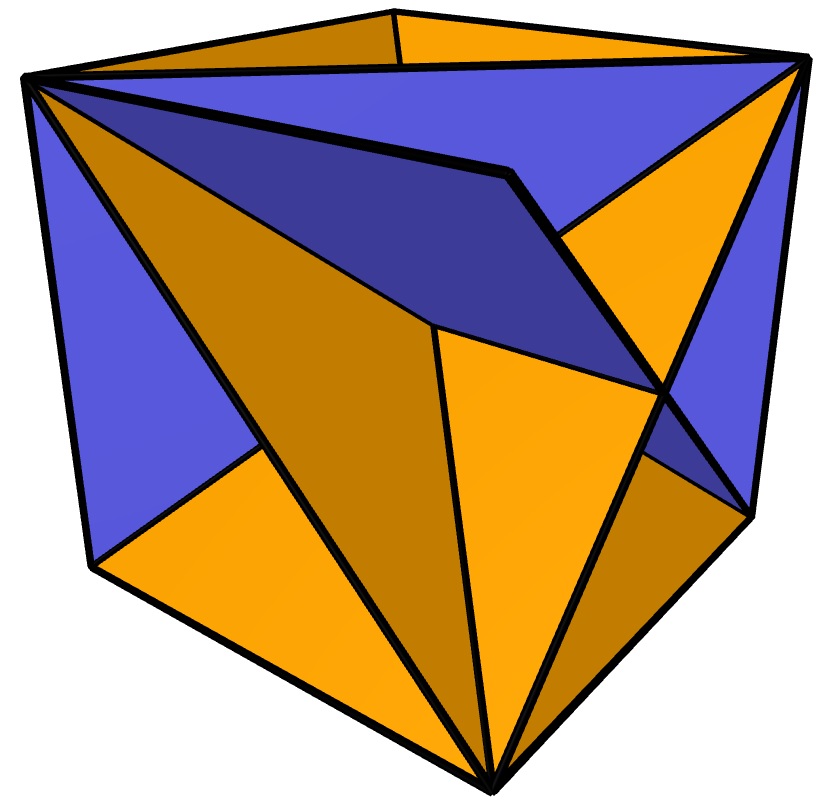}    
\caption{The inversion subarrangement $\A_{3124}$ (in yellow) of $\A_{4132}$.}
\label{fig:A3arrangement_subs}
\end{figure}
Let $r_0'$ denote fundamental chamber of $\A_{3124}$.  We see that $r_0'$ contains three chambers from $\A_{4132}$ and the poset $Q_{\A_{4132},\A_{3124},r_0'}$ is a chain of length $3$.  The same is true for all other chambers of $\A_{4132}$ and hence $\A_{4132}$ is uniform with respect to $\A_{3124}$.
\end{example}


Recall that if $w_0$ is the longest element of $W$, the arrangement $\A_{w_0}$ is the Coxeter arrangement of $W$. Here each chamber is indexed with a permutation $w \in W$ and two chambers $u,w$ are adjacent if and only if $w=su$ for some $s\in S$. Hence, the poset $Q_{\A_{w_0}}$ where $w_0$ is the longest element of $W$ is exactly the (left) weak Bruhat order of $W$. Recall that the weak Bruhat order of $W$ and the strong Bruhat order of $W$ are different poset structures on the same set of elements with the same rank \cite{BB05}. From this we get the following lemma.

\begin{lemma}
\label{lem:maxsame}
Let $w_0$ be the longest element of $W$.  Then $P_{w_0}(q) = R_{w_0}(q)$.  Furthermore, for any $J\subseteq S$, if $u_0$ is the longest element of $W_J$ for some $J \subset S$, then $\A_{w_0}$ is uniform with respect to $\A_{u_0}$. 
\end{lemma}
\begin{proof}
Each chamber of $\A_{w_0}$ is indexed by an element $w \in W$ and each $w \in W$ has a parabolic decomposition $vu$ where $u \in W_J$ and $v \in W^J$. The chambers indexed by $vu$ with common $u \in W_J$ are contained in the same chamber indexed by $u$ in $\A_{u_0}$. For each chamber $u$ in $\A_{u_J}$, the chambers of $\A_{w_0}$ contained in $u$ are only separated by hyperplanes in $\A_{w_0} \setminus \A_{u_0}$. The poset $Q_{\A_{w_0},\A_{u_0}}$ is the left weak Bruhat order of $W^J$.
\end{proof}

We use BP decompositions to develop some tools needed for the recursion on $R_w(q)$ when $w$ is rationally smooth.  Let $J\subseteq S$ and suppose we have a BP decomposition $w=vu$ with respect to $J$.  By Theorem \ref{thm:BP_characterization} part (4), we have that $S(v)\cap J\subseteq D_L(u)$.  In particular, we can write a reduced factorization 
$$u=u_{S(v)\cap J}\cdot u'$$
where $u_{S(v)\cap J}$ is the longest element of $W_{S(v)\cap J}$.  Theorem \ref{thm:BPgen} implies that if $w\in W$ is rationally smooth, then either $w$ or $w^{-1}$ has Grassmannian BP decomposition with respect to some $J$ of the form 
\[v\cdot (u_{S(v)\cap J}\cdot u')\]
For notational simplicity, let $I:=S(v)$.  Given such a decomposition, we decompose 
\[\A_w=\A_0\sqcup \A_1\sqcup \A_2\]
where 
\[\A_2:=\A_w\setminus \A_u, \quad \A_1:=\A_u\setminus\A_0,\quad\text{and}\quad \A_0:= (u')^{-1}\A_{u_{I\cap J}} \]

\begin{proposition}
\label{prop:chamred}
Let $r$ be some chamber inside $\A_1 \sqcup \A_0$. Let $r'$ be the chamber of $\A_0$ that contains $r$. Then the poset $Q_{\A_w, \A_1 \sqcup \A_0, r}$ is isomorphic to $Q_{\A_0 \sqcup \A_2, \A_0, r'}$.
\end{proposition}
\begin{proof}

Once a chamber $r'$ of $\A_0$ is fixed, we will show that any chamber of $\A_0 \sqcup \A_2$ contained in $r'$ intersects every chamber of $\A_1 \sqcup \A_0$ contained in $r'$. In order to show this, we can freely add more hyperplanes to $\A_0, \A_1$ and $\A_2$. So we may assume that $u = u_{I\cap J} u'$ is the longest element of $W_J$ and $v$ is the longest element of $W^J$. 

From Lemma~\ref{lem:maxsame}, each chamber of $\A_0$ is now indexed with a permutation of $W_{I \cap J}$. Fix a chamber $r_x$ labeled with a permutation $x \in W_{I \cap J}$. Each chamber of $\A_0 \sqcup \A_2$ contained in $r_x$ is labeled with a permutation $zx$ where $z \in W^J$. Each chamber of $\A_1 \sqcup \A_0$ contained in $r_x$ is labeled with a permutation $xy$ where $y^{-1} \in W^{I \cap J} \cap W_J$. For any such chamber of $\A_0 \sqcup \A_2$ and $\A_1 \sqcup \A_0$, their intersection will be the chamber of $\A$ that is labeled by $zxy\in W$. 

Let $r_1$ and $r_2$ be two different chambers of $\A$ contained in $r$. They are separated by a hyperplane in $\A_2$. For $i=1,2$, let $r_i'$ be the chamber of $\A_0 \sqcup \A_2$ that contains $r_i$. Then $r_1'$ and $r_2'$ are different chambers, since they are separated by the hyperplane that separates $r_1$ and $r_2$. If $r_1$ and $r_2$ are adjacent, then $r_1'$ and $r_2'$ are adjacent. If $r_1'$ and $r_2'$ are adjacent but $r_1$ and $r_2$ are not, it means there is a hyperplane of $\A_1$ that separates $r_1$ and $r_2$. But that contradicts the fact that $r_1$ and $r_2$ are both contained in the same chamber of $\A_1 \sqcup \A_0$. We conclude that $r_1$ and $r_2$ are adjacent if and only if $r_1'$ and $r_2'$ are adjacent. 
\end{proof}

\begin{example}

Let $w=4132$ and consider the arrangement $\A_{4132}$ from Figure~\ref{fig:A3arrangement_subs}. The BP-decomposition with respect to $J = \{s_1,s_2\}$ is $$w=(s_2s_3)(s_2s_1)$$ where $v = s_2s_3 = 1342 $ and $u = s_2s_1= 3124 $. The set $I=S(v)=\{s_2,s_3\}$ and $u_{I \cap J} = s_2 = 1324$. 

Looking at the inversions of $w$, we get:  $$\A_w=\{x_1 = x_2,\ x_1 = x_3,\ x_1 = x_4,\ x_3 = x_4\}$$ with 
$$\A_u = \{x_1 = x_2,\ x_1 = x_3\},\ \A_2= \{x_1 = x_4,\ x_3 = x_4\},\ \text{and}\ \A_{I \cap J}=\{x_2 = x_3\}.$$  The arrangement $\A_0 = (s_1)^{-1} \A_{I \cap J}=\{x_1 = x_3\}$ and hence $\A_1 = \{x_1=x_2\}$. From  Proposition \ref{prop:chamred}, we have that $Q(\A_{4132},\A_{3124})$ is isomorphic to $Q(\A_{3214},\A_{1324})$ where $s_2s_3s_2=3214$ and $s_2=1324$.  This poset is a chain of length $3$ by Lemma~\ref{lem:maxsame}.


\end{example}

From the above property we immediately get the following tool:
\begin{corollary}
\label{cor:hypmain}
Suppose we have a decomposition $w=v(u_{I\cap J}u')$ as in Proposition \ref{prop:chamred} and assume $\A_{vu_{I \cap J}}$ is uniform with respect to $\A_{u_{I \cap J}}$.  

If $R_{vu_{I \cap J}}(q) = P_{vu_{I \cap J}}(q)$ and $R_u(q) = P_u(q)$, then  $R_w(q) = P_w(q)$.

\end{corollary}
\begin{proof}

If $\A_{vu_{I \cap J}}$ is uniform with respect to $\A_{u_{I \cap J}}$, then Proposition~\ref{prop:chamred} tells us that $\A_w$ is uniform with respect to $\A_u$. Hence $R_w(q)$ is divisible by $R_u(q)$. Moreover, $$\displaystyle\frac{R_w(q)}{R_u(q)} = \frac{R_{vu_{I \cap J}}(q)}{R_{u_{I \cap J}}(q)}.$$ From Lemma~\ref{lem:maxsame}, we have $R_{u_{I \cap J}}(q) = P_{u_{I \cap J}}(q)$. Hence $R_{vu_{I \cap J}}(q) = P_{vu_{I \cap J}}(q)$ and $R_u(q) = P_u(q)$ implies $R_w(q) = P_w(q)$.
\end{proof}

Corollary  \ref{cor:hypmain} allows us to consider only the case where $u$ is the longest element of some $W_I$. 

Let $I$ be the set of simple roots that appear in a reduced word of $v$. We say that $v$ is a \newword{locally-maximal} element in $W^J$ if it is the maximal element of $W_I^{I \cap J}:=W_I\cap W^{I \cap J}$ and $I$ forms a connected subgraph within the Coxeter diagram. Similarly we say that $v$ is in a \newword{local chain} if $W_I^{I \cap J}$ is a chain poset. Notice that in Theorem~\ref{thm:BPgen}, only case when $v$ is not locally-maximal nor a local chain lies in Coxeter groups of type $F_4$ and $B_n$.


\begin{proposition}
\label{prop:hypmain}
Suppose we have a decomposition $w=v(u_{I\cap J}u')$ as in Proposition \ref{prop:chamred}.  

If $v$ is a locally-maximal element or a local chain, then $P_u(q)=R_u(q)$ implies $P_w(q) = R_w(q)$.
\end{proposition}
\begin{proof}

From corollary \ref{cor:hypmain}, it is enough to show $\A_{vu_{I \cap J}}$ is uniform with respect to $\A_{u_{I \cap J}}$ and $R_{vu_{I \cap J}}(q) = P_{vu_{I \cap J}}(q)$.

If $v$ is the longest element of $W^J$, then $vu_{I \cap J}$ is the longest element of $W_I$. In this case, the proposition follows from Lemma~\ref{lem:maxsame}.


When $W_{I}^{I \cap J}$ is a chain, let $v'$ denote the longest element of $W_{I}^{I \cap J}$. Then $w' := v'u_{I \cap J}$ is the longest element of $W_I$. From Lemma~\ref{lem:maxsame}, we have that $R_{u_{I \cap J}}(q) = P_{u_{I \cap J}}(q)$ and $R_{v'u_{I \cap J}}(q) = P_{v'u_{I \cap J}}(q)$. For each chamber $r$ of $\A_{u_{I \cap J}}$, the poset $Q(\A_{w'},\A_{u_{I \cap J}},u)$ is a chain of length $\ell(v')$. In particular, every hyperplane of $\A_{w'} \setminus \A_{u_{I \cap J}}$ intersects the interior of the chamber $r$.

When we go from $\A_{w'=v'u_{I \cap J}}$ to $\A_{vu_{I \cap J}}$, we remove some hyperplanes from $\A_{w'} \setminus \A_{u_{I \cap J}}$. For each chamber $r$ of $\A_{u_{I \cap J}}$, the poset $Q(\A_{v'u_{I \cap J}},\A_{u_{I \cap J}},r)$ is a chain of length $\ell(v')$ minus the number of hyperplanes removed. Hence $\A_{vu_{I \cap J}}$ is uniform with respect to $\A_{u_{I \cap J}}$. Moreover, we have $$R_{vu_{I \cap J}}(q) = (1+\cdots+q^{\ell(v)})\cdot R_{u_{I \cap J}(q)}.$$ The proposition now follows from Lemma~\ref{lem:maxsame}.
\end{proof}

Lastly we analyze two special examples each coming from Coxeter groups of type $F_4$ and $B_n$ which will be needed for our main result. These examples correspond to parts (1b) and (2b) in Theorem \ref{thm:BPgen}.  We start with type $F_4$:

\begin{example}
\label{ex:F4special}
Let $W$ be a Coxeter group of type $F_4$. Let $w=vu$ where $u$ is the longest element of $W_{\{s_1,s_2,s_3\}}$ and $v = s_1s_2s_3s_4$. Then $w=vu$ is a BP decomposition and $$P_w(q) =(1+q+q^2+q^3)\cdot P_{u}(q).$$  
The root system of type $F_4$ lies in $\RR^4$ and the hyperplane arrangement $\A_w$ is the union of
\begin{align*}
\A_u=\{x_1=0,\ &x_2=0,\ x_3=0,\ x_2-x_1=0,\ x_3-x_2=0,\\ 
&x_3-x_1=0,\ x_1+x_2=0,\ x_1+x_3=0,\ x_2+x_3=0\}
\end{align*}
and the hyperplanes 
\begin{align*}
\{x_1+x_2+x_3=x_4,\ -x_1-x_2+x_3=x_4,\ 
-x_1+x_2-x_3=x_4,\ x_1-x_2-x_3=x_4\}.
\end{align*}
Pick any chamber $c$ of $\A_u$ and an arbitrary interior point $z = (z_1,z_2,z_3,z_4)\in c$. Consider the line $l_z$ obtained from $z$ by changing the $z_4$ value from $-\infty$ to $+\infty$. This line is still contained in chamber $c$. Imagine moving through the line $l_z$ by changing the $z_4$ value from $-\infty$ to $+\infty$. The difference between any pair of equations of hyperplanes in $\A_w \setminus \A_u$ is of form $2x_i + 2x_j = 0$. For each pair $i \neq j \leq 3$ whether $x_i+x_j$ is positive or negative is determined by the choice of $c$ since $x_i+x_j=0$ is a hyperplane of $\A_u$. Therefore, the order we cross the hyperplanes of $\A_w \setminus \A_{u}$ is completely determined by $c$.

From this we can conclude that $\A_w$ is uniform with respect to $\A_u$. Moreover, the poset $Q_w$ is obtained from $Q_u$ by a poset product with a chain of length $4$. We get 
$$R_w(q) = (1+q+q^2+q^3 + q^4)\cdot R_{u}(q).$$ 
Since $R_{u}(q) = P_{u}(q)$ from Lemma~\ref{lem:maxsame} and $P_v^{W^J}(q) = (1+ \cdots + q^{\ell(v)})$, we obtain the desired result.
\end{example}

Now we consider the case of $B_n$ with the leaf $s_n$ in its Coxeter diagram.

\begin{lemma}
\label{lem:Bnspecial}
Let $W$ be a type $B_n$ Coxeter group and simple generating set $S=\{s_1,\ldots, s_n\}$ and let $J = S \setminus \{s_n\}$. Let $w = v u$ where $u$ is the longest element of $W_J$ and $v =  s_1 \cdots s_{n-1} s_n$.  Then $\A_w$ is uniform with respect to $\A_u$ and $P_w(q) = R_w(q)$.
\end{lemma}
\begin{proof}
The root system of type $B_n$ lies in $\RR^n$ and the hyperplane arrangement $\A_u$ consists of hyperplanes defined by the following equations:
\begin{enumerate}
    \item $x_i = 0$ for $1 \leq i \leq n-1$,
    \item $x_i - x_j = 0$ for $1 \leq i < j \leq n-1$, and
    \item $x_i + x_j = 0$ for $1 \leq i < j \leq n-1$.
\end{enumerate}
The hyperplane arrangement $\A_w$ is obtained from $\A_u$ with the additional the hyperplanes $x_n=0$ and $x_n + x_i = 0$ for $1 \leq i \leq n-1$. 

Pick any chamber $c$ of $\A_u$ and arbitrary interior point $z = (z_1,\ldots,z_n)\in c$. The chamber $c$ determines a total order on $z_1,\ldots,z_{n-1}$ and $0$ that does not depend on the choice of $z$. Consider the line $l_z$ obtained from $z$ by changing the value of $z_n$ from $-\infty$ to $+\infty$ and that $l_z$ is contained in the chamber $c$. As the value of $z_n$ moves from $-\infty$ to $+\infty$ along $l_z$, the order in which we cross the hyperplanes of $\A_w \setminus \A_u$ is determined by the total order on $z_1,\ldots,z_{n-1}$ and $0$.

Hence $\A_w$ is uniform with respect to $\A_u$. Moreover, the poset $Q_w$ is obtained from $Q_u$ by a poset product with a chain of length $n$. We get $R_w(q) = (1+ \cdots + q^{\ell(v)})\cdot R_u(q) $. Since $R_u(q) = P_u(q)$ from Lemma~\ref{lem:maxsame} and $P_v^{W^J}(q) = (1+ \cdots + q^{\ell(v)})$, we obtain the desired result.
\end{proof}

We are now ready to prove Theorem \ref{thm:hypmain}.  First note that $R_w(q)$ is always palindromic by definition.  So if $w\in W$ is not rationally smooth, $P_w(q)\neq R_w(q)$.  From the above Proposition~\ref{prop:hypmain} and Theorem~\ref{thm:BPgen} we can obtain the following result which completes the proof.

\begin{theorem}\label{thm:hypmain2}
Let $W$ be a Coxeter group of finite Lie-type and let $w\in W$ be rationally smooth. Then $R_w(q) = P_w(q)$. 
\end{theorem}
\begin{proof} 
We use induction on $|S(w)|$. First, if $w=s\in S$, then $R_w(q)=P_w(q)=1+q$. By Theorem \ref{thm:BPgen}, either $w$ or $w^{-1}$ has a Grassmannian BP decomposition $vu$. Furthermore, $v$ is a locally-maximal element or is in a local chain or is in special cases of types $F_4$ or $B_n$. In the first two cases, that is, when $v$ is a locally-maximal element or is in a local chain, then Proposition~\ref{prop:hypmain} allows us to replace $w$ with rationally smooth $u$ where $|S(u)|<|S(w)|$. If we are in the special cases, using Example~\ref{ex:F4special} and Lemma~\ref{lem:Bnspecial} combined with Corollary~\ref{cor:hypmain} allows us the same replacement.
\end{proof}

When $w\in W$ is rationally smooth, it is common for the polynomials $P_w(q) = R_w(q)$ to factor as a product of $q$-integers. If $P_w(q)$ factors into $q$-integers along Grassmannian BP decompositions, we say that $w$ has a \newword{chain BP decomposition} (this name comes from the fact that each poset $[e,v]^J$ is a chain).   By Corollary \ref{cor:BPonA} (\cite{Ga98}), all smooth permutations have chain BP decompositions.  If $w\in W$ has a chain BP decomposition, then the degrees of the $q$-integer factors of $P_w(q)=R_w(q)$ are strongly related to the structure of the corresponding inversion arrangement and are called \newword{exponents} of $w$.  In \cite{SL15}, Slofstra gives an explicit description of these exponents.  For other interesting results in inversion arrangements, we recommend that the reader take a look at \cite{ZIEGLER1993259}, \cite{HULTMAN20111897}, \cite{SL16}, \cite{Woo18}, \cite{Fan20}.

\section{Connections with the geometry of Schubert varieties}\label{S:Geometry}
In this section, we present results from \cite{RS16} which connects BP decompositions to the geometry of Schubert varieties.  Let $G$ be a connected simple Lie group over $\CC$ and fix a Borel subgroup $B$.  Let $W$ denote the Weyl group of $G$ with generating set $S$.  Since $G$ is a finite dimensional Lie group, the Weyl group $W$ is a Coxeter group of finite Lie type (See Figure \ref{fig:Coxeter diagrams}).  For any subset $J\subseteq S$, let $W_J$ the parabolic subgroup of $W$ generated by $J$ and let $P_J:=BW_JB$ denote the corresponding parabolic subgroup of $G$.  The coset space $G/P_J$ is called a \newword{partial flag variety} of $G$ and it is a smooth complex protective homogeneous space.  If $J=\emptyset$, then $P_\emptyset= B$ and $G/B$ is called the \newword{complete or full flag variety} of $G$.  

Consider the natural projection between flag varieties
$$\pi:G/B\rightarrow G/P_J$$
given by $\pi(gB)=gP_J$.  It is not hard to see that the fibers of the this map are $$\pi^{-1}(gP_J)=gP_J/B$$ and that the map $\pi$ gives a $(P_J/B)$-fiber bundle structure on the flag variety $G/B$ with base $G/P_J$.   For any $w\in W^J$, we define the \newword{Schubert variety} as the closure of the $B$-orbit (i.e. Schubert cell)
$$X^J(w):=\overline{BwP_J/P_J}.$$
If $J=\emptyset$, then $P_\emptyset=B$ and we will denote $X(w):=X^{\emptyset}(w)$.  
It is well known that if $w\in W^J$, then  
\[\dim_{\CC}(X^J(w))=\ell(w)\]
and Bruhat order describes the closure relations on Schubert cells
\[X^J(w)=\bigsqcup_{w'\in[e,w]^J} Bw'P_J/P_J.\]
In this section, we consider all Schubert varieties and do not restrict to the cases of smooth or rationally smooth.  If $w=vu$ is a parabolic decomposition with respect to $J$, then $u\in P_J$ and hence $\pi$ restricts to a projection between Schubert varieties
$$\pi:X(w)\rightarrow X^J(v).$$
The question we address is when does the map $\pi$ induce a fiber bundle structure on $X(w)$?  As we will see, the generic fibers of this map are isomorphic to the Schubert variety $X(u)$, however, unlike for $G/B$, the map $\pi$ restricted to a Schubert variety may not fiber bundle.  The following theorem is a geometric realization of BP-decompositions and is proved by Richmond and Slofstra in \cite[Theorem 3.3]{RS16}.

\begin{theorem}\label{thm:fiber_bundle}
Let $w=vu$ be a parabolic decomposition with respect to $J$.  Then the following are equivalent:
\begin{enumerate}
    \item The decomposition $w=vu$ is a $BP$ decomposition with respect to $J$.
    \item The projection $\pi:X(w)\rightarrow X^J(v)$ is a Zarisky-locally trivial with fiber $X(u).$
\end{enumerate}
\end{theorem}

Our goal is to give a detailed proof Theorem \ref{thm:fiber_bundle} following \cite{RS16}.  First, we need several important lemmas about Schubert varieties.  One key property needed in the proof of these results is the following well-known relation for double $B$-orbits for BN-pairs (or Tits systems).

\begin{lemma}\label{lem:B-orbit_relation}
Given $s\in S$ and $u\in W$, we have
\[
BsB\cdot BuB=\begin{cases}
BsuB & \text{if $s\notin D_L(u)$}\\
BuB\cup BsuB & \text{if $s \in D_L(u)$}\\
\end{cases}
\]
\end{lemma}
If $xP_J\in X^J(v)$, then we can write $xP_J=b_0v_0P_J$ for some $b_0\in B$ and $v_0\in[e,v]^J$.  The next lemma (\cite[Lemma 4.6]{RS16}) describes the fibers of the map $\pi$.

\begin{lemma}\label{lem:fibers}
Let $w=vu$ be the parabolic decomposition with respect to $J$ and $\pi:X(w)\rightarrow X^J(v).$  Let $xP_J\in X^J(v)$ and write $x=b_0v_0$ for some $b_0\in B$ and $v_0\in[e,v]^J$.  Then
$$\pi^{-1}(xP_J)=x\bigcup Bu'B/B$$
where the union is over all $u'\in W_J$ such that $v_0u'\leq w$.
\end{lemma}

\begin{proof}
We first look at the fiber over $xP_J$ of the map $\pi:G/B\rightarrow G/P_J.$ Note that $$P_J=B W_J B=\bigcup_{u'\in W_J} BuB$$ and hence the fiber of $xP_J$ in the full flag variety $G/B$ is
$$\pi^{-1}(xP_J)=b_0v_0\bigcup_{u'\in W_J} BuB/B.$$
Restricting the map $\pi$ to the Schubert variety $X(w)$ gives
$$\pi^{-1}(xP_J)=\left(b_0v_0\bigcup_{u'\in W_J} Bu'B/B\right) \cap \left(\bigcup_{w'\leq w}Bw'B/B\right).$$
Since $v_0\in W^J$ and $u'\in W_J$, Lemma \ref{lem:B-orbit_relation} implies $b_0v_0Bu'B\subseteq Bv_0u'B$.  Hence
$$\pi^{-1}(xP_J)=b_0v_0\bigcup Bu'B/B$$ where the union is over all $u'\in W_J$ such that $v_0u'\leq w.$
\end{proof}

The next lemma is from \cite[Proposition 4.7]{RS16}.

\begin{lemma}\label{lem:fiber_equi_dim}
Let $w=vu$ be a parabolic decomposition with respect to $J$.  Then the following are equivalent:
\begin{enumerate}

\item The decomposition $w=vu$ is a BP decomposition.
\item The fibers of the map $\pi:X(w)\rightarrow X^J(v)$ are isomorphic to $X(u)$.
\item The fibers of the map $\pi:X(w)\rightarrow X^J(v)$ are equidimensional.
\end{enumerate}
\end{lemma}
\begin{proof}
Clearly part (2) implies part (3), so we focus on showing part (1) implies part (2) and part (3) implies part (1).  Let $xP_J\in X^J(v)$ and write $x=b_0v_0$ for some $b_0\in B$ and $v_0\in[e,v]^J$.  If $w=vu$ is a BP decomposition, then Theorem \ref{thm:BP_characterization} part (3) implies that $v_0u'\leq w$ if and only if $u'\leq u$.  Lemma \ref{lem:fibers} implies that the fiber
\[\pi^{-1}(xP_J)=b_0v_0\bigcup_{u'\leq u}Bu'B/B=b_0v_0X(u)\]
and hence all fibers of $\pi$ are isomorphic to $X(u)$.

Now suppose all the fibers of $\pi$ are the same dimension.  Then Lemma \ref{lem:fibers} implies the fiber over the identity is $\pi^{-1}(eP_J)=X^J(u')$ where $u'$ denotes the maximal element of $[e,w]\cap W_J$.  Similarly, we have the fiber over $vP_J$ is $\pi^{-1}(vP_J)=vX^J(u)$.  Since the fibers are equidimensional, we have $\ell(u')=\ell(u)$.  But $u\leq u'$ and hence $u=u'$.  Thus $w=vu$ is a BP decomposition.
\end{proof}

What remains to be proved is that when $w=vu$ is a BP-decomposition, then the map $\pi:X(w)\rightarrow X^J(v)$ is locally trivial and hence a $X(u)$-fiber bundle.  We first need the following lemma which is proved in \cite[Lemma 4.8]{RS16}.

\begin{lemma}\label{lem:restricted_supports}
Let $v\in W^J$ and let $I=S(v)$ denote the support set of $v$.  Let $G_I\subseteq G$ denote the Levi subgroup of $P_I$.  Let $P_{I,J}:=G_I\cap P_J$ and $B_I:=G_I\cap B$ denote the corresponding Borel and parabolic subgroups of $G_I$.

Then the inclusion $i:G_I/P_{I,J}\hookrightarrow G/P_J$ induces an isomorphism
\[i:X_I^{I\cap J}(v)\rightarrow X^J(v)\]
where the Schubert variety \[X_I^{I\cap J}(v):=\overline{B_IvP_{I,J}/P_{I,J}}\subseteq G_I/P_{I,J}.\]
\end{lemma}

We can now prove the main theorem of the section.

\begin{proof}[Proof of Theorem \ref{thm:fiber_bundle}]
First observe that if $\pi:X(w)\rightarrow X^J(v)$ is a locally-trivial fiber bundle, then the fibers are equidimensional and hence $w=vu$ is a BP-decomposition by Lemma \ref{lem:fiber_equi_dim}.

Now suppose that $w=vu$ is a BP-decomposition and let $I=S(v)$.  Lemma \ref{lem:fiber_equi_dim} implies the fibers of the map $\pi$ are all isomorphic to $X(u)$ and hence we only need to show local triviality.  Lemma \ref{lem:restricted_supports} states that the inclusion $i:G_I/P_{I,J}\hookrightarrow G/P_J$ restricts to an isomorphism $i:X_I^{I\cap J}(v)\rightarrow X^J(v)$.  The map $G_I\rightarrow G_I/P_{I,J}$ is locally trivial and thus has local sections.  Hence for any $x\in X^J(v)$, there exists a Zariski open neighborhood $U_x\subseteq X^J(v)$ with a local section $s:U_x\rightarrow G_I\subseteq G$.  Define the multiplication map
\[m:U_x\times X(u)\rightarrow G/B\]
by $m(x',y):=s(x')\cdot y$.  We claim that the image of $m$ lies in the Schubert variety $X(w)$.  Let $x'\in U_x\subseteq X^J(v)$ and hence $x'\in Bv_0P_J$ for some $v_0\leq v$.  Thus we can write $s(x')=b_0v_0p_0$ for some $b_0\in B_I:=G_I\cap B$ and $p_0\in P_{I,J}$.  Since $w=vu$ is a BP-decomposition, Theorem \ref{thm:BP_characterization} implies that $I\cap J\subseteq D_L(u)$.  Since $P_{I,J}\subseteq P_{I\cap J}=BW_{I\cap J}B$, Lemma \ref{lem:B-orbit_relation} implies $p_0X(u)=X(u)$.  Hence
\[
m(x',X(u))=b_0v_0p_0X(u)=b_0v_0X(u)\subseteq X(v_0u)\subseteq X(w).
\]
Consider the commuting diagram:
    \begin{equation*}
        \xymatrix{ U_x \times X(u) \ar[r]^{m} \ar[d] & X(w) \ar[d]^{\pi} \\
                    U_x\ \ \ar @{^{(}->}[r] & X^J(v) \\ }
    \end{equation*}
and note that the map $m$ identifies $(x',X(u))$ with the fiber $\pi^{-1}(x')$.  For any $z\in \pi^{-1}(U_x)$, let $g_z:=s(\pi(z))\in G_I$.  Then $z\mapsto (\pi(z), g_z^{-1}z)$ maps $\pi^{-1}(U_x)$ to $U_x\times X(u)$ and is, in fact, the inverse of $m$.   This implies the map $m$ is an algebraic isomorphism and hence $\pi$ is locally trivial.

\end{proof}

One consequence of Theorem \ref{thm:fiber_bundle}, is the following cohomological interpretation of BP decompositions.  For any variety $X$, let $H^*(X)$ denote its singular cohomology with complex coefficients.

\begin{cor}
The decomposition $w=vu$ is a BP decomposition with respect to $J$ if and only if
\begin{equation}\label{eqn:cohom_modules}
H^*(X(w)\simeq H^*(X^J(v))\otimes H^*(X^*(u))
\end{equation}
as $H^*(X^J(v))$-modules.
\end{cor}

\begin{proof}
If $w=vu$ is a BP decomposition, then Equation \eqref{eqn:cohom_modules} follows from Theorem \ref{thm:fiber_bundle} and the Leray-Hirsch theorem.  Conversely, recall that \[P_w(q^2)=\sum_i^{2\ell(w)} \dim(H^i(X(w))\, q^i\quad
\text{and}\quad
P^J_v(q^2)=\sum_i^{2\ell(v)} \dim(H^i(X^J(v))\, q^i.\]
If Equation \eqref{eqn:cohom_modules} holds, then $w=vu$ is a BP decomposition since $P_w(q)=P^J_v(q)\cdot P_u(q)$.
\end{proof}

We give some remarks about the fibers if $\pi$ when $w=vu$ is not necessarily a BP decomposition.  The union describing general fibers in Lemma \ref{lem:fibers} is taken over all $u'$ such that $v_0u'\in [e,w]\cap v_0W_J$.  It is not difficult to see that this collection forms a lower order ideal in $W_J$.  In \cite{MM20, KLS14, OR22}, it is independently shown that these lower order ideals have unique maximal elements and hence are intervals in $W_J$.  This leads to the following corollary.

\begin{cor}\label{cor:coset_interval_structure}
Let $w=vu$ be a parabolic decomposition with respect to $J$ and $\pi:X(w)\rightarrow X^J(v).$  Let $xP_J\in X^J(v)$ and write $x=b_0v_0$ for some $b_0\in B$ and $v_0\in[e,v]^J$.  Then
$$\pi^{-1}(xP_J)=b_0v_0X(u_0)$$
where $u_0$ is the unique maximal element of the set $v_0^{-1}([e,w]\cap v_0W_J)$.

Moreover, if $u'$ denotes the maximal element of the set $[e,w]\cap W_J$, then $u\leq u_0\leq u'$.  If $w=vu$ is a BP decomposition, then $u=u_0=u'$.
\end{cor}

\begin{example}\label{Ex:forintro}
    Let $G=\SL_4(\CC)$.  Geometrically,  we have
    \begin{equation*}
        G/B=\{V_\bullet=(V_1\subset V_2\subset V_3\subset \CC^4)\ |\ \dim V_i=i\}.
    \end{equation*}
    Let $E_\bullet$ denote the flag corresponding to $eB$ and  $w=s_1s_2s_3s_2s_1$.  Then
    \begin{equation*}
        X(w) = \{V_{\bullet}\ |\ \dim(V_2\cap E_2)\geq 1\}.
    \end{equation*}
    We consider the geometric analogues of Examples \ref{example:12321BP} and \ref{example:12321not_BP}.
    
    First, if $J=\{s_1,s_3\}$, then $\pi(V_\bullet)=V_2$ and 
    $$w=vu=(s_1s_3s_2)(s_3s_1)$$
    is a BP decomposition with respect to $J$ as in Example \ref{example:12321BP}.  In particular, the Schubert variety
    \begin{equation*}
        X^J(v)= \{V_2\ |\ \dim(V_2\cap E_2)\geq 1\}
    \end{equation*}
    and the fibre over $V_2$ in the projection $\pi:X^{\emptyset}(w)\rightarrow X^J(v)$ is
    \begin{equation*}
        \pi^{-1}(V_2)=\{(V_1,V_3)\ |\ V_1\subset V_2\subset V_3\}\iso X(u)\iso
        \mathbb{CP}^1\times \mathbb{CP}^1.
    \end{equation*}
    
    If $J=\{s_1,s_2\}$, then $\pi(V_\bullet)=V_3$ and $$w=vu=(s_1s_2s_3)(s_2s_1)$$ is not a BP decomposition as in Example \ref{example:12321not_BP}.  The fiber over $V_3$ is given by
    \begin{align*}
        \pi^{-1}(V_3)&=\{(V_1,V_2)\ |\ V_1\subset V_2\subset V_3
        \text{ and } \dim(V_2\cap E_2)\geq 1\} \\
        & \iso \begin{cases} X(s_2s_1) & \text{if}\quad \dim(V_3\cap E_2)=1\\
                             X(s_1s_2s_1) & \text{if}\quad E_2 \subset V_3
                \end{cases}
    \end{align*}
    Note that the fibres are not equidimensional.
\end{example}

\begin{rem}    
\label{rem:BPnice}
Combinatorially, Corollary \ref{cor:coset_interval_structure} says that if $w=vu$ is a parabolic decomposition with respect to $J$ and $u'$ denotes the maximal element of $[e,w]\cap W_J$, then for every $v_0\in [e,v]^J$, the coset interval 
\[[e,w]\cap v_0W_J\simeq [e,u_0]\]
for some $u\leq u_0\leq u'$.  At the extremes, we have
\[[e,w]\cap W_J\simeq [e,u']\quad\text{and}\quad [e,w]\cap vW_J\simeq [e,u].\]
\end{rem}
\subsection{Relative BP decompositions}
We finish this section by stating a relative version of Theorem \ref{thm:fiber_bundle}.  In this case, we have two parabolic subgroups $P_J\subseteq P_K\subseteq G$ corresponding to subsets $J\subseteq K\subseteq S$.  Consider the projection \[\pi:G/P_J\rightarrow G/P_K\] and we can ask the question: when does the map $\pi$ induced a fiber bundle structure when restricted to the Schubert variety $X^J(w)$?  To answer to this equation, we define the relative version of a BP decomposition.  

\begin{defn}
Let $J\subseteq K\subseteq S$ and $w\in W^J$.  Let $w=vu$ denote the parabolic decomposition with respect to $K$.  We say $w=vu$ is a \textbf{BP decomposition with respect to $(J,K)$} if the Poincare polynomial factors \[P_w^J(q)=P_v^K(q)\cdot P_u^J(q).\]
\end{defn}

Note that if $J=\emptyset$, then this is the usual BP decomposition of $w$ with respect to $K$.  We remark that relative BP decompositions are characterized by a similar list of conditions to those given in Theorem \ref{thm:BP_characterization}.  See \cite[Proposition 4.2]{RS16} for a precise statement.

\begin{theorem}\label{thm:fiber_bundle_relative}
Let $w\in W^J$ and let $w=vu$ be a parabolic decomposition with respect to $(J,K)$.  Then the following are equivalent:
\begin{enumerate}
    \item The decomposition $w=vu$ is a $BP$ decomposition with respect to $(J,K)$.
    \item The projection $\pi:X^J(w)\rightarrow X^K(v)$ is a Zarisky-locally trivial with fiber $X^J(u).$
\end{enumerate}
\end{theorem} 

Theorem \ref{thm:fiber_bundle_relative} is proved in \cite{RS16} and the proof is very similar to that of Theorem \ref{thm:fiber_bundle}.

\begin{rem}
    Theorems \ref{thm:fiber_bundle} and \ref{thm:fiber_bundle_relative} hold for the much larger class of Kac-Moody Schubert varieties.  Kac-Moody groups are infinite dimensional generalizations of Lie groups and include the family of affine Lie groups.  Their Weyl groups are (not necessarily finite) crystalographic Coxeter groups.  While the flag varieties of Kac-Moody groups are also infinite dimensional, their Schubert varieties are finite dimensional.  For more on Kac-Moody flag varieties and their Schubert varieties see \cite{Ku02}.
\end{rem}

\section{Iterated BP decompositions and staircase diagrams}\label{s:iterated_BPs}

In this section, we discuss iterations of BP decompositions for Coxeter groups of finite type.  In particular, if $(W,S)$ is a Coxeter system and $J\subseteq S$, then each subgroup $W_J$ has a unique longest element we denote by $u_J$.  We begin with the following definition.

\begin{defn}
We say a factorization \[w=v_nv_{n-1}\cdots v_1\] is an \textbf{iterated BP decomposition} if $(v_{i+1})(v_i\cdots v_1)$ is a BP decomposition for each $1<i<n$.
\end{defn}

By Theorem \ref{thm:fiber_bundle}, iterated BP decompositions correspond to iterated fiber bundle structures on Schubert varieties.  

\subsection{Staircase diagrams}
In this section we combinatorially characterize iterated BP decompositions by objects called \emph{labelled staircase diagrams}.  Staircase diagrams are certain partially ordered sets over a given graph and were introduced by Richmond and Slofstra in \cite{RS17} with the goal of developing a combinatorial framework to study iterated BP decompositions.  We focus on staircase diagrams over the Coxeter graph of a Coxeter group.  The Coxeter graph is simply the Coxeter diagram of $W$ without the edge labels and we denote this graph by $\Gamma_W$ (See Figure \ref{fig:Coxeter diagrams}). In other words, $\Gamma_W$ is a graph with vertex set $S$ and edge set $\{(s,t)\in S^2 \ |\ m_{st}\geq 3\}$.  Note that the Coxeter groups of types $A_n$ and $B_n/C_n$ all have the same underlying Coxeter graph.  

Before stating the definition of a staircase diagram, we need some terminology.  Given $s,t\in S$, we say $s$ is \newword{adjacent} to $t$ if $(s,t)$ is an edge in $\Gamma_W$.  We say a subset $B\subset S$ is \newword{connected} if the induced subgraph of $B$ in $\Gamma_W$ is connected.  If $\mcD$ is a collection of subsets of $S$ and $s\in S$, we define
$$\mcD_s:=\{B\in \mcD\ |\ s\in B\}.$$
In other words, $\mcD_s$ are the elements in $\mcD$ that contain $s\in S$.

\begin{defn}\label{def:staircase_diagram}
Let $(W,S)$ denote a Coxeter system and let $\mcD$ be a collection of subsets of $S$.  We say a partially ordered set $(\mcD,\prec)$ is a \textbf{staircase diagram} if the following hold:
\begin{enumerate}
\item Every $B\in\mcD$ is connected, and if $B$ covers $B'$, then $B\cup B'$ is connected.
\item The subset $\mcD_s$ is a chain for every $s\in S$.
\item If $s$ is adjacent to $t$, then $\mcD_s\cup \mcD_t$ is a chain, and $\mcD_s$ and $\mcD_t$ are saturated subchains of $\mcD_s\cup \mcD_t$.
\item For every $B\in \mcD$, there exists $s\in S$ (resp. $s'\in S$) such that $B$ is the minimum in $\mcD_s$ (resp. maximum in $\mcD_{s'}$).
\end{enumerate}
\end{defn}

If the generating set $S=\{s_1,\ldots, s_n\}$, then we use interval notation $$[s_i,s_j]:=\{s_i,s_{i+1},\ldots,s_j\}$$ for $i\leq j$.  In type $A_{n}$, we have the Coxeter graph
  \[\begin{tikzpicture}[scale=.4]
    \draw[thick] (0,0)--(2,0)--(4.25,0)--(6.5,0);
    \draw[thick,fill=white] (0,0) circle (.3cm) node[label={[label distance=.1cm]-90:$s_1$}] { };
    \draw[thick,fill=white] (2,0) circle (.3cm) node[label={[label distance=.1cm]-90:$s_2$}]{ };
    \draw[white, fill=white] (4.25,0) circle (.75cm);
    \draw (4.25,0) node {$\cdots$};
    \draw[thick,fill=white] (6.5,0) circle (.3cm) node[label={[label distance=.1cm]-90:$s_{n}$}]{ };
  \end{tikzpicture}\]
An example of staircase diagram in this type is
$$\mcD=\{[s_1,s_3]\prec [s_2,s_4]\prec [s_3,s_5]\succ [s_6]\succ [s_7,s_9]\prec [s_9,s_{10}]\prec[s_{10},s_{11}]\}.$$
In this example, the set $\mcD_{s_3}=\{[s_1,s_3], [s_2,s_4], [s_3,s_5]\}.$   In Figure \ref{fig:staircase1}, we represent this staircase diagram with a picture of uneven steps where ``higher steps" are greater in the partial order:

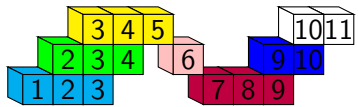
\begin{figure}[htbp]
       \begin{tikzpicture}[scale=0.4]
            \ppAff{12}{
                {0,0,5,5,5,0,0,0,1,1,1},
               {0,6,6,0,0,4,0,2,2,2},
                {7,7,0,0,0,0,3,3,3}}
       \end{tikzpicture}
\caption{Picture of a staircase diagram of type $A_{11}$.}\label{fig:staircase1}
\end{figure}

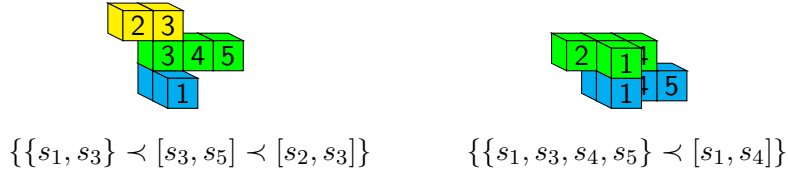
\begin{figure}[htbp]
    \begin{tikzpicture}[scale=0.4]
        \planepartitionD{5}{ {0,0,2,0}, {1,1,1,0}, {0,0,1,1}}
    \end{tikzpicture} \hspace{1.5in}
    \begin{tikzpicture}[scale=0.4]
        \planepartitionD{5}{ {1,1,2,0}, {0,1,2,1}}.
    \end{tikzpicture}
$$\{\{s_1,s_3\}\prec [s_3,s_5]\prec [s_2,s_3]\}\hspace{.5in}\{\{s_1,s_3,s_4,s_5\}\prec [s_1,s_4]\}$$
    \caption{A staircase diagrams of type $D_5$.}
    \label{F:type_D_diagrams}
\end{figure}

Since elements of a staircase diagram are connected, we will refer to them as ``blocks".  Note the blocks may not necessarily be ordered intervals.  In type $D_5$, we have Coxeter graph
\begin{equation*}
  \begin{tikzpicture}[scale=.4]
    \draw[thick] (0,0)--(2,0)--(4,0)--(6,0);
    \draw[thick] (2,0)--(2,2);
    \draw[thick,fill=white] (0,0) circle (.3cm) node[label={[label distance=.1cm]-90:$s_2$}] { };
    \draw[thick,fill=white] (2,0) circle (.3cm) node[label={[label distance=.1cm]-90:$s_3$}]{ };
    \draw[thick,fill=white] (4,0) circle (.3cm) node[label={[label distance=.1cm]-90:$s_{4}$}]{ };
    \draw[thick,fill=white] (6,0) circle (.3cm) node[label={[label distance=.1cm]-90:$s_{5}$}]{ };
    \draw[thick,fill=white] (2,2) circle (.3cm) node[label={[label distance=.1cm]0:$s_{1}$}]{ };
  \end{tikzpicture}
\end{equation*}
with examples of staircase diagrams in Figure \ref{F:type_D_diagrams}.

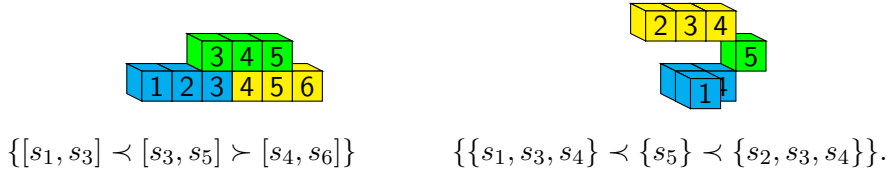
\begin{figure}[htbp]
\begin{tikzpicture}[scale=0.4]
    \ppAff{8}{
        {0,3,3,3,1,1,1,0},
        {0,0,2,2,2,0,0,0},}
\end{tikzpicture} \hspace{1.5in}
\begin{tikzpicture}[scale=0.4]
    \planepartitionD{5}{ {0,1,2},{1},{0,1,1,1}}
\end{tikzpicture}
$$
\{[s_1,s_3]\prec [s_3,s_5]\succ [s_4,s_6]\}\hspace{.5in}
\{\{s_1,s_3,s_4\}\prec \{s_5\}\prec \{s_2,s_3,s_4\}\}.
$$
\caption{Non-examples of staircase diagrams.}\label{F:SD_nonexamples}
\end{figure}

In Figure \ref{F:SD_nonexamples}, we give some non-examples of staircase diagrams.  The first diagram is of type $A_6$ violates parts (3) and (4) of Definition \ref{def:staircase_diagram}.  The second diagram is of type $D_5$ and violates part (2) of Definition \ref{def:staircase_diagram}.  

It is not hard to check that Definition \ref{def:staircase_diagram} is symmetric with respect to the partial order.  Given a staircase diagram $\mcD$, we can define the dual staircase diagram $\flip(\mcD)$ to be the set $\mcD$ with the reverse partial order.  Pictorially this corresponds to ``flipping" the staircase from top to bottom:

\begin{figure}[h]
\begin{tikzpicture}[scale=0.4]
    \ppAff{8}{
        {0,3,0,0,1,1,1,0},
        {0,0,2,2,2,2,0,0},}
\end{tikzpicture}\hspace{1in}
\begin{tikzpicture}[scale=0.4]
    \ppAff{8}{
        {0,0,2,2,2,2,0,0},
        {0,3,0,0,1,1,1,0},}
\end{tikzpicture}
$$\mcD=\{[s_1,s_3]\prec[s_2,s_5]\succ\{s_6\}\}\hspace{.3in}\flip(\mcD)=\{[s_1,s_3]\succ[s_2,s_5]\prec\{s_6\}\}.$$
    \caption{A staircase diagrams $\mcD$ and $\flip(\mcD)$.}
    \label{F:flip_diagram}
\end{figure}
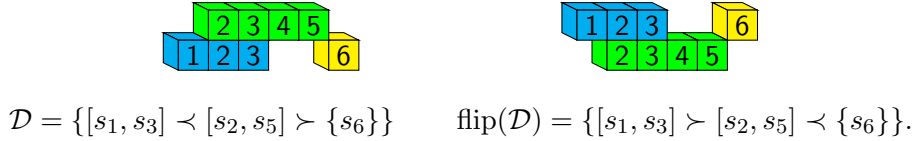

If $\mcD'$ is a saturated subset of $\mcD$, then the induced partial order on $\mcD'$ makes it a staircase diagram.  In this case, we say $\mcD'$ is a \emph{subdiagram} of $\mcD$.
For any $J\subseteq S$, define
$$\displaystyle\mcD_J:=\{B\in\mcD \ | \ J\subseteq B\}.$$
The following lemma describes some combinatorial properties of staircase diagrams.
\begin{lemma}\label{lem:staircase_lemma}
Let $\mcD$ be a staircase diagram of a Coxeter system $(W,S)$.  Then:
\begin{enumerate}
\item For any $J\subseteq S$, the set $\mcD_J$ is a chain in $\mcD$.
\item If $B,B'\in \mcD$, then $B\nsubseteq B'$.
\item If $B,B'\in \mcD$ and $B\cup B'$ is connected, then $B$ and $B'$ are comparable.
\end{enumerate}
\end{lemma}

\begin{proof}
Part (1) follows from the fact that $\mcD_J$ is the intersection of $\mcD_s$ where $s\in J$ and each $\mcD_s$ is a chain. 

For part (2), select $s, s'\in S$ such that $B$ is the maximal and minimal block of $\mcD_s$ and $\mcD_{s'}$ respectively.  Then $\mcD_{\{s,s'\}}$ consists only of $B$.  If $B\subseteq B'$, then $B'\in\mcD_{\{s,s'\}}$ and hence $B=B'$.  

For part (3), if $B\cup B'$ is connected, then there exist $s\in B$ and $t\in B'$ such that $s$ is adjacent to $t$.  Thus $B,B'$ belong to the chain $\mcD_s\cup \mcD_t$ and hence $B,B'$ are comparable.

\end{proof}

\subsection{Labellings of staircase diagrams}

\begin{figure}[htbp]
\begin{tikzpicture}[scale=0.4]
    \ppAff{9}{
        {0,3,3,0,0,1,1,1,0},
        {0,0,2,9,9,2,2,0,0},}
\end{tikzpicture}
\caption{The set $J_R([s_2,s_6])=\{s_2,s_3,s_6\}$ highlighted in green.}\label{fig:J_rset}
\end{figure}
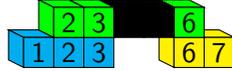

Staircase diagrams provide the framework for building iterated BP decompositions.  Let $\mcD$ be a staircase diagram.  For any $B\in\mcD$, define the sets

\[J_R(B):=B\cap \left(\bigcup_{B'\prec B} B'\right)\quad\text{and}\quad J_L(B):=B\cap \left(\bigcup_{B'\succ B} B'\right).\]
Pictorially, we can think of the set $J_R(B)$ as the elements of $B$ that are ``covered below" by other blocks in $\mcD$ and $J_L(B)$ as the elements of $B$ that are ``covered above".  For example, if $$\mcD=\{[s_1,s_3]\prec [s_2,s_6]\succ [s_6,s_7]\}$$
then $$J_R([s_2,s_6])=\{s_2,s_3,s_6\}$$
which we highlight in Figure \ref{fig:J_rset}.

We define a labelling of a staircase diagram which assigns a Coxeter group element to each block in $\mcD$.  For any $J\subseteq S$, we let $u_J$ denote the longest element of $W_J$.

\begin{defn}
Let $\mcD$ be a staircase diagram on a Coxeter system $(W,S)$.  We say a function
\[\lambda:\mcD\rightarrow W\]
is a \textbf{labelling} of $\mcD$ if for every $B\in \mcD$, we have
\begin{enumerate}
\item $J_R(B)\subseteq D_R(\lambda(B))$,
\item $J_L(B)\subseteq D_L(\lambda(B))$, and
\item $S(\lambda(B)u_{J_R(B)})=B=S(u_{J_L(B)}\lambda(B))$.
\end{enumerate}
We denote a labeled staircase diagram by the pair $(\mcD,\lambda)$.
\end{defn}

\begin{example}
The function $\lambda:\mcD\rightarrow W$ given by $\lambda(B)=u_B$ is a labelling of $\mcD$.  This labelling is called the maximal labelling of $\mcD$.
\end{example}

Note that while staircase diagrams of type $A_n$ and $B_n/C_n$ are the same, labelled staircase diagrams are different since they depend on the group $W$ and not just the underlying graph $\Gamma_W$.
The definition of a labelling is compatible with the dual of staircase diagram.  For any labeled staircase diagram $(\mcD,\lambda)$, define the inverse labelling
\[\lambda^{-1}:\flip(\mcD)\rightarrow W\]
by $\lambda^{-1}(B):=\lambda(B)^{-1}$.  It is easy to check that $(\flip(\mcD),\lambda^{-1})$ is also a labelled staircase diagram.  The condition $J_R(B)\subseteq D_R(\lambda(B))$ implies that $\lambda(B)u_{J_R(B)}$ is the minimal right coset representative of $\lambda(B)$ in $W^{J_R(B)}$.   Similarly, we have that $u_{J_L(B)}\lambda(B)$ is a minimal left coset representative of $\lambda(B)$.  These coset representatives play an important role in the next definition, so for any labelled staircase diagram $(\mcD,\Lambda)$ and $B\in\mcD$, we define
\[\overline{\lambda}(B):=\lambda(B)u_{J_R(B)}.\]
\begin{defn}\label{defn:Lambda_staircase}
Given a labeled staircase diagram $(\mcD,\lambda)$ define
\[\Lambda(\mcD,\lambda):=\overline{\lambda}(B_n)\overline{\lambda}(B_{n-1})\cdots\overline{\lambda}(B_1)\]
where $B_1,\ldots,B_n$ is some linear extension of the poset $\mcD$.  If the labelling $\lambda$ is clear from the context, then we will denote $\Lambda(\mcD)=\Lambda(\mcD,\lambda)$.
\end{defn}

By part (3) of Lemma \ref{lem:staircase_lemma}, if $B_i$ and $B_j$ are not comparable, then $\overline{\lambda}(B_i)$ and $\overline{\lambda}(B_j)$ have commuting supports and hence commute as elements of $W$.  This implies that $\Lambda(\mcD)$ is independent of choice of linear extension and is well defined.

\begin{example}
Let $\mcD=\{[s_1,s_3],[s_5,s_6],[s_2,s_5]\}$ in type $A_6$.
$$
\begin{tikzpicture}[scale=0.4]
    \ppAff{8}{
        {0,3,3,0,1,1,1,0},
        {0,0,2,2,2,2,0,0},}
\end{tikzpicture}
$$
Then
$$\overline{\lambda}([s_1,s_3])=s_1s_2s_3s_1s_2s_1,\quad \overline{\lambda}([s_5,s_6])=s_5s_6s_5,$$
$$\overline{\lambda}([s_2,s_5])=(s_3s_2s_4s_3s_5s_4\text{\redt{$s_5s_2s_3s_2$}})(\text{\redt{$s_2s_3s_2s_5$}})=s_3s_2s_4s_3s_5s_4$$
and
$$\Lambda(\mcD)=(s_3s_2s_4s_3s_5s_4)(s_5s_6s_5)(s_1s_2s_3s_1s_2s_1).$$
Highlighted in red is the element $u_{\{s_2,s_3,s_5\}}$ since $J_R([s_2,s_5])=\{s_2,s_3,s_5\}.$
\end{example}
Note that if $\lambda:\mcD\rightarrow W$ is a maximal labelling, then $\overline{\lambda}(B)$ is the maximal element $W_B\cap W^{J_R(B)}$
We define the support of $\mcD$ to be the set \[S(\mcD):=\bigcup_{B\in\mcD} B.\]  Note that if $\lambda:\mcD\rightarrow W$ is a labeling, then $S(\mcD)=S(\Lambda(\mcD))$.
Furthermore, since the support set
\[S(\overline{\lambda}(B_{i-1})\cdots\overline{\lambda}(B_1))=B_1\cup\cdots\cup B_{i-1}\]
is disjoint with  $B\setminus J_R(B_i)$, the product
\[\overline{\lambda}(B_i)\cdot(\overline{\lambda}(B_{i-1})\cdots\overline{\lambda}(B_1))\]
is a parabolic decomposition with respect to $B_1\cup\cdots\cup B_{i-1}$.
We will show that this decomposition is in fact a BP decomposition and thus the factorization of $\Lambda(\mcD)$ in Definition \ref{defn:Lambda_staircase} corresponds to an iterated BP decomposition.  The next lemma gives several properties on how the Coxeter theoretic data of the element $\Lambda(\mcD)$ is extracted from the combinatorial data of the staircase diagram $\mcD$.

\begin{lemma}\label{lem:labeled_SD_properties}
Let $(\mcD,\lambda)$ be a labeled staircase diagram.  Then the following are true:

\begin{enumerate}
\item $\Lambda(\mcD)^{-1}=\Lambda(\flip(\mcD), \lambda^{-1})$.

\smallskip

\item The right descents of $\Lambda(\mcD)$ consist of all $s\in S(\mcD)$ that satisfy:
\begin{enumerate}
\item $\min(\mcD_s)\preceq\min(\mcD_t)$ for all $t$ adjacent to $s$ and
\item $s$ is a right descent of $\lambda(\min(\mcD_s))$.
\end{enumerate}

\smallskip

\item The left descents of $\Lambda(\mcD)$ consists of all $s\in S(\mcD)$ that satisfy:
\begin{enumerate}
\item $\max(\mcD_s)\succeq\max(\mcD_t)$ for all $t$ adjacent to $s$ and
\item $s$ is a left descent of $\lambda(\max(\mcD_s))$.
\end{enumerate}

\smallskip

\item Let $\mcD'$ be a lower order ideal in $\mcD$ and let $\mcD'':=\mcD\setminus\mcD.$  Then
\[\Lambda(\mcD)=(\Lambda(\mcD'')u_K)\cdot \Lambda(\mcD')\]
is a parabolic decomposition with respect to $S(\mcD')$ where \[K=\{s\in S\ |\ \min(\mcD''_s)\neq \min(\mcD_s)\}.\]
\end{enumerate}
\end{lemma}

The proof of Lemma \ref{lem:labeled_SD_properties} is technical, so we refer the reader to \cite{RS17} for more details.  Observe that part (3) follows from parts (1) and (2).

\begin{theorem}\label{thm:iterated_BP}
Let $\mcD$ be a staircase diagram with a linear extension $B_1,\ldots,B_n$.  For $i\geq 2$, let $\mcD^i$ denote the subdiagram
\[\mcD^i:=\{B_1,\ldots,B_{i-1}\}.\]
If $\lambda$ is a labelling of $\mcD$, then
\[\Lambda(\mcD^{i+1})=\overline{\lambda}(B_i)\cdot\Lambda(\mcD^i)\]
is a BP decomposition with respect to $S(\mcD^i)$.
\end{theorem}

\begin{proof}[Proof of Theorem \ref{thm:iterated_BP}]
Lemma \ref{lem:labeled_SD_properties} part (4) implies $\overline{\lambda}(B_i)\cdot\Lambda(\mcD^i)$ is a parabolic decomposition, so it suffices to show that the decomposition satisfies the BP condition in Theorem \ref{thm:BP_characterization} part (4).  First observe that
\[S(\overline{\lambda}(B_i))\cap S(\mcD^i)=B_i\cap S(\mcD^i)=J_R(B_i).\]
Thus $\overline{\lambda}(B_i)\cdot\Lambda(\mcD^i)$ is a BP decomposition if and only if $J_R(B_i)\subseteq D_L(\Lambda(\mcD^i))$.  Let $s\in J_R(B_i)$.  We use the characterization given in Lemma \ref{lem:labeled_SD_properties} part (3) to show that $s\in D_L(\Lambda(\mcD^i))$ .   Suppose that $t\in S(\mcD^i)$ is adjacent to $s$.   Observe that if $B_j$ the predecessor of $B_i$ in the chain $\mcD_s$, then $B_j=\max(\mcD^i_s)$.  By definition of staircase diagram, $\mcD_s$ is a saturated subchain of the chain $\mcD_s\cup\mcD_t$.  Since $B\preceq B_i$ for all $B\in\mcD^i$, it follows that $\max(\mcD^i_t)\preceq B_j$.  By Lemma \ref{lem:labeled_SD_properties}, it remains to show that $s\in D_L(\lambda(B_j))$.  Since $s\in B_i\cap B_j$, we have $s\in J_L(B_j)$ and, by the definition of a labelling, $J_L(B_j)\subseteq\lambda(B_j)$.  Thus $s\in D_L(\Lambda(\mcD^i))$ which completes the proof.
\end{proof}

\subsection{Complete BP decompositions}\label{S:completeBP}

In this section we discuss a special class of decompositions called complete BP decompositions.  We start with the following definition which was introduced in Section \ref{S:rat_smooth_intro}.

\begin{defn}
A BP composition $w=vu$ with respect to $J$ is a \textbf{Grassmannian BP decomposition} if $|J|=|S(w)|-1$.  In other words, $J$ is maximal proper subset of $S(w)$.
\end{defn}

Geometrically, Grassmannian BP decompositions correspond to projections $\pi: G/B\rightarrow G/P$ where $P$ is taken to be a maximal parabolic.  In the classical type $A$ setting this partial flag variety $G/P$ corresponds to a Grassmannian variety.  If $w=vu$ is a Grassmannian BP decomposition, then Theorem \ref{thm:fiber_bundle} implies the Schubert variety $X(w)$ is an $X(u)$-fiber bundle over the Grassmannian Schubert variety $X^J(v)$.  Note that the decompositions that arise in Theorems \ref{thm:smooth_perm} and \ref{thm:BPgen} are Grassmannian BP decompositions.

\begin{defn}
Let $n=|S(w)|$.  We say \[w=v_nv_{n-1}\cdots v_1\] is a \textbf{complete BP decomposition} if $(v_{i+1})(v_i\cdots v_1)$ is a Grassmannian BP decomposition for each $1\leq i<n$.
\end{defn}

Complete BP decompositions are iterated BP decompositions where the number of non-trivial factors is maximized in the sense that each iteration adds exactly one additional generator to the support set of $w$.  For example, in type $A_3$, we have that \[w=(s_1s_2s_3)(s_1s_2)(s_1)\] is a complete BP decomposition of the longest element.  Note that these decompositions are not unique.  For $w$ above, the decomposition
\[w=(s_2s_1s_3s_2)(s_1)(s_3)\]
is also a complete BP decomposition.  The goal of this section is to classify which elements $w\in W$ that have complete BP decompositions.  The key to this classification is the notion of nearly-maximal elements.

\begin{defn}\label{def:nearly_max}
We say an element $w\in W$ is \textbf{nearly-maximal} if there is a Grassmannian BP decomposition $w=vu$ such that $S(u)\subset S(v)$.

Furthermore, we say a labelled staircase diagram $(\mcD,\lambda)$ is \textbf{nearly-maximal} if each $B\in\mcD$, $\lambda(B)$ is nearly-maximal.
\end{defn}

If $w=vu$ is nearly-maximal, then \[S(u)\subseteq S(v)\cap J\subseteq D_L(u)\] 
and hence $S(u)=D_L(u)$.  This implies that $u$ is the maximal element of $W_J$.  Geometrically, this corresponds to the fiber $X(u)$ being isomorphic to the flag variety $P_J/B$.  Not all Grassmannian BP decompositions satisfy the nearly-maximal condition.  For example, in type $A_4$, \[w=(s_1s_2)(s_1s_3s_4)\]
is Grassmannian BP decomposition with respect to $J=\{s_1,s_3,s_4\}$, but $w$ is not nearly maximal.  Note that the maximal labelling of a staircase diagram is nearly-maximal.  The importance of nearly-maximal labelings is that they can used to construct complete BP decompositions.  In fact, this construction will yield the following bijection:

\begin{theorem}\label{thm:complete_BP_bijection}
Let $W$ be a Coxeter group.  Then the map $(\mcD,\lambda)\mapsto\Lambda(\mcD)$ defines a bijection between staircase diagrams over $W$ with a nearly-maximal labelling $\lambda$, and elements of $W$ with a complete BP decomposition.
\end{theorem}

\begin{proof}[Outline of the proof of Theorem \ref{thm:complete_BP_bijection}]
First note by Theorem \ref{thm:iterated_BP} and Definition \ref{def:nearly_max}, if $\lambda$ is a nearly maximal labelling of $\mcD$, then $\Lambda(\mcD)$ has a complete BP-decomposition and thus the map $(\mcD,\lambda)\mapsto\Lambda(\mcD)$ is well defined.

To show that the map is injective, suppose we have two nearly-maximal labelled staircase diagrams $(\mcD_1,\lambda_1)$ and $(\mcD_2,\lambda_2)$ such that $\Lambda(\mcD_1)=\Lambda(\mcD_2)$.  Choose $s\in S$ such that $\Lambda(\mcD_i)=vu$ is a BP decomposition with respect to $J=S\setminus\{s\}$.  It can be shown that $B:=S(v)$ is a maximal block of $\mcD_i$ and hence, by induction on the number of blocks, $\mcD_1=\mcD_2$.  To show that $\lambda_1=\lambda_2$, note that, by Lemma \ref{lem:labeled_SD_properties} part (1), the parabolic decomposition of $\Lambda(\mcD_i)^{-1}$ with respect to $B$ is given by
 \[\Lambda(\mcD_i)^{-1}=v'\cdot \lambda_i(B)^{-1}\]
 for some $v'$ and thus $\lambda_1(B)=\lambda_2(B)$.  We also have
 \[J_R(B,\mcD_i)=B\cap S(\mcD_i\setminus\{ B\})=B\cap S(\lambda_i(B)\cdot \Lambda(\mcD_i))\]
 and hence $J_R(B,\mcD_1)=J_R(B,\mcD_2)$.
This implies $\overline{\lambda}_1(B)=\overline{\lambda}_2(B)$.  Since $\Lambda(\mcD_1)=\Lambda(\mcD_2)$ and $\overline{\lambda}_1(B)=\overline{\lambda}_2(B)$, we have that the induced labelling on lower order ideals satisfies $\Lambda(\mcD_1\setminus\{B\})=\Lambda(\mcD_2\setminus\{B\})$.  By induction on $|\mcD_i|$, we have that the labellings $\lambda_1=\lambda_2$.

To show that the map is surjective, suppose $x\in W$ has a complete BP decomposition $x=v_n\cdots v_1$.  By induction, suppose that $(\mcD,\lambda)$ is a nearly-maximal labelled staircase diagram such that $\Lambda(\mcD)=v_{n-1}\cdots v_1$.  Define the staircase diagram
\[\tilde\mcD:=\mcD^0\cup\{S(v_n)\}\quad\text{where}\quad \mcD^0:=\{B\in\mcD\ |\ B\nsubseteq S(v_n)\}\]
with the added covering relations $\max(\mcD^0_s)\prec S(v_n)$ for every $s\in S(\mcD^0)$ contained in, or adjacent to $S(v_n)$.  It can be shown that $\tilde\mcD$ satisfies Definition \ref{def:staircase_diagram} of staircase diagram.  Finally, define the labelling $\tilde\lambda:\tilde\mcD\rightarrow W$ by
\[\tilde\lambda(B):=\begin{cases} \lambda(B) & \text{if $B\in\mcD^0$}\\  v_n\cdot u_{S(v_n)\cap S(\mcD)} &\text{if $B=S(v_n)$}.\end{cases}\]
Again, it can shown that $\tilde\lambda$ is a nearly maximal labelling of $\tilde\mcD$ such that $\Lambda(\tilde\mcD,\tilde\lambda)=x$.  This completes the proof.
\end{proof}

Next we apply Theorem \ref{thm:complete_BP_bijection} to rationally smooth elements of Coxeter groups of finte Lie-type.  The following rephrasing of Theorem \ref{thm:BPgen}.

\begin{theorem}\label{thm:BPgen_ver_2}
Let $w\in W$ be rationally smooth with $|S(w)|\geq 2$.  Then either $w$ or $w^{-1}$ has a Grassmannain BP decomposition $vu$ with respect to $J=S(w)\setminus\{s\}$ such that $s$ is a leaf in the Coxeter diagram of $W_{S(w)}$ and $vu_{S(v)\cap J}$ is nearly maximal.
\end{theorem}

The proof of Theorem \ref{thm:BPgen_ver_2} follows from checking that the list of elements given in Theorem \ref{thm:BPgen} all satisfy the definition of nearly maximal given in Definition \ref{def:nearly_max}.  We remark that there exist nearly maximal elements that are not rationally smooth.  Hence Theorem \ref{thm:BPgen_ver_2} is slightly weaker statement than Theorem \ref{thm:BPgen}.  In \cite{RS17}, Richmond and Slofstra define the stronger condition of ``almost-maximal" to make these theorems equivalent.  Our next goal is to give an outline of a proof of Theorem \ref{thm:BPexists} which states that rationally smooth elements always have Grassmannian BP decompositions.

\begin{theorem}\label{thm:smooth_Grass_BP}
Let $w\in W$ be (rationally) smooth.  Then there exists a Grassmannian BP decomposition $w=vu$ with respect to some maximal proper subset $J=S(w)\setminus\{s\}$.

Moreover, $u$ is (rationally) smooth and $v$ is (rationally) smooth with respect to $J$.
\end{theorem}

\begin{proof}[Outline of the proof of Theorem \ref{thm:BPexists}/ \ref{thm:smooth_Grass_BP}]
First note that if $w=vu$ is a BP decomposition, then Theorem \ref{thm:fiber_bundle} implies that if $X(w)$ is (rationally) smooth, then both $X(u)$ and $X^J(v)$ are also (rationally) smooth.  

Recall that Theorem \ref{thm:BPgen} states that if $w$ is rationally smooth, then either $w$ or $w^{-1}$ has a Grassmannian BP decomposition with respect to $J=S(w)\setminus\{s\}$ for some leaf $s\in S(w)$ in the Coxeter diagram of $W_{S(w)}$.  If $w$ has such a BP decomposition, then the theorem is proved.  Now suppose $w^{-1}$ has such a BP decomposition and hence we can write $w=uv$ where $u\in W_J$ and $v^{-1}\in W^J$ and $w^{-1}=v^{-1}u^{-1}$ is a Grassmannian BP decomposition with respect to $J=S(w)\setminus\{s\}$.  Since $w$ is (rationally) smooth, we have that $w^{-1}$ is (rationally) smooth and hence $u$, $u^{-1}$ are also (rationally) smooth.  Since $|S(u)|<|S(w)|$, we can inductively assume that there exists a Grassmannian BP decomposition $u=v'u'$ with respect to some maximal proper set $J'=J\setminus\{s'\}$.  It can shown that $s'\in J$ can be selected appropriately so that
\[w=v'(u'u)\] is a Grassmannian BP decomposition with respect to $S(w)\setminus\{s'\}$.
\end{proof}

\begin{cor}\label{cor:smooth_complete_BP}
If $w\in W$ is (rationally) smooth, then $w$ has a complete BP-decomposition.  In particular, there exists a staircase diagram $\mcD$ over $W$ and nearly-maximal labelling $\lambda$ such that $\Lambda(\mcD)=w$.
\end{cor}

We say a nearly maximal labelling $\lambda:\mcD\rightarrow W$ is (rationally) smooth if $\Lambda(\mcD)$ is (rationally) smooth.  In fact, if $\lambda$ is rationally smooth, then for each $B\in \mcD$, the element $\overline{\lambda}(B)$ must correspond to one of the elements in the list found in Theorem \ref{thm:BPgen}.  In particular, if $W$ is simply laced, then $\lambda$ must be the maximal labelling.  This implies the following corollary.

\begin{cor}\label{cor:smooth_staircase_bijection}
Let $W$ be a simply laced of finite type.  Then there is a bijection between staircase diagrams over $\Gamma_W$ and smooth elements of $W$.
\end{cor}

\begin{proof}
Let $\mcD$ be a staircase diagram and let $\lambda:\mcD\rightarrow W$ denote the maximal labeling.  Then the Schubert variety $X(\Lambda(\mcD))$ is an iterated fiber bundle of smooth Schubert varieties and hence smooth.  Conversely, if $X(w)$ is smooth, then Theorem \ref{thm:complete_BP_bijection} and Corollary \ref{cor:smooth_complete_BP} imply there is a unique smoothly labelled staircase diagram $(\mcD,\lambda)$ such that $\Lambda(\mcD)=w$.  Since $W$ is simply laced, Theorem \ref{thm:BPgen} implies $\lambda$ is the maximal labelling.
\end{proof}

\subsection{Enumerating smooth Schubert varieties}\label{S:Enumeration}
An application of Theorem \ref{thm:complete_BP_bijection} and Corollary \ref{cor:smooth_staircase_bijection} is that we can enumerate smooth Schubert varieties by counting staircase diagrams.  We give an overview of this enumeration in type $A$.  Recall that the Coxeter graph of type $A_n$ is a path on $n$ vertices:
\begin{equation*}
  \begin{tikzpicture}[scale=.4]
    \draw[thick] (0,0)--(2,0)--(4.25,0)--(6.5,0);
    \draw[thick,fill=white] (0,0) circle (.3cm) node[label={[label distance=.1cm]-90:$s_1$}] { };
    \draw[thick,fill=white] (2,0) circle (.3cm) node[label={[label distance=.1cm]-90:$s_2$}]{ };
    \draw[white, fill=white] (4.25,0) circle (.75cm);
    \draw (4.25,0) node {$\cdots$};
    \draw[thick,fill=white] (6.5,0) circle (.3cm) node[label={[label distance=.1cm]-90:$s_{n}$}]{ };
  \end{tikzpicture}
\end{equation*}
We will denote this graph by $\Gamma_n$.

\begin{theorem}\label{thm:smooth_enumeration_typeA}
Let $a_n$ denote the number of staircase diagrams over $\Gamma_n$ (equivalently, the number of smooth permutations in $\mfS_{n+1}$) and define the generating function
\[A(x):=\sum_{n=0}^{\infty} a_n\, t^n.\]
Then
\[A(x)=\frac{1-5t+4t^2+t\sqrt{1-4t}}{1-6t+8t^2-4t^3}.\]
\end{theorem}
A proof of Theorem \ref{thm:smooth_enumeration_typeA} first appeared in an unpublished paper by Haiman \cite{Ha92}.  The first published proof of Theorem \ref{thm:smooth_enumeration_typeA} is due to Bousquet-M\'elou and Butler in \cite{BMB07}.  In this section, we provide an alternate proof using staircase diagrams from \cite{RS17} and \cite{RS18}.
We first focus on diagrams that are chains.  

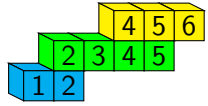
\begin{figure}  
\begin{tikzpicture}[scale=0.4]
        \ppAff{7}{
            {0,0,0,0,1,1},
            {0,2,2,2,2},
            {3,3,3,0}}
    \end{tikzpicture}
\caption{An increasing staircase diagram of type $A_6$.}\label{Fig:increasing_diagram}
\end{figure}

We say a staircase diagram $\mcD$ is \newword{increasing} over $\Gamma_n$ if $\mcD$ is fully supported (i.e. $S(\mcD)=\{s_1,\ldots,s_n\})$ and if for every $B,B'\in \mcD$ such that $s_i\in B$ and $s_j\in B'$ with $i<j$, we have $B\preceq B'$.   Pictorially, increasing staircase diagrams are represented by a sequence of blocks that are ``going up" from left to right with no gaps.  For example, $\mcD=\{[s_1,s_2]\prec [s_2,s_5]\prec [s_4,s_6]\}$ is increasing over $\Gamma_6$ as in Figure \ref{Fig:increasing_diagram}.   We say that $\mcD$ is \newword{decreasing} over $\Gamma_n$ if $\flip(\mcD)$ is increasing over $\Gamma_n$.

\begin{lemma}\label{lem:staircase_catalan}
The number of increasing staircase diagrams over $\Gamma_n$ is the $n$-th Catalan number.
\end{lemma}

\begin{proof}
We show that increasing staircase diagrams over $\Gamma_n$ are in bijection with Dyck paths.  Let \[\mcD=\{B_1\prec B_1\prec \cdots \prec B_m\}\] be such a diagram.  For each $B_i\in\mcD$ define
\begin{equation*}
    r(B_i):=\#\{s\in B_i\setminus B_{i-1}\}\qquad\text{and}\qquad
    u(B_i):=\#\{s\in B_i\setminus B_{i+1}\}
\end{equation*}
where we set $B_0=B_{m+1}=\emptyset.$  Let $P(\mcD)$ denote the lattice path in $\ZZ^2$ from $(0,0)$ to $(n,n)$ which takes $r(B_1)$ steps to the right, then $u(B_1)$ steps going up, followed by $r(B_2)$ steps to the right, then $u(B_2)$ steps going up and so forth (See Example \ref{Ex:Dyckpath_bijection}).  It is easy to check that $P(\mcD)$ is a Dyck path that stays below the diagonal in $\ZZ^2$.  One can also check that the map $P$ is invertible and hence a bijection.
\end{proof}

\begin{example}\label{Ex:Dyckpath_bijection}
    Consider the staircase diagram $\mcD=({s_1}\prec[s_2,s_5]\prec[s_4,s_6])$ on $\Gamma_6$.  The sequence of pairs $(r(B_i),u(B_i))$ is $((1,1),(4,2),(1,3))$ and corresponding Dyck path $P(\mcD)$ is given in Figure \ref{fig:SD_to_Dyckpath}.
\begin{figure}[htbp]
    \begin{tikzpicture}[scale=0.4]
        \ppAff{7}{
            {0,0,0,0,0,1},
            {0,2,2,2,2},
            {3,3,3,0}}
        \draw[<->, shift={(2.5,2)},black,line width=.25mm](0,0)--(2,0);
        \draw[shift={(6,-1)},step=1.0,black] (0,0) grid (6,6);
        \draw[shift={(6,-1)},red, line width=.75mm] (0,0)--(1,0)--(1,1)--(5,1)--(5,3)--(6,3)--(6,6);
    \end{tikzpicture}
    \caption{Bijection between increasing staircase diagrams and Dyck paths.}\label{fig:SD_to_Dyckpath}
\end{figure}
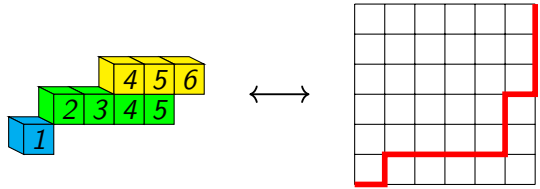
\end{example}

We use the enumeration of increasing staircase diagrams as the starting point to enumerate general staircase diagrams of type $A$.  The next step is to decompose staircase diagram with connected support into a smaller staircase diagram and an increasing/decreasing ``diagram" as follows:

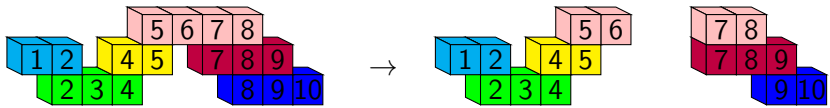
\begin{figure}[h]
\begin{equation*}
        \begin{aligned}
            \begin{tikzpicture}[scale=0.4]
            \ppAff{11}{
                {6,6,6,0,0,0,2,2,2,0},
                {0,5,5,5,0,3,3,0,1,1},
                {0,0,4,4,4,4}}
            \end{tikzpicture}
        \end{aligned} \quad \rightarrow\quad
        \begin{aligned}
            \begin{tikzpicture}[scale=0.4]
            \ppAff{7}{
                {0,0,2,2,2,0},
                {0,3,3,0,1,1},
                {4,4}}
            \end{tikzpicture}
        \end{aligned} \qquad
        \begin{aligned}
        \begin{tikzpicture}[scale=0.4]
            \ppAff{11}{
                {6,6,0},
                {0,5,5,5},
                {0,0,4,4}}
        \end{tikzpicture}
        \end{aligned}
\end{equation*}
\caption{Decomposition of fully supported staircase diagrams}\label{F:type_A_staircase_decomp}
\end{figure}

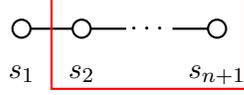
\begin{figure}[htbp]
  \begin{tikzpicture}[scale=.4]
    \draw[thick] (0,0)--(2,0)--(4.25,0)--(6.5,0);
    \draw[thick,fill=white] (0,0) circle (.3cm) node[label={[label distance=.1cm]-90:$s_1$}] { };
    \draw[thick,fill=white] (2,0) circle (.3cm) node[label={[label distance=.1cm]-90:$s_2$}]{ };
    \draw[white, fill=white] (4.25,0) circle (.75cm);
    \draw (4.25,0) node {$\cdots$};
    \draw[thick,fill=white] (6.5,0) circle (.3cm) node[label={[label distance=.1cm]-90:$s_{n+1}$}]{ };
    \draw[thick, red] (1,1)--(1,-2)--(7.5,-2)--(7.5,1)--(1,1);
  \end{tikzpicture}
\caption{Embedding of $\Gamma_n$ in $\Gamma_{n+1}$ }\label{F:graph_embedding_type_A}
\end{figure}

Note that second part of the decomposition in Figure \ref{F:type_A_staircase_decomp} may not be a valid staircase diagram which leads to the following definition.  First, we set $\Gamma_n\subseteq \Gamma_{n+1}$ as a subgraph by removing the vertex $s_1$ (See Figure \ref{F:graph_embedding_type_A}).

We say a $\mcD$ is an increasing (decreasing) \textbf{broken staircase diagram} over $\Gamma_n$ if we can write \[\mcD=\{B\cap [s_2,s_n] \ |\ B\in \mcD'\}\]
where $\mcD'$ is some increasing (decreasing) staircase diagram over $\Gamma_{n+1}$.  


\begin{lemma}\label{lem:staircase_catalan_broken}
Let $b_n$ denote the number of increasing (equivalently decreasing) broken staircase diagrams over $\Gamma_n$.  Then $b_n=c_{n+1}-c_n$ where $c_n$ denotes the $n$-th Catalan number.
\end{lemma}

\begin{proof}
Let $\mcD=\{B_1\prec B_2\prec\cdots\prec B_k\}$ be an increasing staircase diagram over $\Gamma_{n+1}$ and let
$$\mcB(\mcD):=\{B\cap[s_2,s_n]\ |\ B\in\mcD\}$$
denote the corresponding broken staircase diagram over $\Gamma_n$.  By Lemma \ref{lem:staircase_catalan}, the number of increasing staircase diagrams over $\Gamma_{n+1}$ is $c_{n+1}$.  We prove the lemma by determining the pre-images of the map $\mcB$.  First note that since $\mcD$ is increasing, we have that $s_1\in B_i$ if and only if $i=1$.  Hence $B_i\cap[s_2,s_n]=B_i$ unless $i=1$ and the pre-image of $B$ is determined by the changes on $B_1$.  Now if $B_1\cap[s_2,s_n]\subset B_2\cap [s_2,s_n]$, then $\mcB(\mcD)$ is uniquely determined by $\mcD$ as in Figure \ref{F:broken_case1}.

\begin{figure}[h]
    \begin{tikzpicture}[scale=0.4]
        \begin{scope}[shift={(-2,0)}]
        \ppAff{6}{
                {0,0,1,1,0},
                {0,2,2,2,0},
                {3,3,0,0,0}}
        \end{scope}
        \draw[->, shift={(0,1.5)},black,line width=.25mm](0,0)--(1,0);
        \begin{scope}[shift={(7,0)}]
        \ppAff{6}{
                {0,0,1,1,1},
                {0,2,2,2,0},
                {3,3,0,0,0}}
        \draw[dashed, red,line width=0.4mm] (-3,-1)--(-3,3);
        \end{scope}
        \end{tikzpicture}
        \caption{The broken diagram $\mcB(\mcD)$ determined by $\mcD$.}
        \label{F:broken_case1}
\end{figure}
Otherwise, $\mcB(\mcD)$ has two pre-images as in Figure \ref{F:broken_case2}.
\begin{figure}[h]
    \begin{tikzpicture}[scale=0.4]
        \begin{scope}[shift={(-2,0)}]
        \ppAff{6}{
                {0,0,1,1,0},
                {0,2,2,0,0},
                {3,3,0,0,0}}
        \end{scope}
        \draw[->, shift={(0,1.5)},black,line width=.25mm](0,0)--(1,0);
        \begin{scope}[shift={(7,0)}]
        \ppAff{6}{
                {0,0,1,1,1},
                {0,2,2,0,0},
                {3,3,0,0,0}}
        \draw[dashed, red,line width=0.4mm] (-3,-1)--(-3,3);
        \end{scope}
        \draw (10,1.5) node {\text{or}};
        \begin{scope}[shift={(16,0)}]
        \ppAff{6}{
                {0,0,0,0,4},
                {0,0,1,1,0},
                {0,2,2,0,0},
                {3,3,0,0,0}}
        \draw[dashed, red,line width=0.4mm] (-3,-1)--(-3,3);
        \end{scope}
        \end{tikzpicture}

        \caption{Two possibilities for $\mcD$ given $\mcB(\mcD)$.}
        \label{F:broken_case2}
\end{figure}
Broken staircase diagrams over $\Gamma_n$ with two pre-images under the map $\mcB$ can be identified with increasing staircase diagrams over $\Gamma_n$ via the second pre-image in Figure $\ref{F:broken_case2}$.  The Lemma now follows from Lemma \ref{lem:staircase_catalan}.
\end{proof}

\begin{proof}[Proof of Theorem \ref{thm:smooth_enumeration_typeA}]
We first note that the generating function for Catalan numbers is given by
\begin{equation}\label{eqn:Catalan}
\Cat(t):=\sum_{n=0}^{\infty} c_n\, t^n=\frac{1-\sqrt{1-4t}}{2t}
\end{equation}
and by Lemma \ref{lem:staircase_catalan}, the Catalan number $c_n$ denotes the number of increasing staircase diagrams over $\Gamma_n$.  Recall the $b_n$ denotes the number of increasing broken staircase diagrams over $\Gamma_n$ and let \[\Inc(t):=\sum_{n=0}^{\infty} b_n\, t^n.\]
Lemma \ref{lem:staircase_catalan_broken} implies that
\begin{equation}\label{eqn:Catalan_broken}
t+t\Inc(t)=\Cat(t)-t\Cat(t).
\end{equation}
Now suppose $\mcD$ is a fully supported staircase diagram over $\Gamma_n$.  Then either $\mcD$ is increasing on $\Gamma_n$ or $\mcD$ decomposes, as in Figure \ref{F:type_A_staircase_decomp}, into a smaller fully supported staircase diagram and a broken staircase diagram (note that this second case includes decreasing diagrams).
Let $\bar a_n$ denote number of fully supported staircase diagrams over $\Gamma_n$ and define
\[\overline A(t):=\sum_{n=0}^{\infty} \bar a_n\, t^n.\]  We now have
\begin{equation}\label{eqn:fully_supported}
    \overline A(t)=\Cat(t)+\overline{A}(t)\cdot \Inc(t).
\end{equation}
Finally, any staircase diagram is a disjoint union of fully supported staircase diagrams and hence
\begin{equation}\label{eqn:all_type_A}
A(t)=\frac{1+\overline{A}(t)}{1-t-t\overline{A}(t)}.
\end{equation}
The theorem follows from combining Equations \eqref{eqn:Catalan}, \eqref{eqn:Catalan_broken}, \eqref{eqn:fully_supported}, and \eqref{eqn:all_type_A}.
\end{proof}

One advantage of using staircase diagrams to enumerate smooth elements is that the techniques can extended to calculate generating functions for smooth and rationally smooth elements of other families of Coxeter groups.  Define the generating series
\[B(t):=\sum_{n=0} b_n\ t^n, \ C(t):=\sum_{n=0} c_n\ t^n,\
 D(t):=\sum_{n=3} d_n\ t^n,\ \text{and}\ BC(t):=\sum_{n=0} bc_n\ t^n\]
where $b_n,c_n,d_n$ denote number of smooth elements of type $B_n, C_n, D_n$ respectively and $bc_n$ denotes the number of rationally smooth elements of type $B_n/C_n$.   The following theorem is one of the main results of \cite{RS17}.

\begin{theorem}\label{thm:main_gen_series}
    Let $W(t) := \sum_n w_n\, t^n$ denote one of the above generating series,
    where $W = A$, $B$, $C$, $D$ or $BC$.  Then
    \begin{equation*}
        W(t)=\frac{P_W(t)+Q_W(t)\sqrt{1-4t}}{(1-t)^2(1-6t+8t^2-4t^3)}
    \end{equation*}
    where $P_W(t)$ and $Q_W(t)$ are polynomials given in Table \ref{TBL:polys}.
\end{theorem}

\begin{table}[h]
    \centering
    \begin{tabular}{ccc}
    \toprule
        Type  & $P_W(t)$  & $Q_W(t)$  \\
        \midrule
        $A$& $(1-4t)(1-t)^3$              & $t(1-t)^2$     \\
        $B$& $(1-5t+5t^2)(1-t)^3 $        & $(2t-t^2)(1-t)^3$     \\
        $C$& $1-7t+15t^2-11t^3-2t^4+5t^5$ & $t-t^2-t^3+3t^4-t^5 $   \\
        $D$& $(-4t+19t^2+8t^3-30t^4+16t^5)(1-t)^2$   & $(4t-15t^2+11t^3-2t^5)(1-t)$   \\
        $BC$& $1-8t+23t^2-29t^3+14t^4$    & $2t-6t^2+7t^3-2t^4$  \\
    \bottomrule
    \end{tabular}
    \caption{Polynomials in Theorem \ref{thm:main_gen_series}.}
    \label{TBL:polys}
\end{table}

The proof of Theorem \ref{thm:main_gen_series} involves enumerating staircase diagrams similar fashion to the proof of Theorem \ref{thm:smooth_enumeration_typeA}.  For type $D$, we can apply Corollary \ref{cor:smooth_staircase_bijection}.  Since types $B$ and $C$ are not simply-laced, we need to consider (rationally) smooth labellings of staircase diagrams that are not the maximal labelling.  These additional labellings are characterized by Theorem \ref{thm:BPgen} parts (1a) and (1b).


\section{BP decompositions and pattern avoidance}\label{S:BP_pattern_avoidance}

In this section we give an overview of how permutation pattern avoidance is related to BP decompositions.  Here we will only consider permutation groups (type $A$).  Recall that $\mfS_n$ is permutation group on $[n]=\{1,\ldots, n\}$.  The permutation group $\mfS_n$ is generated by the set of simple transpositions $S=\{s_1,\ldots, s_{n-1}\}$ where $s_i$ denotes the transposition swapping $i$ and $(i+1)$ and with the relations
\[s_i^2=(s_is_j)^2=(s_is_{i+1})^3=e\quad\text{for all $|i-j|>1$.} \]
Any $w\in \mfS_n$ has a unique expression in one-line notation $w=w(1)\cdots w(n)$.  We use matrices to represent permutations with nodes marking the points $(w(i),i)$ using the convention that $(1,1)$ marks the upper left corner.  For example, $w=3241$ corresponds to the matrix:

\[\begin{tikzpicture}[scale=0.4]
\draw[step=1.0,black] (0,0) grid (3,3);
\fill (0,1) circle (7pt);
\fill (1,2) circle (7pt);
\fill (2,0) circle (7pt);
\fill (3,3) circle (7pt);
\end{tikzpicture}\]

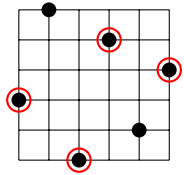
\begin{figure}[htbp]
\begin{tikzpicture}[scale=0.4]
\draw[step=1.0,black] (0,0) grid (5,5);
\fill (0,2) circle (7pt);\draw[thick, red] (0,2) circle (11pt);
\fill (1,5) circle (7pt);
\fill (2,0) circle (7pt);\draw[thick, red] (2,0) circle (11pt);
\fill (3,4) circle (7pt);\draw[thick, red] (3,4) circle (11pt);
\fill (4,1) circle (7pt);
\fill (5,3) circle (7pt);\draw[thick, red] (5,3) circle (11pt);
\end{tikzpicture}
\caption{$w=416253$ contains the pattern $3412$.}\label{fig:contains:3412}
\end{figure}

Let $u\in \mfS_k$ and $w\in \mfS_n$.  We say $w$ \textbf{contains the pattern} $u$ if there exists a subsequence $(i_1<\cdots<i_k)$ such that $w(i_1)\cdots w(i_k)$ has the same relative order as $u(1)\cdots u(k)$.  If no such sequence exists, we say that $w$ \textbf{avoids the pattern} $u$.  For example, in Figure \ref{fig:contains:3412}, we see that $w=416253$ contains the pattern $3412$, but avoids the pattern $1234$.

Permutation pattern avoidance has been used to characterize many geometric properties of Schubert varieties of type $A$.  A survey of these results can be found at \cite{AB16}.  Most notably, Lakshmibai and Sandhya prove that a Schubert variety $X(w)$ is smooth if and only if $w$ avoids the patterns 3412 and 4231 in \cite{LS90}.  Combining this result with Corollary \ref{cor:smooth_complete_BP}, we have the following theorem:

\begin{theorem}\label{thm:pattern_smooth_complete_BP}
  If the permutation $w$ avoids the patterns $3412$ and $4231$, then $w$ has a complete BP decomposition.

\[\begin{tikzpicture}[scale=0.4]
\draw[step=1.0,black] (0,0) grid (3,3);
\fill (0,1) circle (7pt);
\fill (1,0) circle (7pt);
\fill (2,3) circle (7pt);
\fill (3,2) circle (7pt);
\draw (1.5,-2) node {$3412$};
\end{tikzpicture}
\hspace{1in}
\begin{tikzpicture}[scale=0.4]
\draw[step=1.0,black] (0,0) grid (3,3);
\fill (0,0) circle (7pt);
\fill (1,2) circle (7pt);
\fill (2,1) circle (7pt);
\fill (3,3) circle (7pt);
\draw (1.5,-2) node {$4231$};
\end{tikzpicture}\]

\end{theorem}

The geometric version of Theorem \ref{thm:pattern_smooth_complete_BP} states that smooth Schubert varieties of type $A$ are iterated fiber bundles of Grassmannian varieties.  This geometric result was proved by Ryan in \cite{Ry87}.  Wolper gives an analogous result for Schubert varieties over algebraically closed fields in characteristic zero in \cite{Wo89}.   Note that it is not necessary for $w$ to avoid 3412 and 4231 for $w$ to have a complete BP decomposition.  In fact, if $w=\text{4231}$, then
\[w=(s_1s_3s_2)(s_3)(s_1)\]
is a complete BP decomposition of $w$.  The following theorem from \cite{AR18} is a precise pattern avoidance characterization of permutations that have complete BP decompositions.

\begin{theorem}\label{thm:pattern_complete_BP}
  The permutation $w$ avoids the patterns $3412$, $52341$, and $635241$ if and only if $w$ has a complete BP decomposition.
\[\vcenteredinclude{\begin{tikzpicture}[scale=0.4]
\draw[step=1.0,black] (0,0) grid (3,3);
\fill (0,1) circle (7pt);
\fill (1,0) circle (7pt);
\fill (2,3) circle (7pt);
\fill (3,2) circle (7pt);
\end{tikzpicture}}
\hspace{0.5in}
\vcenteredinclude{\begin{tikzpicture}[scale=0.4]
\draw[step=1.0,black] (0,0) grid (4,4);
\fill (0,0) circle (7pt);
\fill (1,3) circle (7pt);
\fill (2,2) circle (7pt);
\fill (3,1) circle (7pt);
\fill (4,4) circle (7pt);
\end{tikzpicture}}
\hspace{0.5in}
\vcenteredinclude{\begin{tikzpicture}[scale=0.4]
\draw[step=1.0,black] (0,0) grid (5,5);
\fill (0,0) circle (7pt);
\fill (1,3) circle (7pt);
\fill (2,1) circle (7pt);
\fill (3,4) circle (7pt);
\fill (4,2) circle (7pt);
\fill (5,5) circle (7pt);
\end{tikzpicture}}\]
\[3412\hspace{0.8in} 52341\hspace{0.875in} 635241\]
\end{theorem}

\subsection{Split pattern avoidance}\label{ss:split_patterns}

The proof of Theorem \ref{thm:pattern_complete_BP} relies on the idea of split pattern avoidance which is used to characterize Grassmannian BP decompositions of permutations with respect to $J=S\setminus\{s_r\}$ for any $s_r\in S$. 

A \newword{split pattern} $w=w_1|w_2\in \mfS_n$ is a divided permutation with \[w_1=w(1)\cdots w(j)\quad\text{and}\quad w_2=w(j+1)\cdots w(n)\] for some $1\leq j<n$.

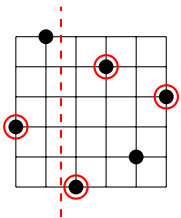
\begin{figure}[htbp]
\begin{tikzpicture}[scale=0.4]
\draw[step=1.0,black] (0,0) grid (5,5);
\fill (0,2) circle (7pt);\draw[thick, red] (0,2) circle (11pt);
\fill (1,5) circle (7pt);
\fill (2,0) circle (7pt);\draw[thick, red] (2,0) circle (11pt);
\fill (3,4) circle (7pt);\draw[thick, red] (3,4) circle (11pt);
\fill (4,1) circle (7pt);
\fill (5,3) circle (7pt);\draw[thick, red] (5,3) circle (11pt);
\draw[dashed,thick, red] (1.5,-1)--(1.5,6);
\end{tikzpicture}
\caption{The permutation $416253$ contains $3|412$ with respect to position $r=2$.}\label{fig:split_pattern_416253}
\end{figure}

\begin{defn}
Let $k\leq n$ and $r<n$.  Let $w\in \mfS_n$ and let \[u=u(1)\cdots u(j)|u(j+1)\cdots u(k)\] denote a split pattern.  We say $w$ \textbf{contains the split pattern $u$ with respect to position $r$} if there exists a sequence $(i_1<\cdots<i_k)$ such that

\begin{enumerate}
\item $w(i_1)\cdots w(i_k)$ has the same relative order as $u$.
\item $i_j\leq r<i_{j+1}.$
\end{enumerate}
Otherwise, we say the permutation $w$ \textbf{avoids the split pattern $u$ with respect to position $r$}.
\end{defn}
In other words, $w$ contains $u=u_1|u_2$ if it contains $u$ in the usual sense of pattern containment, but with the extra condition that $u_1$ appears to the right of the $r$-th position and $u_2$ to the left of the $r$-th position in the one-line notation of $w$.  For example $w=416253$ contains the split pattern $3|412$ with respect to positions $r=1,2$ but avoids $3|412$ with respect to $r=3,4,5$ (See Figure \ref{fig:split_pattern_416253}).

The next theorem is from \cite[Theorem 1.1]{AR18} and completely characterizes Grassmannian BP decompositions in terms of split pattern avoidance.

\begin{theorem}\label{thm:split_patterns_BP}
  Let $r<n$ and $w\in \mfS_n$.  Then $w$ has a Grassmannian BP decomposition with respect to $J=S\setminus\{s_r\}$ if and only if $w$ avoids the split patterns $3|12$ and $23|1$ with respect to position $r$.
$$\begin{tikzpicture}[scale=0.4]
\draw[step=1.0,black] (0,0) grid (2,2);
\fill (0,0) circle (7pt);
\fill (1,2) circle (7pt);
\fill (2,1) circle (7pt);
\draw[dashed,thick, red] (0.5,-1)--(0.5,3);
\draw (1,-2) node {$3|12$};\end{tikzpicture}\hspace{1in}
\begin{tikzpicture}[scale=0.4]
\draw[step=1.0,black] (0,0) grid (2,2);
\fill (0,1) circle (7pt);
\fill (1,0) circle (7pt);
\fill (2,2) circle (7pt);
\draw[dashed,thick, red] (1.5,-1)--(1.5,3);
\draw (1,-2) node {$23|1$};\end{tikzpicture}
$$
\end{theorem}

\begin{proof}[Outline of Proof of Theorem \ref{thm:split_patterns_BP}]
We use Theorem \ref{thm:BP_characterization} part (4) which states that a parabolic decomposition $w=vu$ with respect to $J$ is a BP decomposition if and only if $S(v)\cap J\subseteq D_L(u)$.  The next lemma gives an explicit description of these ideas in terms of the one-line notation of permutations.  We leave the proof as an exercise.  

\begin{lemma}\label{lem:perm_parabolic}
Let $w\in \mfS_n$ and $r<n$ and write
\[w=w_1|w_2=w(1)\cdots w(r)|w(r+1)\cdots w(n).\]

Let $w=vu$ denote the parabolic decomposition with respect to $J=S\setminus\{s_r\}$.  Then the following are true:

\begin{enumerate}
\item $v=v_1|v_2$ where $v_1$ and $v_2$ respectively consist the entries in $w_1$ and $w_2$ arranged in increasing order and
\[S(v)=\{s_k\in S \ |\ v(r+1)\leq k<v(r)\}.\]
\item $u=u_1|u_2$ where $u_1$ and $u_2$ are respectively the unique permutations on $\{1,\ldots,r\}$ and $\{r+1,\ldots,n\}$ with relative orders of $w_1$ and $w_2$ and \[D_L(u)=\{s_k\in S \ |\ u^{-1}(k+1)<u^{-1}(k)\}.\]
\end{enumerate}
\end{lemma}

The description of the decent set in part (2) of Lemma \ref{lem:perm_parabolic} is equivalent to saying that $s_k$ is a left descent of $u$ if and only if the node in the $k$-th row is to the right of the node in the $(k+1)$-th row in the permutation matrix of $u$.  The proof of Theorem \ref{thm:split_patterns_BP} follows from showing that avoiding the split patterns $3|12$ and $23|1$ with respect to position $r$ is equivalent to $S(v)\setminus\{s_r\}\subseteq D_L(u)$ using Lemma \ref{lem:perm_parabolic}.  We illustrate this connection with the following examples:
\begin{example}
Let $w=17264|5938$ and note that $w$ avoids $3|12$ and $23|1$ with respect to position $r=5$.  If $w=vu$ is the parabolic decomposition with respect to $J=S\setminus\{s_5\}$, then $v=12467|3589$ and $u=15243|7968$ as seen in Figure \ref{fig:splitBPex1}.  Lemma \ref{lem:perm_parabolic} says that 
\[S(v)\setminus\{s_5\}=\{s_3,s_4,s_6\}\quad\text{and}\quad D_L(u)=\{s_3,s_4,s_6,s_7\}\]
and hence $S(v)\setminus\{s_5\}\subseteq D_L(u)$.

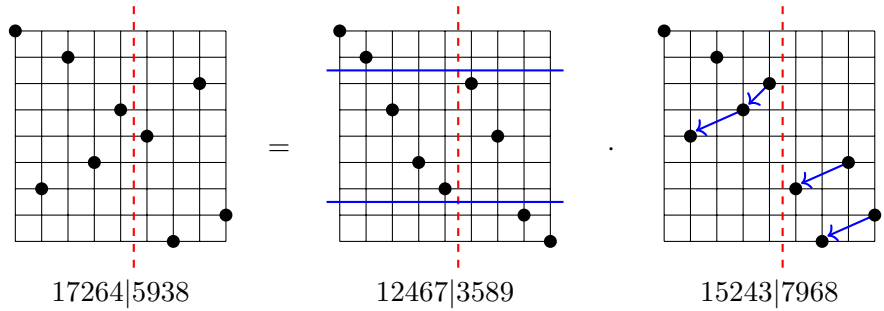
\begin{figure}[htbp]
\begin{tikzpicture}[scale=0.35]
\draw[step=1.0,black] (0,0) grid (8,8);
\fill (0,8) circle (7pt);
\fill (1,2) circle (7pt);
\fill (2,7) circle (7pt);
\fill (3,3) circle (7pt);
\fill (4,5) circle (7pt);
\fill (5,4) circle (7pt);
\fill (6,0) circle (7pt);
\fill (7,6) circle (7pt);
\fill (8,1) circle (7pt);
\draw[dashed,thick, red] (4.5,-1)--(4.5,9);
\draw (4,-2) node {$17264|5938$};
\end{tikzpicture}
\begin{tikzpicture}[scale=0.35]
\draw (0,3.5) node {$=$};
\draw (-1,0) node {};\draw (1,0) node {};
\end{tikzpicture}
\begin{tikzpicture}[scale=0.35]
\draw[step=1.0,black] (0,0) grid (8,8);
\fill (0,8) circle (7pt);
\fill (1,7) circle (7pt);
\fill (2,5) circle (7pt);
\fill (3,3) circle (7pt);
\fill (4,2) circle (7pt);
\fill (5,6) circle (7pt);
\fill (6,4) circle (7pt);
\fill (7,1) circle (7pt);
\fill (8,0) circle (7pt);
\draw[dashed,thick, red] (4.5,-1)--(4.5,9);
\draw[thick, blue] (-0.5,6.5)--(8.5,6.5);
\draw[thick, blue] (-0.5,1.5)--(8.5,1.5);
\draw (4,-2) node {$12467|3589$};
\end{tikzpicture}
\begin{tikzpicture}[scale=0.35]
\draw (0,3.5) node {$\cdot$};
\draw (-1,0) node {};\draw (1,0) node {};
\end{tikzpicture}
\begin{tikzpicture}[scale=0.35]
\draw[step=1.0,black] (0,0) grid (8,8);
\draw[thick, ->, blue] (4,6)--(3.2,5.2);
\draw[thick, ->, blue] (3,5)--(1.2,4.2);
\draw[thick, ->, blue] (7,3)--(5.2,2.2);
\draw[thick, ->, blue] (8,1)--(6.2,0.2);
\fill (0,8) circle (7pt);
\fill (1,4) circle (7pt);
\fill (2,7) circle (7pt);
\fill (3,5) circle (7pt);
\fill (4,6) circle (7pt);
\fill (5,2) circle (7pt);
\fill (6,0) circle (7pt);
\fill (7,3) circle (7pt);
\fill (8,1) circle (7pt);
\draw[dashed,thick, red] (4.5,-1)--(4.5,9);
\draw (4,-2) node {$15243|7968$};
\end{tikzpicture}
\caption{The parabolic decomposition $w=vu$ with respect to $J=S\setminus\{s_5\}$.  The rows between the blue lines corresponds to $S(v)$ and the blue arrows denote $D_L(u)$.}\label{fig:splitBPex1}
\end{figure}

If we take the parabolic decomposition of $w=1726|45938$ with respect to $J=S\setminus\{s_4\}$, then $w$ contains $3|12$ with respect to $r=4$. In this case $v=1267|34589$ and $u=1423|57968$ (See Figure \ref{fig:splitBPex2}).  Lemma \ref{lem:perm_parabolic} says that 
\[S(v)\setminus\{s_4\}=\{s_3,s_5,s_6\}\quad\text{and}\quad D_L(u)=\{s_3,s_5,s_7\}\]
and hence $S(v)\setminus\{s_5\}\nsubseteq D_L(u)$.

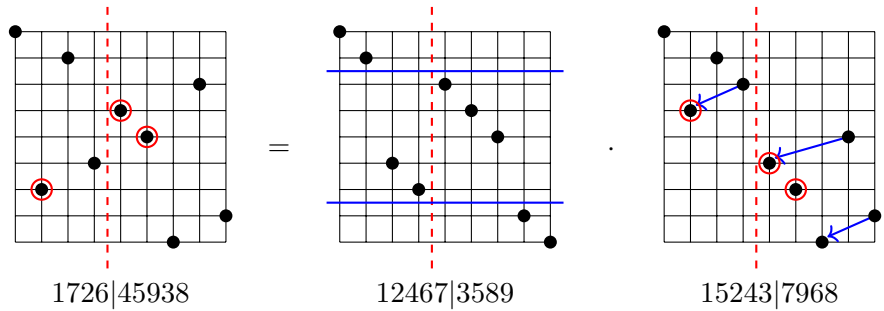
\begin{figure}[htbp]
\begin{tikzpicture}[scale=0.35]
\draw[step=1.0,black] (0,0) grid (8,8);
\fill (0,8) circle (7pt);
\fill (1,2) circle (7pt);\draw[thick, red] (1,2) circle (11pt);
\fill (2,7) circle (7pt);
\fill (3,3) circle (7pt);
\fill (4,5) circle (7pt);\draw[thick, red] (4,5) circle (11pt);
\fill (5,4) circle (7pt);\draw[thick, red] (5,4) circle (11pt);
\fill (6,0) circle (7pt);
\fill (7,6) circle (7pt);
\fill (8,1) circle (7pt);
\draw[dashed,thick, red] (3.5,-1)--(3.5,9);
\draw (4,-2) node {$1726|45938$};
\end{tikzpicture}
\begin{tikzpicture}[scale=0.35]
\draw (0,3.5) node {$=$};
\draw (-1,0) node {};\draw (1,0) node {};
\end{tikzpicture}
\begin{tikzpicture}[scale=0.35]
\draw[step=1.0,black] (0,0) grid (8,8);
\fill (0,8) circle (7pt);
\fill (1,7) circle (7pt);
\fill (2,3) circle (7pt);
\fill (3,2) circle (7pt);
\fill (4,6) circle (7pt);
\fill (5,5) circle (7pt);
\fill (6,4) circle (7pt);
\fill (7,1) circle (7pt);
\fill (8,0) circle (7pt);
\draw[dashed,thick, red] (3.5,-1)--(3.5,9);
\draw[thick, blue] (-0.5,6.5)--(8.5,6.5);
\draw[thick, blue] (-0.5,1.5)--(8.5,1.5);
\draw (4,-2) node {$12467|3589$};
\end{tikzpicture}
\begin{tikzpicture}[scale=0.35]
\draw (0,3.5) node {$\cdot$};
\draw (-1,0) node {};\draw (1,0) node {};
\end{tikzpicture}
\begin{tikzpicture}[scale=0.35]
\draw[step=1.0,black] (0,0) grid (8,8);
\draw[thick, ->, blue] (3,6)--(1.2,5.2);
\draw[thick, ->, blue] (7,4)--(4.2,3.2);
\draw[thick, ->, blue] (8,1)--(6.2,0.2);
\fill (0,8) circle (7pt);
\fill (1,5) circle (7pt);\draw[thick, red] (1,5) circle (11pt);
\fill (2,7) circle (7pt);
\fill (3,6) circle (7pt);
\fill (4,3) circle (7pt);\draw[thick, red] (4,3) circle (11pt);
\fill (5,2) circle (7pt);\draw[thick, red] (5,2) circle (11pt);
\fill (6,0) circle (7pt);
\fill (7,4) circle (7pt);
\fill (8,1) circle (7pt);
\draw[dashed,thick, red] (3.5,-1)--(3.5,9);
\draw (4,-2) node {$15243|7968$};
\end{tikzpicture}
\caption{The parabolic decomposition $w=vu$ with respect to $J=S\setminus\{s_4\}$.}\label{fig:splitBPex2}
\end{figure}
\end{example}
\end{proof}

We remark the an explicit formula for the number of permutations $w\in \mfS_n$ which avoid $3|12$ and $23|1$ with respect to a given position $r$ is calculated by Grigsby and Richmond in \cite{GR24}.   The connection between Theorem \ref{thm:pattern_complete_BP} and \ref{thm:split_patterns_BP} is the following proposition.

\begin{prop}\label{prop:split_pattern_complete}
If $w\in \mfS_n$ avoids the patterns $3412$, $52341$, and $635241$, then there exists $s_r\in S(w)$ where $w$ avoids the split patterns $3|12$ and $23|1$ with respect to position $r$.
\end{prop}

One can prove Proposition \ref{prop:split_pattern_complete} by contradiction and we refer the reader to \cite{AR18} for more details.

\begin{proof}[Outline of the proof of Theorem \ref{thm:pattern_complete_BP}]
Proposition \ref{prop:split_pattern_complete} implies that if $w\in \mfS_n$ avoids the patterns $3412$, $52341$, and $635241$, then $w$ has a Grassmannian BP decomposition $w=vu$ with respect to some $r<n$. It can be shown that $u$ also avoid these patterns and hence we can iterate this process yielding a complete BP decomposition of $w$.  For more details, see \cite{AR18}.
\end{proof}

If $w$ avoids $3412$, $52341$, and $635241$, then we can we construct complete BP decompositions of $w$ by finding positions such that $w$ avoids $3|12$ and $23|1$ and then iterating the process on factor $u$.

\begin{example}
Observe that $w=513462$ avoids $3412$, $52341$, and $635241$ and hence $w$ has a complete BP decomposition.  Complete BP decompositions of $w$ correspond to sequences of ``splittings" along lines that avoid $3|12$ and $23|1$.  For example, we can split $w$ along the sequence of positions $(3,2,4,1,5)$ as in Figure \ref{fig:split_sequence}.

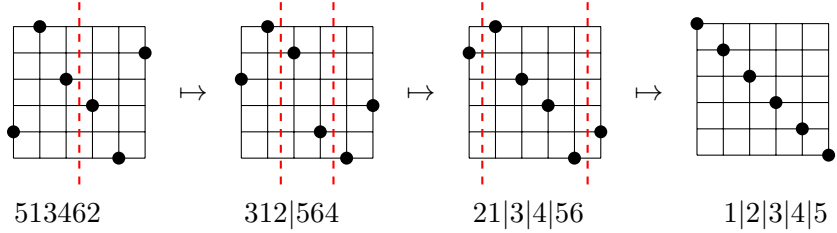
\begin{figure}
$$
\vcenteredinclude{\begin{tikzpicture}[scale=0.35]
\draw[step=1.0,black] (0,0) grid (5,5);
\fill (0,1) circle (7pt);
\fill (1,5) circle (7pt);
\fill (2,3) circle (7pt);
\fill (3,2) circle (7pt);
\fill (4,0) circle (7pt);
\fill (5,4) circle (7pt);
\draw[dashed,thick, red] (2.5,-1)--(2.5,6);
\end{tikzpicture}}\hspace{0.1in}\mapsto\hspace{0.1in}
\vcenteredinclude{\begin{tikzpicture}[scale=0.35]
\draw[step=1.0,black] (0,0) grid (5,5);
\fill (0,3) circle (7pt);
\fill (1,5) circle (7pt);
\fill (2,4) circle (7pt);
\fill (3,1) circle (7pt);
\fill (4,0) circle (7pt);
\fill (5,2) circle (7pt);
\draw[dashed,thick, red] (1.5,-1)--(1.5,6);
\draw[dashed,thick, red] (3.5,-1)--(3.5,6);
\end{tikzpicture}}\hspace{0.1in}\mapsto\hspace{0.1in}
\vcenteredinclude{\begin{tikzpicture}[scale=0.35]
\draw[step=1.0,black] (0,0) grid (5,5);
\fill (0,4) circle (7pt);
\fill (1,5) circle (7pt);
\fill (2,3) circle (7pt);
\fill (3,2) circle (7pt);
\fill (4,0) circle (7pt);
\fill (5,1) circle (7pt);
\draw[dashed,thick, red] (0.5,-1)--(0.5,6);
\draw[dashed,thick, red] (4.5,-1)--(4.5,6);
\end{tikzpicture}}\hspace{0.1in}\mapsto\hspace{0.1in}
\vcenteredinclude{\begin{tikzpicture}[scale=0.35]
\draw[step=1.0,black] (0,0) grid (5,5);
\fill (0,5) circle (7pt);
\fill (1,4) circle (7pt);
\fill (2,3) circle (7pt);
\fill (3,2) circle (7pt);
\fill (4,1) circle (7pt);
\fill (5,0) circle (7pt);
\end{tikzpicture}}
$$
\[513462\hspace{0.75in} 312|564\hspace{0.7in} 21|3|4|56\hspace{0.73in} 1|2|3|4|5\]
    \caption{Decomposing $w=513462$ along splits that avoid $3|12$ and $23|1$.}
    \label{fig:split_sequence}
\end{figure}
\end{example}

Next we use Theorem \ref{thm:split_patterns_BP} to give a new proof of Gasharov's Theorem that Poincar\'e polynomials of smooth permutations are products of $q$-integers.  The following proposition is the ``forward" direction of Theorem \ref{thm:smooth_perm}.

\begin{prop}\label{prop:smooth_perm2}
Let $w\in \mfS_n$.  If $w$ avoids $3412$ and $4231$, then either $w$ or $w^{-1}$ has a BP decomposition $vu$ with respect to $J=S\setminus\{s_{n-1}\}$ where
\[P_w(q)=[\ell(v)+1]_q\cdot P_u(q)\] and $u\in W_J\simeq \mfS_{n-1}$ also avoids $3412$ and $4231$.
\end{prop}

\begin{proof}
We prove the proposition by contradiction.  Let $w\in \mfS_n$ and assume $w$ avoids $3412$ and $4231$.  For the sake of contradiction, suppose that both $w$ and $w^{-1}$ do not have BP decompositions with respect to $J=S\setminus\{s_{n-1}\}$.  Theorem \ref{thm:split_patterns_BP} implies that both $w$ and $w^{-1}$ contain the split pattern $23|1$ with respect to position $r=n-1$.
\begin{figure}[htbp]
\begin{tikzpicture}[scale=0.4]
\draw[step=1.0,black] (0,0) grid (2,2);
\fill (0,0) circle (7pt);
\fill (1,2) circle (7pt);
\fill (2,1) circle (7pt);
\draw[dashed,thick, red] (-1,0.5)--(3,0.5);
\end{tikzpicture}
\hspace{1in}
\begin{tikzpicture}[scale=0.4]
\draw[step=1.0,black] (0,0) grid (2,2);
\fill (0,1) circle (7pt);
\fill (1,0) circle (7pt);
\fill (2,2) circle (7pt);
\draw[dashed,thick, red] (1.5,-1)--(1.5,3);
\end{tikzpicture}
\caption{The horizontal and vertical split pattern $23|1$.}\label{fig:ver_hor_split_231}
\end{figure}
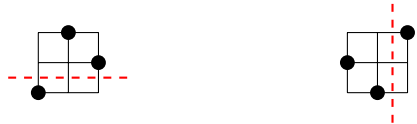
Since $w^{-1}$ corresponds to the transpose of $w$, we consider a ``horizontal" analogue of split pattern containment in Figure \ref{fig:ver_hor_split_231}.  Note that it is not possible for either $w$ or $w^{-1}$ to contain the other split pattern, $3|12$, with respect to position $r=n-1$.  Let $w(d)=n$ and $w(n)=e$ and consider the matrix diagram of $w$ where we mark the nodes $(d,n)$ and $(n,e)$ as in Figure \ref{fig:ver_hor_split_231part2}.  These nodes divide the matrix into four regions of which we label three of them $A,B,$ and $C$.
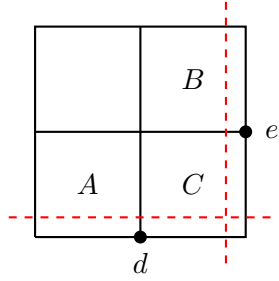
\begin{figure}[htbp]
$$
\begin{tikzpicture}[scale=0.35]
\draw[thick] (0,4) -- (8,4);
\draw[thick] (4,0) -- (4,8);
\fill (4,0) circle (7pt);
\fill (8,4) circle (7pt);
\draw[thick] (0,0) -- (0,8) -- (8,8) -- (8,0) -- (0,0);
\draw (2,2) node {$A$};\draw (6,6) node {$B$}; \draw (6,2) node {$C$};
\draw (4,-1) node {$d$};\draw (9,4) node {$e$};
\draw[dashed,thick, red] (7.25,-1)--(7.25,9);
\draw[dashed,thick, red] (-1,0.75)--(9,0.75);
\end{tikzpicture}
$$
\caption{The matrix of $w$ with $w(d)=n$ and $w(n)=e$ and regions $A,B,C$.}\label{fig:ver_hor_split_231part2}
\end{figure}

We have two cases to consider when containing the split pattern $23|1$ both vertically and horizontally with respect to position $r=n-1$.  First, if either region $A$ or $B$ contain no nodes, then region $C$ must contain 2 increasing nodes which implies that $w$ contains 4231.  Otherwise, each of regions $A$ and $B$ must contain at least one node which implies $w$ contains 3412 (See Figure \ref{fig:ver_hor_split_231part3}).

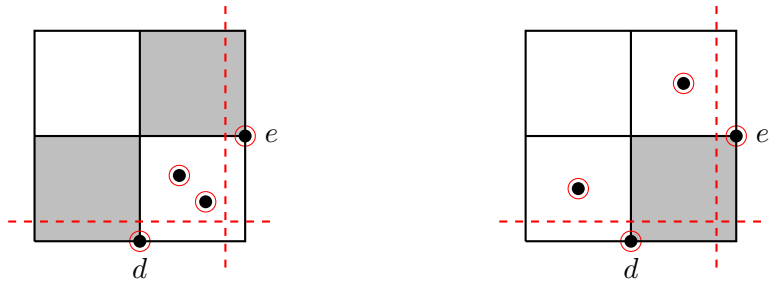
\begin{figure}[htbp]
\begin{tikzpicture}[scale=0.35]
\fill[lightgray] (0,0) rectangle (4,4);
\fill[lightgray] (4,4) rectangle (8,8);
\draw[thick] (0,4) -- (8,4);
\draw[thick] (4,0) -- (4,8);
\fill (4,0) circle (7pt);\draw[red] (4,0) circle (11pt);
\fill (8,4) circle (7pt);\draw[red] (8,4) circle (11pt);
\fill (5.5,2.5) circle (7pt);\draw[red] (5.5,2.5) circle (11pt);
\fill (6.5,1.5) circle (7pt);\draw[red] (6.5,1.5) circle (11pt);
\draw[thick] (0,0) -- (0,8) -- (8,8) -- (8,0) -- (0,0);
\draw (4,-1) node {$d$};\draw (9,4) node {$e$};
\draw[dashed,thick, red] (7.25,-1)--(7.25,9);
\draw[dashed,thick, red] (-1,0.75)--(9,0.75);
\end{tikzpicture}
\hspace{1in}
\begin{tikzpicture}[scale=0.35]
\fill[lightgray] (4,0) rectangle (8,4);
\draw[thick] (0,4) -- (8,4);
\draw[thick] (4,0) -- (4,8);
\fill (4,0) circle (7pt);\draw[red] (4,0) circle (11pt);
\fill (8,4) circle (7pt);\draw[red] (8,4) circle (11pt);
\fill (2,2) circle (7pt);\draw[red] (2,2) circle (11pt);
\fill (6,6) circle (7pt);\draw[red] (6,6) circle (11pt);
\draw[thick] (0,0) -- (0,8) -- (8,8) -- (8,0) -- (0,0);
\draw (4,-1) node {$d$};\draw (9,4) node {$e$};
\draw[dashed,thick, red] (7.25,-1)--(7.25,9);
\draw[dashed,thick, red] (-1,0.75)--(9,0.75);
\end{tikzpicture}
\caption{The matrix of $w$ containing $4231$ and $3412$.}\label{fig:ver_hor_split_231part3}
\end{figure}

In either case, we have a contradiction and hence at least one of $w$ or $w^{-1}$ has a BP decomposition $vu$ with respect to $J$.  The fact that $u$ is smooth follows from Lemma \ref{lem:perm_parabolic} part (2).  Since $s_{n-1}$ is a leaf in the Coxeter diagram of $\mfS_n$, the interval $[e,v]^J$ is chain and hence $P_v^J(v)=[\ell(v)+1]_q$.  This completes the proof.
\end{proof}

\subsection{Related results on pattern avoidance}
In this section by state two analogues of the following theorem which summarizes various characterizations of smooth permutations.

\begin{theorem}\label{thm:smooth2}
Let $w\in \mfS_n$.  Then the following are equivalent.
\begin{enumerate}
\item $w$ avoids $3412$ and $4231$.
\item $X(w)$ is an iterated fiber-bundle of Grassmannian varieties.
\item The interval $[e,w]$ is rank symmetric.
\end{enumerate}
\end{theorem}

Theorem \ref{thm:smooth2} follows from the combined works of Lakshmibai-Sandhya \cite{LS90}, Ryan \cite{Ry87}, and Carrell \cite{Ca94}.   Note that Theorem \ref{thm:pattern_complete_BP} can be viewed as analogue of the equivalence of parts (1) and part (2) in Theorem \ref{thm:smooth2} where we replace part (2) with an iterated fiber-bundle of Grassmannian Schubert varieties.  Each Grassmannian has a co-dimension one Schubert variety which is unique in the sense that, as a Weil divisor, it generates the Picard group of the Grassmannian.   We call this variety a Grassmannian Schubert divisor.  The following theorem is another analogue of the equivalence of parts (1) and (2) and is proved by Azam in \cite{Az23}.

\begin{theorem}\label{thm:divisor}
Let $w\in \mfS_n$.  Then the following are equivalent.
\begin{enumerate}
    \item $w$ avoids $3412, 52341, 52431,$ and $53241$.
    \item $X(w)$ is an iterated fiber bundle of Grassmannian varieties or Grassmannian Schubert divisors.
\end{enumerate}
\end{theorem}

The class of permutations in Theorem \ref{thm:divisor} is larger than the smooth class of permutations, but within the class of permutations that have complete BP decompositions.  Note that Grassmannian Schubert divisors are almost always singular varieties.  One consequence of Theorem \ref{thm:divisor} is that the generating function for permutations that avoid $3412, 52341, 52431,$ and $53241$ can be calculated using labelled staircase diagrams.  This calculation uses ``Catalan type" objects similar to those used to prove Theorem \ref{thm:smooth_enumeration_typeA}.  For more details see \cite{Az23}.

The next theorem is a analogue of the equivalence of parts (1) and (3) in Theorem \ref{thm:smooth2}.  Given a poset $P$, the dual poset $P^*$ is obtained by reversing the partial order.  We say $P$ is \newword{self-dual} if $P\simeq P^*$ as posets.  It is easy to check that any graded self-dual poset is rank symmetric.  However the converse may not be true. The next theorem is proved by Gaetz and Gao in \cite{GG20}.

\begin{theorem}\label{thm:polished}
Let $w\in \mfS_n$.  Then the following are equivalent.
\begin{enumerate}
\item $w$ avoids $3412, 4231, 34521, 45321, 54123,$ and $54312$.  
\item The interval $[e,w]$ is self-dual (as a poset).
\end{enumerate}
\end{theorem}

The authors refer to permutations characterized in Theorem \ref{thm:polished} as ``polished" permutations since the condition of self-duality on the interval $[e,w]$ is sufficient for smoothness, but not necessary.

\subsection{Affine permutations}
In this section we discuss applications of BP decompositions to the group of affine permutations denoted $\widetilde\mfS_n$.  An \newword{affine permutation} is a bijection $w:\ZZ\rightarrow \ZZ$ such that 
\begin{enumerate}
    \item $w(i+n)=w(i)+n$ for all $i\in \ZZ$ and 
    \item $\displaystyle \sum_{i=1}^n w(i)=\frac{n(n+1)}{2}$.
\end{enumerate}
Note that a regular permutation extends to an affine permutation by applying part (1) above to the one-line notation sequence $w(1)\cdots w(n)$. 
 Similarly, any affine permutation is uniquely determined by the ``window" of values 
\[\cdots w(-1), w(0),[w(1),w(2),\cdots, w(n)],w(n+1),w(n+2),\cdots\]
by the same extension.  For example $[4,2,3,1]$, $[8,1,-2,3]$, and $[-7, 7,6,4]$ are all examples of affine permutations in $\widetilde\mfS_4$.

The group of affine permutations is an infinite Coxeter group with generating set $S=\{s_0,s_1,\ldots,s_{n-1}\}$ and Coxeter graph:
\begin{equation*}
  \begin{tikzpicture}[scale=.5]
    \draw[thick] (0,0)--(2,0)--(4.25,0)--(6.5,0)--(3.25,1.5)--(0,0);
    \draw[thick,fill=white] (0,0) circle (.3cm) node[label={[label distance=.1cm]-90:$s_1$}] { };
    \draw[thick,fill=white] (2,0) circle (.3cm) node[label={[label distance=.1cm]-90:$s_2$}]{ };
    \draw[white, fill=white] (4.25,0) circle (.65cm);
    \draw (4.25,0) node {$\cdots$};
    \draw[thick,fill=white] (6.5,0) circle (.3cm) node[label={[label distance=.1cm]-90:$s_{n-1}$}]{ };
    \draw[thick,fill=white] (3.25,1.5) circle (.3cm) node[label={[label distance=.1cm]90:$s_0$}]{ };
  \end{tikzpicture}
\end{equation*}
Affine permutations are referred to as Coxeter groups of affine type $A$.  Note that all maximal parabolic subgroups of $\widetilde\mfS_n$ are isomorphic to the finite permutation group $\mfS_n$.

As with finite permutations, (rational) smoothness is closely tied to pattern avoidance and was studied by Billey and Crites in \cite{BC12}.  We say an affine permutation $w$ \newword{contains the (finite) pattern} $u\in\mfS_k$ if there is a sequence $(i_1<\cdots <i_k)$ such that $w(i_1)\cdots w(i_k)$ has the same relative order as $u$.   Note that the sequence  $(i_1<\cdots <i_k)$ does not necessarily have to be contained in the integers $[n]=\{1,\ldots,n\}$.  If $w$ does not contain $u$, we say it \newword{avoids the pattern} $u$.  The following thoerem if proved by Billey and Crites in \cite{BC12}

\begin{theorem}\label{thm:affine_perms1}
    The affine permutation $w\in \widetilde\mfS_n$ is rationally smooth if and only if one of the following hold:
    \begin{enumerate}
        \item $w$ avoids the patterns $3412$ and $4231$ or 
        \item $w$ is a twisted spiral permutation (see \cite[Section 2.5]{BC12}).
    \end{enumerate}
\end{theorem}
It is shown by Mitchell in \cite{BM10}, that if $w$ is a twisted spiral permutation, then the Schubert variety $X(w)$ is not smooth.  Hence smooth is not equivalent to rationally smooth for affine permutations.  One technical result used to prove Theorem \ref{thm:affine_perms1} is the following analogue of Theorem \ref{thm:BPgen}.
\begin{proposition}\label{prop:affine_BP1}
    If $w \in\widetilde\mfS_n$ avoids $3412$ and $4231$, then either $w$ or $w^{-1}$ has a Grassmannian BP decomposition $vu$ where both $v$ and $u$ belong to proper parabolic subgroups of $\widetilde\mfS_n$.
\end{proposition}
Billey and Crites show that, for the BP decomposition $vu$ found in Proposition \ref{prop:affine_BP1}, the Poincar\'e polynomial $P_v^J(q)$ is a $q$-binomial and hence palindromic.  They also show that $u$ is a smooth (finite) permutation.  So Theorem \ref{thm:BPexists} implies $P_w(q)$ is palindromic.  
The next theorem was partially conjectured in \cite{BC12} and proved by Richmond and Slofstra in \cite{RS18}.
\begin{theorem}\label{thm:affine_perms2}
  Let $w\in \widetilde\mfS_n$.  Then the following are equivalent:
    \begin{enumerate}
        \item $X(w)$ is smooth.
        \item $w$ avoids the patterns $3412$ and $4231$.
        \item $w$ has a Grassmannian BP decomposition $vu$ with respect to some $J$ where both $v$ and $u$ belong to proper parabolic subgroups of $\widetilde\mfS_n$.  
        Furthermore $v$ is the maximal element of $W_{S(v)}\cap W^J$ and $u$ is a smooth permutation in the $W_{S(u)}$.
    \end{enumerate}  
\end{theorem}

The proof of Theorem \ref{thm:affine_perms2} is similar to the proof of Theorem \ref{thm:smooth_Grass_BP}.  We remark that part (3) of Theorem \ref{thm:affine_perms2} implies that for an affine permutation $w\in\widetilde\mfS_n$, the affine Schubert variety $X(w)$ is smooth if and only if it is an iterated fiber bundle of Grassmannian varieties. 

Theorem \ref{thm:affine_perms2} also implies an analogue of Corollary \ref{cor:smooth_staircase_bijection} on staircase diagrams of affine type $A$ which we state below.  Since $\widetilde\mfS_n$ is an infinite Coxeter groups, we say a staircase diagram is \newword{spherical} if for each $B\in \mcD$, the parabolic subgroup $W_B$ is a finite Coxeter group.  The next statement is from \cite[Thoerem 3.3]{RS18}.
\begin{corollary}\label{cor:affine_SD_bijection}
    The maximal labelling gives a bijection between spherical staircase diagrams over the Coxeter graph of $\widetilde\mfS_n$ and smooth affine permutations in $\widetilde\mfS_n$.
\end{corollary}
One immediate consequence of Corollary \ref{cor:affine_SD_bijection} is that the number of affine permutations that avoid $3412$ and $4231$ in $\widetilde\mfS_n$ is finite.  This fact also follows directly from results in \cite{BC12}.  For single patterns, Crites proved in \cite{Cr10} that the number affine permutations in $\widetilde\mfS_n$ avoiding a pattern $u$ is finite if and only if $u$ contains $321$.  

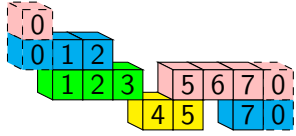
\begin{figure}[htbp]
\begin{tikzpicture}[scale=0.4]
    \ppAfff{9}{
        {1,1,0,3,3,0,0,0,0},
        {4,4,4,4,0,2,2,2,0},
        {0,0,0,0,0,0,1,1,1},
        {0,0,0,0,0,0,0,0,4}}
\end{tikzpicture}
\caption{Staircase diagram of affine type $A$}\label{fig:affine_staircase}
\end{figure}

Staircase diagrams of affine type $A$ can be thought of as staircase diagrams of finite type $A$ that ``loop" back on themselves since the Coxeter graph is a cycle.  Figure \ref{fig:affine_staircase} gives an example of an affine staircase diagram.  For more details see \cite{RS18}.
As with Theorems \ref{thm:smooth_enumeration_typeA} and \ref{thm:main_gen_series} we can use staircase diagrams to enumerate smooth affine permutations.  The following is proved in \cite[Theorem 1.1]{RS18}.
\begin{theorem}
    Let $\widetilde A(t):=\sum \tilde a_n t^n$ where $\tilde a_n$ denotes the number of smooth affine permutations in $\widetilde\mfS_n$.  Then 
    \begin{equation*}
        \widetilde A(t) = \frac{P(t) - Q(t) \sqrt{1-4t}}{(1-t)(1-4t)\left(1-6t+8t^2-4t^3\right)} 
    \end{equation*}
    where
    \begin{equation*}
        P(t) = (1-4t)\left(2-11t+18t^2-16t^3+10t^4-4t^5\right) 
    \end{equation*}
    and
    \begin{equation*}
        Q(t) = (1-t)(2-t)\left(1-6t+6t^2\right). 
    \end{equation*}
\end{theorem}

\section{Future directions}\label{S:future}

We state some open questions and possible future directions for the study of Coxeter groups in relation to BP-decompositions.

\begin{question}
    While rational smoothness for Coxeter groups of finite Lie type have been extensively studied.  Characterizations of rationally smooth elements for arbitrary Coxeter groups are relatively unknown.  For example, if $w$ is rationally smooth, does $w$ have a Grassmannian BP decomposition?  Does Theorem \ref{thm:hypmain} hold for inversion hyperplane arrangements of rationally smooth elements in arbitrary Coxeter groups?  We remark that Richmond and Slofstra study rationally smooth elements in Coxeter groups that avoid certain rank 3 parabolic subgroups in \cite{RS12}.
\end{question}

\begin{question}
Let $W$ be a Coxeter group and for $u\leq v\in W$, define the Poincar\'e polynomial of the interval
\[P_{u,w}(q):=\sum_{z\in[u,w]} q^{\ell(z)-\ell(u)}.\]
If $u=e$, then this is the usual Poincar\'e polynomial $P_w(q)$.  For example, if $u=s_2$ and $v=s_2s_1s_3s_2$, then Figure \ref{fig:Bruhat_interval_uv} shows the interval $[u,v]$ and  
\[P_{u,v}(q)=1+4q+4q^2+q^3.\]
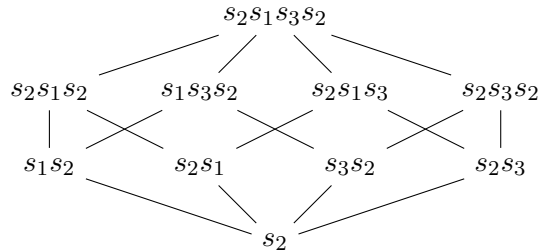
\begin{figure}[htbp]
\begin{tikzpicture}[scale=0.5]
  \node(max) at (0,4) {$s_2s_1s_3s_2$};
  \node (a1) at (-2,2) {$s_1s_3s_2$};
  \node (a2) at (-6,2) {$s_2s_1s_2$};
  \node (a3) at (6,2)  {$s_2s_3s_2$};
  \node (a4) at (2,2)  {$s_2s_1s_3$};
  \node (b1) at (-6,0) {$s_1s_2$};
  \node (b2) at (2,0) {$s_3s_2$};
  \node (b3) at (-2,0)  {$s_2s_1$};
  \node (b5) at (6,0)  {$s_2s_3$};
  \node (c1) at (-0,-2) {$s_2$};
  \draw(a3)--(b2)--(a1);\draw (a4)--(b3)--(c1)--(b2);
  \draw (b5) -- (a3) -- (max) --(a2) --(b3);
  \draw (c1) -- (b1) -- (a1) -- (max) -- (a4) -- (b5);
  \draw (b1)--(a2);
  \draw (c1) -- (b5);
\end{tikzpicture}
\caption{Bruhat interval of $[s_2,s_2s_1s_3s_2]$.}\label{fig:Bruhat_interval_uv}
\end{figure}
We ask is under what conditions does the polynomial $P_{u,v}(q)$ factor nicely?  and if so, does the interval $[u,v]$ also decompose as a poset?  Is there a generalization of the characterization Theorem \ref{thm:BP_characterization} for arbitrary intervals $[u,v]$?  We remark that the poset structure of the interval $[u,v]$ has connections to Kazhdan-Lustig theory and Richardson varieties \cite{KL79,KL80}.
\end{question}


\begin{question}
Let $w=vu$ be a parabolic decomposition with respect to $J$.  In Remark~\ref{rem:BPnice} we describe the coset intervals $[e,w]\cap v_0W^J$ for $v_0 \in [e,v]^J$ and show they are poset isomorphic to $[e,u_0]$ for some $u\leq u_0\leq u'$ where $u'$ is the maximal element of $[e,w]\cap W_J$.  

We ask if there is nice description of the set of all $u_0\in[u,u']$ that appear for some $v_0\in[e,v]^J$.  If $u=u'$, then $w=vu$ is a BP decomposition by Theorem \ref{thm:BP_characterization}.  Note that not every element of $[u,u']$ may appear in this set.
\end{question}

\begin{question}
    In Sections \ref{s:iterated_BPs} and \ref{S:BP_pattern_avoidance}, we see several cases of BP-decompositions and staircase diagrams used to enumerate classes of (rationally) smooth elements.  To what extend can these structures help with enumerating other classes of Coxeter group elements.  For example, can we use BP decompositions to calculate the generating series for the number of permutations that avoid $3412$, $52341$, and $635241$ from Theorem \ref{thm:pattern_complete_BP}?
\end{question}

\bibliographystyle{plain}
\bibliography{palindromic}

\end{document}